\documentclass[11pt,preprint]{article}

\usepackage{amsmath, latexsym, amsfonts, amssymb, amsthm}
\usepackage{graphicx, color,hyperref,dsfont,epsfig,caption,wrapfig,subfig}
\usepackage{pstricks}
\usepackage{movie15}
\hypersetup{
    linktoc=page,
    linkcolor=red,          % color of internal links
    citecolor=blue,        % color of links to bibliography
    filecolor=blue,      % color of file links
    urlcolor=cyan,
    colorlinks=true           % color of external links
}

\numberwithin{equation}{section}
\newcounter{exa}
\newtheorem{theorem}{Theorem}[section]
\newtheorem{lemma}[theorem]{Lemma}
\newtheorem{proposition}[theorem]{Proposition}
\newtheorem{corollary}[theorem]{Corollary}
\newtheorem{rem}[theorem]{Remark}
\newtheorem{remark}[theorem]{Remark}
\newtheorem{definition}[theorem]{Definition}

\newtheorem{example}[exa]{Example}

\theoremstyle{definition}

\renewcommand{\tilde}{\widetilde}          % wider `tilde'
\DeclareMathSymbol{\leqslant}{\mathalpha}{AMSa}{"36} % nicer `smaller or equal'
\DeclareMathSymbol{\geqslant}{\mathalpha}{AMSa}{"3E} % nicer `larger or equal'
\DeclareMathSymbol{\eset}{\mathalpha}{AMSb}{"3F}     % nicer `emptyset'
\renewcommand{\leq}{\;\leqslant\;}                   % redef. of < or =
\renewcommand{\geq}{\;\geqslant\;}                   % redef. of > or =
\newcommand{\R}{\mathbb{R}}
\newcommand{\Z}{\mathbb{Z}}
\newcommand{\N}{\mathbb{N}}
\newcommand{\Q}{\mathbb{Q}}
\newcommand{\E}{\mathds{E}}
\newcommand{\Pb}{\mathds{P}}
\newcommand{\ind}{\mathds{1}}

\def\S{\mathbb{S}}
\def\T{\mathbb{T}}

\def\bi{\begin{itemize}}
\def\ei{\end{itemize}}
\def\ni{\noindent}
\def\bnum{\begin{enumerate}}
\def\enum{\end{enumerate}}
\def\LB{\mathcal{B}} %Liouville BM
\def\<#1{\langle #1 \rangle}

\def\r{\mathbf{r}}
\def\t{\mathbf{t}}

\def\p{\mathbf{p}}

\title{Liouville Brownian motion at criticality}
\author{  R\'emi Rhodes \footnote{Universit{\'e} Paris-Dauphine, Ceremade, F-75016 Paris, France. Partially supported by grant ANR-11-JCJC  CHAMU} \footnotetext[2]{Partially supported by grant ANR-11-JCJC  CHAMU}  \and 
 Vincent Vargas \footnote{Ecole Normale Sup\'erieure, DMA, 45 rue d'Ulm,  75005 Paris, France. Partially supported by grant ANR-11-JCJC  CHAMU} }

\date{\vspace{-5ex}}

\begin{document}

\maketitle
 
%\vspace{1cm}
 
\begin{abstract}
In this paper, we construct the  Brownian motion of  Liouville Quantum Gravity with central charge $c=1$ (more precisely we restrict to the corresponding free field theory). Liouville quantum gravity with $c=1$ corresponds  to two-dimensional string theory and is the conjectural scaling limit of  large planar maps weighted with a $O(n=2)$ loop model or a $Q=4$-state Potts model  embedded in a two dimensional surface in a conformal manner.

Following \cite{GRV1}, we start by constructing the critical LBM from one fixed point $x\in\R^2$ (or $x\in\S^2$), which amounts to changing the speed of a standard planar Brownian motion depending on the local behaviour of the critical Liouville measure $M'(dx)=-X(x)e^{2X(x)}\,dx$ (where $X$ is a Gaussian Free Field, say on $\S^2$). Extending this construction simultaneously to all points in $\R^2$ requires a fine analysis of the potential properties of the measure $M'$. This allows us to construct a strong Markov process with continuous sample paths  living on the support of $M'$, namely a dense set of Hausdorff dimension $0$. We finally construct the associated Liouville semigroup, resolvent, Green function, heat kernel and Dirichlet form.

In passing, we extend to quite a general setting the construction of the critical Gaussian multiplicative chaos that was initiated in \cite{Rnew7,Rnew12} and also establish new capacity estimates for the critical Gaussian multiplicative chaos. 
\end{abstract}
%\vspace{0.5cm}
\footnotesize

%\noindent{\bf Short Title.} 

%\vspace{0.5cm}

\noindent{\bf Key words or phrases:} Gaussian multiplicative chaos, critical Liouville quantum gravity, Brownian motion, heat kernel, potential theory.

%\noindent{\bf MSC 2000 subject classifications: 60J60, 60G15, 60J35}
\normalsize

%\newpage
%%\vspace{0.1cm}
%%
\tableofcontents

%%%%%%%%%%%%%%%%%%%%%%%%%%%%%%%%%%%%%%%%%%%%%%%%%%
%%%%%%%%%%%%%%%%%%%%%%%%%%%%%%%%%%%%%%%%%%%%%%%%%%

%%%%%%%%%%%%%%%%%%%%%%%%%%%%%%%%%%%%%%%%%%%%%%%%%%
\section{Introduction}
%%%%%%%%%%%%%%%%%%%%%%%%%%%%%%%%%%%%%%%%%%%%%%%%%%
\subsection{Physics motivations}
%%%%%%%%%%%%%%%%%%
Liouville Quantum Field Theory is a two dimensional conformal field theory which plays an important part in two dimensional models of Euclidean quantum gravity. Euclidean  quantum gravity is an attempt to quantize general relativity based on Feynman's functional integral and on the Einstein-Hilbert action principle.  More precisely, one couples a Conformal Field Theory  (CFT)  with central charge $c$  to gravity. A famous example is the coupling of $c$ free scalar matter fields to gravity,  leading to an interpretation of such a specific theory of $2d$-Liouville Quantum Gravity   as a bosonic string  theory  in $c$ dimensions \cite{Pol}. 
 
It is shown in \cite{Pol,cf:KPZ,cf:Da}) that the coupling of the CFT with gravity  can be factorized as a tensor product: the random metric governing the geometry of the space that the CFT lives on   is independent of the CFT and  roughly takes on the form \cite{Pol,cf:KPZ,cf:Da} (we consider an Euclidean background metric for simplicity):
\begin{equation}\label{tensor}
g(x)=e^{\gamma X(x)}dx^2,
\end{equation}
where the fluctuations of the field $X$ are governed by the Liouville action and the parameter $\gamma$ is related to the central charge of the CFT via the famous result in \cite{cf:KPZ} (for $c\leq 1$)
\begin{equation}\label{KPZgammac}
\gamma=\frac{1}{\sqrt{6}}(\sqrt{25-c}-\sqrt{1-c}).
\end{equation}
Therefore $\gamma\in[0,2]$. The reader is referred to \cite{DKRV} for further details and a rigorous construction of the Liouville field and to  \cite{david-hd,cf:Da,cf:DuSh,DFGZ,DistKa,bourbaki,GM,glimm,Kle,cf:KPZ,Nak,Pol} for further insights on $2d$-Liouville quantum gravity. For these reasons and though it may be an interesting field in its own right, Liouville quantum field theory (governing the metric in $2d$-Liouville quantum gravity)  is an important object in theoretical physics.

Here we will restrict to the case when the   the cosmological constant in the Liouville action is set to $0$, turning the field $X$ in \eqref{tensor} into a Free Field, with appropriate boundary conditions.   In the subcritical case $\gamma<2,c<1$, the geometry of the metric tensor \eqref{tensor} (with $X$ a free field) is mathematically investigated in \cite{cf:DuSh,GRV1,GRV2,review,spectral}   and the famous KPZ scaling relations \cite{cf:KPZ} are rigorously proved in \cite{Rnew4,cf:DuSh,Rnew10} in a geometrical framework (see also \cite{David-KPZ} for a non rigorous heat kernel derivation). This paper  focuses on the coupling of a CFT with a $c=1$ central charge to gravity or equivalently $2d$ string theory: in this case, one has $\gamma=2$ by relation \eqref{KPZgammac}.  For an excellent review on $2d$ string theory, we refer to Klebanov's lecture notes \cite{Kle}. As expressed by Klebanov in \cite{Kle}: "Two-dimensional string theory is the kind of toy model which possesses a remarkably simple structure but at the same time incorporates 
some of the physics of string theories embedded in higher dimensions".  Among the $\gamma\leq 2$ theories, the case $\gamma=2$ probably possesses the richest structure, inherited from its specific status of phase transition. 
For instance, the construction of the volume form, denoted by $M'$, associated to the metric tensor \eqref{tensor} is investigated in \cite{Rnew7,Rnew12} where it is proved that it takes on the unusual form: 
\begin{equation}
M'(dx)=-X(x)e^{2X(x)}dx,
\end{equation}
which also coincides with a proper renormalization of $e^{2X(x)}dx$. The reader may also consult  \cite{complex} for a construction of non trivial conformally invariant gravitationally dressed vertex operators, the so-called {\bf tachyonic fields}, and \cite{Kle} for more physics insights. In this paper, we will complete this picture by constructing the Brownian motion (called critical Liouville Brownian motion, critical LBM for short), semi-group, resolvent, Dirichlet form, Green function and heat kernel of the metric tensor \eqref{tensor} with $\gamma=2$.  

We further point out that Liouville quantum gravity is conjecturally related to randomly triangulated random surfaces (see \cite{DKRV} for precise conjectures) weighted by discrete critical statistical physics models.  For  $c=1$, these models  include one-dimensional matrix models (also called ``matrix quantum mechanics'' (MQM)) \cite{BKZ,GM,GZ,GrossKleban,GubserKleban,KKK,Kleb2,Parisi,Polch,sugino}, the  so-called $O(n)$ loop model on a random planar lattice for $n=2$ \cite{kostov91,kostov92,Kostov:houches,KS}, and the $Q$-state Potts model on a random lattice for $Q=4$ \cite{BBM,Daul,Ey}.  The  critical LBM  is therefore conjectured to be the scaling limit of random walks on large planar maps weighted with a $O(n=2)$ loop model or a $Q=4$-state Potts model, which are embedded in a two dimensional surface in a conformal manner as explained in \cite{DKRV}. 
For an introduction to the above mentioned $2d$-statistical models, see, e.g., \cite{Nienhuis}. We further mention \cite{curien,eynard} for recent advances on this topic in the context of pure gravity, i.e. with no coupling with a CFT.

To complete this overview, we point out that the notions of diffusions or heat kernel are at the core of  physics literature about Liouville quantum gravity (see \cite{amb,amb2,calg1,calg2,david-hd,David-KPZ,wata} for instance): for more on this, see the subsection on the associated distance.

%A heat kernel derivation of the KPZ formula is  obtained in \cite{David-KPZ}. Fractal structure of quantum space-time is also investigated in %\cite{amb,amb2,david-hd,wata} via diffusions and heat kernel properties, obtaining relations about the fractal dimensions of quantum space-%time.  

\subsection{Strategy and results} 
%%%%%%%%%%%%% 
Basically, our approach of the metric tensor $e^{2X(x)}\,dx^2$ (where $X$ is a Free Field) relies on the construction of the associated Brownian motion $\LB$, called the critical LBM. Standard results of $2d$-Riemannian geometry tell us that the law of this Brownian motion is a time change of a standard planar Brownian motion $B$ (starting from $0$):
\begin{equation}
\LB_t^x=x+B_{\langle \LB^x\rangle_t}
\end{equation}
where the quadratic variations $\langle \LB^x\rangle$ are formally given by:
\begin{equation}\label{naive}
 \langle \LB^x\rangle_t=F'(x,t)^{-1},\quad F'(x,t)=\int_0^te^{2X(x+B_r)}\,dr.
 \end{equation}
Put in other words, we should integrate the weight $e^{2X(x)}$ along the paths of the Brownian motion $x+B$ to construct a mapping $t\mapsto F'(x,t)$. The inverse of this mapping corresponds to the quadratic variations of $\LB^x$. Of course, because of the irregularity of the field $X$, giving sense to \eqref{naive} is not straightforward and one has to apply a renormalization procedure: one has to apply a cutoff to the field $X$ (a procedure that smoothes up the field $X$) and pass to the limit as the cutoff is removed. The procedure is rather standard in this context. Roughly speaking, one introduces an approximating field $X_\epsilon$ where the parameter $\epsilon$ stands for the extent to which one has regularized the field $X$ (we have $X_\epsilon \rightarrow X$ as $\epsilon$ goes to $0$).  One then defines
\begin{equation}\label{renormF}
F^\epsilon(x,t)=\int_0^te^{2X_\epsilon(x+B_r)}\,dr
\end{equation} and one looks for a suitable deterministic renormalization $a(\epsilon)$ such that the family $a(\epsilon)F^\epsilon(x,t)$ converges towards a non trivial object as $\epsilon\to 0$.  In the subcritical case $\gamma<2$, the situation is rather well understood as the family $a(\epsilon)$ roughly corresponds to 
$$a(\epsilon)\simeq\exp(-\frac{\gamma^2}{2}\E[X_\epsilon(x)^2])$$ (the dependence of the point $x$ is usually fictive and may easily get rid of) in such a way that
\begin{equation}\label{renorm}
a(\epsilon)\int_0^te^{\gamma X_\epsilon(x+B_r)}\,dr
\end{equation}
 converges towards a non trivial limit. We will call this renormalization procedure standard: it has been successfully applied to construct random measures of the form $e^{\gamma X(x)}dx$  \cite{cf:Kah,cf:RoVa,cf:RoVa1}  (the reader may consult \cite{review} for an overview  on Gaussian multiplicative chaos theory). Though the choice of the cutoff does not affect the nature of the limiting object, a proper choice of the cutoff turns the expression \eqref{renorm} into a martingale. This is convenient to handle the convergence of this object.

At criticality ($\gamma=2$), the situation is conceptually more involved.  It is known \cite{Rnew7} that the standard 
renormalization procedure of the volume form yields a trivial object. Logarithmic corrections in the choice of the family $a(\epsilon)$ are necessary  (see \cite{Rnew7,Rnew12}) and the limiting measure that we get is the same as that corresponding to a metric tensor of the form $-X(x)e^{2X(x)}$. Observe that it is not straightforward to see at first sight that such   metric tensors coincide, or even are positive. This subtlety is at the origin of some misunderstandings in the physics literature where the two forms of the tachyon field $e^{2X(x)}$ and $-X(x)e^{2X(x)}$ appear without making perfectly clear that they coincide. 

The case of the LBM at criticality obeys this rule too. We will prove that non standard logarithmic corrections are necessary to make the change of time $F^\epsilon$ converge and they produce the same limiting change of times as that corresponding to a metric tensor of the form $-X(x)e^{2X(x)}$. This summarizes the almost sure convergence of $F^\epsilon$ for one given fixed point $x\in\R^2$. Yet, if one wishes to define a proper Markov process, one has to go one step further and establish that, almost surely,  $F^\epsilon(x,\cdot)$ converges simultaneously for all possible starting points $x\in\R^2$: the place of the "almost sure" is important and gives rise to difficulties that are conceptually far different from the construction from one given fixed point. In \cite{GRV1}, it is noticed that this simultaneous convergence is possible as soon as the volume form $M_\gamma(dx)$ associated to the metric tensor  \eqref{tensor} is regular enough so as to make the mapping
$$x\mapsto \int_{\R^2}\ln_+\frac{1}{|x-y|}M_\gamma(dx)$$ continuous. When $\gamma<2$, multifractal analysis shows that the measure $M_\gamma$ possesses a power law decay of the size of balls and this is enough to ensure the continuity of this mapping. In the critical case $\gamma=2$, the situation is more complicated because the measure $M'$ is rather wild: for instance the Hausdorff dimension of its support is zero \cite{KNW}.
Furthermore, the decay of the size of balls investigated in \cite{KNW} shows that continuity and even finiteness of the mapping
\begin{equation}\label{intro:lnF}
x\mapsto \int_{\R^2}\ln_+\frac{1}{|x-y|}M'(dx)
\end{equation} is unlikely to hold on the whole of $\R^2$. Yet we will show that we can have a rather satisfactory control of the size of balls for all $x$ belonging to a set of full $M'$-measure, call it $S$:
$$\forall x\in S,\quad \sup_{r\in]0,1]}M'(B(x,r))(-\ln r)^p<+\infty$$ for some $p$ large enough.
In particular this estimate shows that the expression \eqref{intro:lnF} is finite for every point $x\in S$. Also, this estimate answers a question raised in \cite{KNW} (a similar estimate was proved by the same authors \cite{BKNSW} in the related case of the discrete multiplicative cascades). We will thus construct the change of time $F'(x,\cdot)$  simultaneously for all $x\in S$. What happens on  the complement of $S$ does not matter that much since it is a set with null $M'$ measure. Yet we will extend the change of time $F'(x,\cdot)$ to the whole of $\R^2$.

Once $F'$ is constructed, potential theory \cite{fuku} tells us that the LBM at criticality is a strong Markov process, which preserves the critical measure $M'$. We will then define the semigroup, resolvent, Green function and heat kernel associated to the LBM at criticality.  

We stress that, once all the pieces of the puzzle are glued together, this LBM at criticality appears as a rather weird mathematical object. It is a Markov process with continuous sample paths living on a very thin set, which is dense in $\R^2$ and has Hausdorff dimension $0$. And in spite of this rather wild structure, the LBM at criticality is regular enough to possess  a (weak form of) heat kernel. Beyond the possible applications in physics, we feel that the study of such an object is a fundamental and challenging mathematical problem, which is far from being settled in this paper.

\subsection{Discussion about the associated distance} 
%%%%%%%%%%%%%%%%%%%%%%%%%%%%%%  
An important question related to this work is the existence of a distance associated to the metric tensor \eqref{tensor} for $\gamma\in [0,2]$. The physics literature contains several suggestions to handle this problem, part of which are discussed along the forthcoming speculative lines. Also since we will base part of our discussion on the results established in \cite{spectral}, we will focus on the non critical case $\gamma<2$ though a similar discussion should hold in the critical case. In this context, we denote $F(x,t)=\int_0^te^{ \gamma X(x+B_r)}\,dr$ the associated additive functional along the Brownian paths.  We know that for all $t>0$ there exists a Liouville heat kernel $\p^X_t(x,y)$. Many papers in the physics literature have argued that the Liouville heat kernel should have the following representation which is classical in the context of Riemannian geometry
\begin{equation}\label{representationheat}
\p^X_t(x,y) \sim \frac{e^{-d(x,y)^2/t}}{t}
\end{equation}   
where $d(x,y)$ is the associated distance, i.e. the "Riemannian distance" defined by \eqref{tensor}. From the representation \eqref{representationheat}, physicists \cite{amb,amb2,wata} have derived many fractal and geometrical properties of Liouville quantum gravity. In particular, the paper \cite{wata} established an intriguing formula for the dimension $d_H$ of Liouville quantum gravity which can be defined by the heuristic
\begin{equation*}    
M_\gamma(\lbrace y; \: d(x,y) \leq r \rbrace) \sim r^{d_H}
\end{equation*}
where $M_\gamma$ is the associated volume form. Note that the meaning of the above definition is not obvious since $M_\gamma$ is a multifractal random measure. 

Along the same lines, a recent physics paper \cite{David-KPZ}  establishes an interesting heat kernel derivation of the KPZ formula. The idea behind the paper is that, if relation \eqref{representationheat} holds, then one can extract the metric from the heat kernel by using the Mellin-Barnes transform given by
\begin{equation*}  
\int_0^{\infty} t^{s-1}\p^X_t(x,y) dt.
 \end{equation*}
Indeed, a standard computation gives the following equivalent for $s \in ]0,1[$
\begin{equation}\label{heuristicDavid}
\int_0^\infty  t^{s-1} \frac{e^{-d(x,y)^2/t}}{t} \: dt \underset{d(x,y) \to 0}{\sim}  \frac{C_s}{d(x,y)^{2(1-s)}}
\end{equation}
where $C_s$ is some positive constant. Though it is not clear at all that such a relation holds rigorously, it gives at least a way of defining a notion of capacity dimension for which a KPZ relation has been heuristically derived in \cite{David-KPZ}. Thanks to the relation \eqref{heuristicDavid}, the authors claim that this yields a geometrical version of the KPZ equation which does not rely on the Euclidean metric. Recall that the rigorous geometrical derivations of KPZ in \cite{cf:DuSh,Rnew10} rely on the measure $M_\gamma$ and imply working with Euclidean balls\footnote{At the time of publishing this manuscript, this  heat kernel based KPZ formula has been rigorously proved in \cite{heatKPZ}.}. 

\begin{rem}
In fact, after the present work and based on results in the field of fractal diffusions, it is argued in \cite{MRVZ} that one should replace the heuristic \eqref{representationheat} by the following heuristic
\begin{equation*}
\p^X_t(x,y) \sim \frac{e^{-\left ( \frac{d(x,y)^{d_H}}{t} \right )^{\frac{1}{d_H-1}}}    }{t}
\end{equation*}
where $d_H$ is the dimension of Liouville quantum gravity. Hence, we should also rather get the relation
\begin{equation*}
\int_0^{\infty} t^{s-1}\p^X_t(x,y) dt \underset{d(x,y) \to 0}{\sim}  \frac{C_s}{d(x,y)^{d_H (1-s)}}.
\end{equation*}
  
\end{rem}

 %%%%%%%%%%%%%%%%%%%%%%%%%%%%%%%%%%%%%%%%%%%%%%%%%%%%%%%%%%%%%%%%%%%%%%%%%%%%%%%%%%%%%%%%%%%
\subsubsection*{Acknowledgements}
%%%%%%%%%%%%%%%%%%%%%%%%%%%%%%%%%%%%%%%%%%%%%%%%%%%%%%%%%%%%%%%%%%%%%%%%%%%%%%%%%%%%%%%%%%%
The authors  wish to thank Antti Kupiainen, Miika Nikula and Christian Webb for many very  interesting discussions which have helped a lot in understanding  the specificity of the critical case and Fran\c{c}ois David who always takes the time to answer their questions with patience and kindness.

%%%%%%%%%%%%%%%%%%%%%%%%%%%%%%%%%%%%%%%%%%%%%%%%%%%%%%%%%%%%%%%%%%%%%%%%%%%%%%%%%%%%%%%%%%%%%%%%%%%%%%%
\section{Setup}\label{sec.setup}
%%%%%%%%%%%%%%%%%%%%%%%%%%%%%%%%%%%%%%%%%%%%%%%%%%
%%%%%%%%%%%%%%%%%%%%%%%%%%%%%%%%%%%%%%%%%%%%%%%%%%
 
In this section, we draw up the framework to construct the Liouville Brownian motion at criticality on the whole plane $\R^2$. Other geometries are possible and discussed at the end of the paper.
 
\subsection{Notations}\label{index}
%%%%%%%%%%%%%%%%%%%%

In what follows, we will consider Brownian motions $B$ or $\bar{B}$ on $\R^2$ (or other geometries) independent of the underlying Free Field. We will denote by $\E^Y$ or $\Pb^Y$ expectations and probability with respect to  a field $Y$. For instance, $\E^X$ or $\Pb^X$ (resp. $\E^B$ or $\Pb^B$) stand for expectation and probability with respect to the log-correlated field $X$ (resp. the Brownian motion $B$). 
For $d\geq 1$, we consider the space $C(\R_+,\R^d)$ of continuous functions from $\R_+$ into $\R^d$ equipped with the topology of uniform convergence over compact subsets of $\R_+$.

\subsection{Representation of a log-correlated field}
%%%%%%%%%%%%%%%
In this section we introduce the log-correlated Gaussian fields $X$ on $\R^2$ that we will work with throughout this paper.   One may consider other geometries as well, like the sphere $\S^2$ or the torus $\T^2$ (in which case an adaptation of the setup and proofs is straightforward).    We will represent them via a white noise decomposition.

 We consider a family of  centered  Gaussian processes $((X_{\epsilon}(x))_{x \in \R^2})_{\epsilon> 0}$ with covariance structure given, for $1\geq \epsilon>\epsilon'> 0$,  by:
\begin{equation} \label{corrX}
K_\epsilon(x,y)=\E[X_{\epsilon}(x)X_{\epsilon'}(y)]= \int_1^{\frac{1}{\epsilon }}\frac{k(u,x,y)}{u}\,du
\end{equation}
for some  family $(k(u,\cdot,\cdot))_{u\geq 1}$ of covariance kernels satisfying:
\begin{description}
\item[A.1] $k$ is nonnegative, continuous. 
\item[A.2] $k$ is locally Lipschitz on the diagonal, i.e. $\forall  R>0$, $\exists C_R>0$, $\forall |x|\leq R$, $\forall u\geq 1$, $\forall y\in\R^2$ 
%$$ |y-x|\leq R/u\Rightarrow  |k(u,x,x)-k(u,x,y)|\leq C_Ru|x-y|$$
$$  |k(u,x,x)-k(u,x,y)|\leq C_Ru|x-y|$$
\item[A.3] $k$ satisfies the integrability condition: for each compact set $S$, $$\sup_{x\in S,y\in\R^2}\int_{\frac{1}{|x-y|}}^{\infty}\frac{k(u,x,y)}{u}\,du<+\infty.$$ 
\item[A.4] the mapping $H_\epsilon(x)=\int_1^{ \frac{1}{\epsilon}}  \frac{k(u,x,x)}{u}\,du-\ln\frac{1}{\epsilon}$ converges pointwise as $\epsilon\to 0$ and for all compact set $K$, 
$$\sup_{x,y\in K}\sup_{\epsilon\in ]0,1]}\frac{|H_\epsilon(x)-H_\epsilon(y)|}{|x-y|}<+\infty.$$
\item[A.5] for each compact set $K$, there exists a constant $C_K>0$ such that
$$k(u,x,y)\geq k(u,x,x)(1-C_Ku^{1/2}|x-y|^{1/2})_+$$ for all $x\in K$ and $y\in\R^2$.
 for all $u\geq 1$, $x\in K$ and $y\in\R^2$.
\end{description}
Such a construction of Gaussian processes is carried out in  \cite{Rnew1,sohier}  in the translation invariant case. Furthermore, [A.2] implies the following relation that we will use throughout the paper: for each compact set $S$, there exists a constant $c_S>0$ (only depending on $k$) such that for all $y\in S$,  $\epsilon\in (0,1]$ and $w\in B(y,\epsilon)$, we have 
\begin{equation}\label{property:k}
\ln\frac{1}{\epsilon}-c_S\leq K_\epsilon(y,w)\leq \ln\frac{1}{\epsilon}+c_S.
\end{equation}

We denote by $\mathcal{F}_\epsilon$ the sigma algebra generated by $\{X_u(x);\epsilon\leq u,x\in\R^2\}$.  
\subsection{Examples}
%%%%%%%%%%%%%%%%%%%%%%%%%%%%%%%%%%
We explain first a Fourier white noise decomposition of log-correlated translation invariant fields as this description appears rather naturally in physics. Consider a nonnegative even function $\varphi$  defined on $\R^2$ such that $\lim_{|u|\to \infty}|u|^2\varphi(u)=1$. We consider the kernel  
\begin{equation}\label{Fourierrepresentation}
K(x,y)=\int_{\R^2}e^{i\langle u,x-y\rangle}\varphi(u)\,du.
\end{equation}
We consider the following cut-off approximations
\begin{equation}\label{FWN}
K_\epsilon(x,y)=\frac{1}{2\pi}\int_{B(0,\epsilon^{-1})}e^{i\langle u,x-y\rangle}\varphi(u)\,du.
\end{equation}
The kernel $K$ can be seen as the prototype of kernels of log-type in dimension $2$. It has obvious counterparts in any dimension. The cut-off approximation is quite natural, rather usual in physics (sometimes called the ultraviolet cut-off) and has well known analogues on compact manifolds (in terms of series expansion along eigenvalues of the Laplacian for instance). If $X$ has covariance given by \eqref{Fourierrepresentation}, then $X$ has the following representation
\begin{equation*}
X(x)= \int_{\R^2}  e^{i \langle u . x \rangle } \sqrt{\varphi(u)} (W_1(du)+iW_2(du))
\end{equation*} 
where $W_1(du)$ and $W_2(du)$ are independent Gaussian distributions. The distributions $W_1(du)$ and $W_2(du)$ are functions of the field $X$ (since they are the real and imaginary parts of the Fourier transform of $X$). The law of $W_1(du)$ is $W(d u)+ W(-d u)$ where $W$ is a standard white noise and the law of $W_2(du)$ is $\tilde{W}(d u)- \tilde{W}(-d u)$ where $\tilde{W}$ is also standard white noise. One can then consider the following family with covariance \eqref{FWN} and which fits into our framework as it corresponds to adding independent fields
\begin{equation*}
X_\epsilon(x)= \int_{B(0,\epsilon^{-1}) }  e^{i \langle u . x \rangle} \sqrt{\varphi(u)} (W_1(du)+iW_2(du))
\end{equation*}  
Notice also that the approximations $X_\epsilon$ are functions of the original field $X$ since $W_1(du)$ and $W_2(du)$ are functions of the field $X$.

Notice that $K_\epsilon$ can be rewritten as ($S$ stands for the unit sphere and $ds$ for the uniform probability measure on $S$)
$$K_\epsilon(x,y)=\int_0^{\frac{1}{\epsilon}}\frac{k(r,x,y)}{r}\,dr,\quad \text{ with }k(r,x,y)=r^2\int_S\varphi (r s)\cos(r\langle x-y,s\rangle)\,ds.$$
Let us simplify a bit the discussion by assuming that $\varphi$ is isotropic. In that case, it is plain to check that assumptions [A.1-5] are satisfied (in the slightly extended context of integration over $[0, \frac{1}{\epsilon}]$ instead of $[1, \frac{1}{\epsilon}]$ but this is harmless as $K_{1}(x,y)$ is here a very regular Gaussian kernel).

\begin{example}{\bf Massive Free Field (MFF).}  \label{ex:MFF}
The whole plane MFF is a centered Gaussian distribution with covariance kernel given by the Green function of the operator $2\pi (m^2-\triangle)^{-1}$ on $\R^2$, i.e. by:
\begin{equation}\label{MFFkernel}
\forall x,y \in \R^2,\quad G_m(x,y)=\int_0^{\infty}e^{-\frac{m^2}{2}u-\frac{|x-y|^2}{2u}}\frac{du}{2 u}.
\end{equation}
  The real $m>0$ is called the mass. This kernel is of $\sigma$-positive type in the sense of Kahane \cite{cf:Kah} since we integrate a continuous function of positive type with respect to a positive measure. It is furthermore  a star-scale invariant kernel (see \cite{Rnew1,sohier}): it can be rewritten as 
\begin{equation}\label{MFF3}
G_m(x,y)=\int_{1}^{+\infty}\frac{k_m(u(x-y))}{u}\,du.
\end{equation}
 for some continuous covariance kernel $k_m(z)=\frac{1}{2}\int_0^\infty e^{-\frac{m^2}{2v}|z|^2-\frac{v}{2}}\,dv$ and therefore satisfies the assumptions [A.1-5].

One may also choose the Fourier white noise \eqref{FWN} decomposition with $\varphi(u)=\frac{1}{|u|^2+m^2}$ or the semigroup covariance structure 
$$\E[X_{\epsilon}(x)X_{\epsilon'}(y)]= \pi \int^{\infty}_{\max(\epsilon ,\epsilon')^2} p(u,x,y)e^{-\frac{m^2}{2}u}\,du,$$
which also satisfies assumptions [A.1-5] (modulo a change of variable $u=1/v^2$ in the above expression: see \cite[section D]{Rnew12}).
\end{example}

\begin{example}{\bf Gaussian Free Field (GFF).}  \label{ex:GFF}
Consider a bounded open domain $D$ of $\R^2$. Formally, a GFF on $D$ is a Gaussian distribution with covariance kernel given by the Green function of the Laplacian on $D$ with prescribed boundary conditions. We describe here the case of Dirichlet boundary conditions. The Green function is then given by the formula:
\begin{equation}\label{bounGreen}
G_D (x,y)=  \pi \int_{0}^{\infty}p_D(t,x,y)dt
 \end{equation}
where $p_D$  is the (sub-Markovian) semi-group of a Brownian motion $B$ killed upon touching the boundary of $D$, namely for a Borel set $A\subset D$
$$\int_A p_D(t,x,y)\,dy=P^x(B_{t} \in A, \; T_D > t) $$ with $T_D=\text{inf} \{t \geq 0, \; B_t\not \in D \}$. The most  direct way to construct a cut-off family of the GFF on $D$ is then to consider  a white noise $W$ distributed on $D\times \R_+$ and to define:
 $$X(x)=\sqrt{\pi}\int_{D\times \R_+}p_D(\frac{s}{2},x,y)W(dy,ds).$$
One can check that $\E[X(x)X(x')]=\pi\int_0^{\infty} p_D(s,x,x')\,ds=G_D(x,x') $. The corresponding cut-off approximations are given by:
$$X_\epsilon(x)=\sqrt{\pi}\int_{D\times [\epsilon^2,\infty[}p_D(\frac{s}{2},x,y)W(dy,ds).$$
They have the following covariance structure  
\begin{equation}\label{eq:GFFcutoff}
\E[X_{\epsilon}(x)X_{\epsilon'}(y)]= \pi \int^{\infty}_{\max(\epsilon ,\epsilon')^2} p_D(u,x,y) \,du,
\end{equation}
which also satisfies assumptions [A.1-5] (on every subdomain of $D$ and modulo a change of variable $u=1/v^2$ in the above expression: see also \cite[section D]{Rnew12}).
\end{example}

For some technical reasons, we will sometimes also consider either of the following assumptions:
\begin{description}
\item[A.6]  $k(v,x,y)=0$ for $|x-y|\geq Dv^{-1}(1+2\ln v)^\alpha$ for some constants $D,\alpha>0$, 
\item[A.6']  $(k(v,x,y))_v$ is the family of kernels presented in examples \ref{ex:MFF} or \ref{ex:GFF}.
\end{description}

\subsection{Regularized Riemannian geometry}
%%%%%%%%%%%%%%%%%%%%%%%%%%%%%%%%%%
We would like to consider  a Riemannian metric tensor on $\R^2$ (using conventional notations in Riemannian geometry) of the type $ e^{ 2 X_\epsilon(x)}dx^2$, where $dx^2$ stands for the standard Euclidean metric on $\R^2$. Yet, as we will see, such an object has no suitable limit as $\epsilon$ goes to $0$. So, for future renormalization purposes, we rather consider:
$$g_\epsilon(x)dx^2=\sqrt{-\ln \epsilon}\,\epsilon^2 \, e^{ 2 X_\epsilon(x)}dx^2.$$

\subsubsection{Volume form}
%%%%%%%%%%%%%%%%
The Riemannian volume on the manifold $(\R^2,g_\epsilon)$ is given by:
\begin{equation}\label{defchaos}
M_{ \epsilon}(dx)=\sqrt{-\ln \epsilon}\,\epsilon^2 \, e^{ 2 X_\epsilon(x)}\,dx,
\end{equation}
where $dx$ stands for the Lebesgue measure on $\R^2$, and will be called $\epsilon$-regularized critical  measure. The study of the limit of the random measures $(M_{ \epsilon}(dx))_\epsilon$ is carried out in \cite{Rnew7,Rnew12} in a less general context.  It is based on the study of  the limit of  the family $M_{\epsilon}'(dx)$ defined by:
\begin{align*}
M_{\epsilon}'(dx)& := (2\,\ln\frac{1}{\epsilon}-X_{\epsilon}(x))e^{2X_{\epsilon}(x)-2\ln\frac{1}{\epsilon}}dx .
\end{align*}
  
We will extend the results in  \cite{Rnew7} and prove
\begin{theorem}\label{mainderiv}
%$\bullet$ For each bounded open set $A\subset \R^2$, the martingale $(M'_\epsilon(A))_{\epsilon>0}$ converges almost surely towards a  positive random variable denoted by $M'(A)$, such that $M'(A)>0$ almost surely.\\
Almost surely, the (locally signed) random measures $(M'_\epsilon(dx))_{\epsilon> 0}$ converge  as $\epsilon\to 0$ towards a positive random measure $M'(dx)$ in the sense of weak convergence of measures. This limiting measure has full support and is atomless. 
\end{theorem}

%Regarding the Seneta-Heyde norming, we need to introduce some further technical constraints. We follow \cite{Rnew12} and consider the further assumptions:
%\begin{description}
%\item[A.6] for each $u\geq 1$, $(x,y)\mapsto k(u,x,y)$ is of class $C^1$.
%\item[A.7] there exists some constant $C>0$ and $\alpha>1$ such that 
%$$\forall u\geq 1,\forall x,y\in\R^2,\quad |k(u,x,y )|+|\partial_xk(u,x,y)| \leq ue^{-|u(x-y)|^{1/\alpha}}.$$
%\item[A.8] there exists a nonnegative function $g(u,x,y)$, continuous on $\{(u,x,y);u\geq1,x\not = y\}$ such that for each $u\geq 1$, $$\int_{\R^2}g(u,x,z)g(u,z,y)\,dz=k(u,x,y)$$ and $\forall u\geq 1$
%$$\sup_{|x-y|\geq 2u}\big(g(u,x,y)+|\partial_x g(u,x,y)|\big)<+\infty\quad \text{ and }\int_1^{\infty}v^{-1}\int_{}.$$
%\end{description}

Concerning the Seneta-Heyde norming, we have 
\begin{theorem}\label{mainderiv2}
Assume [A.1-5] and either [A.6] or [A.6']. We have the convergence in probability in the sense of weak convergence of measures:  
$$ M_{\epsilon}(dx)\to\sqrt{\frac{2}{\pi}}M'(dx),\quad \text{ as }\epsilon\to 0.$$
\end{theorem}

The proof of Theorem \ref{mainderiv2} is carried out in \cite[section D]{Rnew12} (in fact, it is assumed in \cite{Rnew12} that the family of kernels $(k(v,x,y))_v$ is translation invariant but adapting the proof is straightforward and thus left to the reader).

Beyond its conceptual importance, the Seneta-Heyde norming is crucial to establish, via Kahane's convexity inequalities \cite{cf:Kah},  the study of moments carried out in \cite{Rnew7,Rnew12} and obtain
\begin{proposition}
Assume [A.1-5] and either [A.6] or [A.6']. For each bounded Borel set $A$ and $q\in]-\infty,1[$, the random variable $M'(A)$ possesses moments of order $q$. Furthermore, if $A$ has non trivial Lebesgue measure and $x\in\R^2$:
$$\E[M'(\lambda A+x)^q]\simeq C(q,x) \lambda^{\xi_{M'}(q)}$$ where $\xi_{M'}$ is the power law spectrum of $M'$:
$$\forall q<1,\quad \xi_{M'}(q)=4q-2q^2.$$
\end{proposition}  

Another important result about the modulus of continuity of the measure $M'$ is established in \cite{KNW}. We stress here that we pursue the discussion at a heuristic level since the result in \cite{KNW} is not  general enough to apply in our context (to be precise, it is valid for a well chosen family of kernels $(k(v,x,y))_v$ in dimension $1$ in order to get nice scaling relations for the associated measure $M'$). Anyway, we  expect this result to be true in greater generality and we will not use it in this paper: we are more interested in its conceptual significance. So, by analogy with \cite{KNW},  the measure $M'$ is expected to possess  the following modulus of continuity of "square root of log" type: for all $\gamma<1/2$, there exists a random variable $C$ almost surely finite such that
\begin{equation}\label{modulus}
\forall \text{ ball }B\subset B(0,1),\quad M'(B)\leq C(\ln(1+|B|^{-1}))^{-\gamma}.
\end{equation}
Furthermore, the Hausdorff dimension of the carrier of $M'$ is $0$. By analogy with the results that one gets in the context of multiplicative cascades \cite{BKNSW}, one also expects  that the above theorem \ref{modulus} cannot be improved. In particular, the measure $M'$ does not possess a modulus of continuity   better than a log unlike the subcritical situation explored in \cite{GRV1}, where this property turned out to be crucial for the construction of the LBM as a whole Markov process. This remark is at the origin of the further complications arising in our paper (the critical case) in comparison with \cite{GRV1,GRV2} (the subcritical case).

\subsubsection{Regularized Liouville Brownian motion}
%%%%%%%%%%%%%%%%%%%%%%%%%%%%%%%%%%%%%%%%%%%%%%%%
The main concern of this paper will be the Brownian motion associated with the metric tensor $g_\epsilon$: following standard formulas of Riemannian geometry, one can associate to the Riemannian manifold $(\R^2,g_\epsilon)$ a Brownian motion $\mathcal{B}^\epsilon$:
\begin{definition}[$\epsilon$-regularized critical Liouville Brownian motion, ${\rm LBM}_\epsilon$ for short]\label{d.elbm}
For any fixed $\epsilon>0$, we define the following diffusion on $\R^2$. For any $x\in \R^2$, 
\begin{equation}\label{def.LBM1}
\LB^{\epsilon,x}_t= x+\,\int_0^t (-\ln \epsilon)^{-1/4} \, \epsilon^{-1} e^{-  X_\epsilon(\LB^{\epsilon,x}_u)} \,d\bar B_u\,.
\end{equation}%\begin{equation}\label{def.LBM1}
%\left\lbrace
%\begin{array}{l}
%\LB^{\epsilon,x}_{t=0}=x \\
%d\LB^{\epsilon,x}_t= (-\ln \epsilon)^{-1/4}e^{- X_\epsilon(\LB^{n,x}_t)+  \E[X_\epsilon(\LB^{\epsilon,x}_t)^2]} \,d\bar B_t
%\end{array} \right.
%\end{equation}
where $\bar B_t$ is a standard two-dimensional Brownian motion. 
\end{definition}
We stress that the fact that there is no drift term in the definition of the Brownian motion is typical from a scalar metric tensor in dimension $2$. By using the Dambis-Schwarz Theorem, one can define the law of the ${\rm LBM}_\epsilon$ as 
\begin{definition}\label{d.elbm2} For any $\epsilon>0$ fixed and $x\in\R^2$,
\begin{equation}\label{}
\LB^{\epsilon,x}_t=x + B_{\<{\LB^{\epsilon,x}}_t}\,,
\end{equation}
where $(B_r)_{r\geq 0}$ is a two-dimensional Brownian motion independent of the field $X$ and  the {\bf quadratic variation} $\<{\LB^{\epsilon,x}}$ of $\LB^{\epsilon,x}$ is defined as 
\begin{equation}\label{eq.qvdambis}
\<{\LB^{\epsilon,x}}_t := \inf \{s\geq 0\, :\, (-\ln\epsilon)^{1/2}\,\epsilon^{2}\int_{0}^s e^{ 2  X_\epsilon(x+B_u)} \,du \geq  t\} \,.
\end{equation}
\end{definition}

%Observe that \eqref{eq.qvdambis} amounts to saying that the increasing process $\langle \mathcal{B}^{\epsilon,x}\rangle:\R_+\to\R_+$ satisfies the differential equation:
%\begin{equation}\label{EDbeps}
%\langle \mathcal{B}^{\epsilon,x}\rangle_t=(-\ln\epsilon)^{-1/2} \int_0^t e^{-2 X_\epsilon(x+B_{\langle \mathcal{B}^{\epsilon,x}\rangle_u})+2 \ln\frac{1}{\epsilon}}\,du.
%\end{equation}
It will thus be useful to define the following quantity on $\R^2 \times \R_+$:
\begin{equation}\label{def.fn}
F^\epsilon(x,s) =\epsilon^{2} \int_0^s e^{ 2 X_\epsilon(x+B_u)} \,du, 
\end{equation}
in such a way that the process $\langle \mathcal{B}^{\epsilon,x}\rangle$ is entirely characterized by:
\begin{equation}\label{car.fn}
(-\ln\epsilon)^{1/2}F^\epsilon(x,\langle\mathcal{B}^{\epsilon,x}\rangle_t)=t.
\end{equation}

Several standard facts can be deduced from the smoothness of $X_\epsilon$. For each fixed $\epsilon>0 $, the  ${\rm LBM}_\epsilon$  a.s. induces a Feller diffusion on $\R^2$, and thus a semi-group $(P_t^{\epsilon})_{t\geq 0}$, which is symmetric w.r.t the volume form $M_\epsilon$.
 
%\begin{proposition}\label{prop.fellern}
%Let $\epsilon>0 $ be fixed. 
%Clearly, since $x\in \R^2\mapsto X_\epsilon(x)$ is a.s. (in $X$) a continuous function, the above $\epsilon$-regularized Liouville Brownian motion $\mathcal{B}^\epsilon$ a.s. induces a Feller diffusion on $\R^2$.  Let us denote by $(P_t^{\epsilon})_{t\geq 0}$ its semi-group. Also, $\mathcal{B}^\epsilon$ is reversible with respect to the Riemannian volume $M_{\epsilon}$, which is therefore invariant for $\mathcal{B}^\epsilon$.
%\end{proposition}
We will be mostly interested in establishing the convergence in law of the ${\rm LBM}_\epsilon$ as $\epsilon\to 0$.  Basically, studying the convergence of the ${\rm LBM}_\epsilon$ thus boils down to establishing the convergence of its quadratic variations $\langle\mathcal{B}^{\epsilon,x}\rangle$.

%%%%%%%%%%%%%%%%%%%%%%%%%%%%%%%%%%%%%%%%%%%%%%%%%%%%%%%%%%%%%%%%%%%%%%%%%%%%%%%%%%%%%%%%%%%%%%%%%%%%%%%
\section{Critical LBM starting from one fixed point}\label{sec.starLBM}
%%%%%%%%%%%%%%%%%%%%%%%%%%%%%%%%%%%%%%%%%%%%%%%%%%
%%%%%%%%%%%%%%%%%%%%%%%%%%%%%%%%%%%%%%%%%%%%%%%%%%
The first section is devoted to the convergence of the $\epsilon$-regularized LBM when starting from one fixed point, say $x\in\R^2$. As in the case of the convergence of measures \cite{Rnew7,Rnew12}, the critical situation here is technically   more complicated  than in the subcritical case \cite{GRV1}, though conceptually similar.  The first crucial step of the construction consists in establishing the convergence towards $0$ of the family of functions $(F^\epsilon(x,\cdot))_\epsilon$ and then in
computing the first order expansion of the maximum of the field $X_\epsilon$ along the Brownian path up to time $t$, and more precisely  to prove that
\begin{equation}
\max_{s\in [0,t]}X_\epsilon(x+B_s)-2\ln\frac{1}{\epsilon}\to-\infty,\quad \text{as }\epsilon\to 0. 
\end{equation}
This is mainly the content of subsection \ref{ss.convfield}, after some preliminary lemmas in subsection \ref{prelim}. Then our strategy will mainly be  to adapt the ideas related to convergence of critical measures \cite{Rnew7,Rnew12}.

 We further stress that, in the case of measures (see \cite{Rnew7}), the content of subsection \ref{ss.convfield} is established thanks to comparison with multiplicative cascades measures and Kahane's convexity inequalities. In our context, no equivalent result has been established in the context of multiplicative cascades in such a way that we have to carry out a direct proof.
%%%%%%%%%%%%%%%%%%%%%%%%%%%%%%%%%%%%%%%%%%%%%%%%%%%%%%%
\subsection{Preliminary results about properties of Brownian paths}\label{prelim}
%%%%%%%%%%%%%%%%%%%%%%%%%%%%%%%%%%%%%%%%%%%%%%%%%%%%%%%
Let us consider a standard Brownian motion $B$ on the plane $\R^2$ starting at some given point $x\in\R^2$. Let us consider the occupation measure $\mu_t$ of the Brownian motion up to time $t>0$ and let us define the function
\begin{equation}\label{def:h}
\forall \epsilon \in]0,1],\quad h(\epsilon)= \ln\frac{1}{\epsilon}\ln\ln\ln\frac{1}{\epsilon}.
\end{equation}
The following result is proved in \cite{legall2}  
\begin{theorem}\label{th:legall}
There exists a deterministic constant $c>0$ such that $\Pb^B$-almost surely, the set 
$$E=\{z\in \R^2;\limsup_{\epsilon\to 0} \frac{\mu_t(B(z,\epsilon))}{\epsilon^2h(\epsilon)}=c \}$$ has full $\mu_t$-measure.
\end{theorem}

We will need an extra elementary result about the structure of Brownian paths:
\begin{lemma}\label{browniancap}
For every $p>2$, we have almost surely:
$$\int_{\R^2\times\R^2}\frac{1}{|x-y|^2\ln\big(\frac{1}{|x-y|}+2\big)^p} \mu_t(dx)\mu_t(dy)<+\infty.$$
\end{lemma}

\vspace{2mm}
\noindent {\it Proof of Lemma \ref{browniancap}.} We have:
\begin{align*}
\E^B\Big[\int_{\R^2\times\R^2} \frac{1}{|x-y|^2\ln\big(\frac{1}{|x-y|}+2\big)^p} &\mu_t(dx)\mu_t(dy)\Big] 
=   \E^B\Big[\int_{[0,t]^2}\frac{1}{|B_r-B_s|^2\ln\big(\frac{1}{|B_r-B_s|}+2\big)^p} \,drds \Big]\\
=&\int_{[0,t]^2}\E^B\Big[\frac{1}{|r-s||B_1|^2\ln\big(\frac{1}{|r-s|^{1/2}|B_1|}+2\big)^p}  \Big]\,drds.
\end{align*}
Let us compute for $a\leq t$ the quantity $\E^B\Big[\frac{1}{a |B_1|^2\ln\big(\frac{1}{a^{1/2}|B_1|}+2\big)^p}  \Big]$. By using the density of the Gaussian law, we get:
\begin{align*}
\E^B\Big[& \frac{1}{a|B_1|^2\ln\big(\frac{1}{a^{1/2}|B_1|}+2\big)^p}  \Big]=\frac{1}{2\pi}\int_{\R^2}\frac{1}{a|u|^2\ln\big(\frac{1}{a^{1/2}|u|}+2\big)^p} e^{-\frac{|u|^2}{2}}\,du 
\leq   \int_{0}^\infty\frac{1}{ar\ln\big(\frac{1}{a^{1/2}r}+2\big)^p} e^{-\frac{r^2}{2}}\,dr\\
\leq &\frac{1}{a}\int_{0}^\infty\frac{1}{u\ln\big(\frac{1}{u}+2\big)^p} e^{-\frac{u^2}{2a}}\,du\\
\leq & \frac{1}{a}\int_{0}^{a^{1/4}}\frac{1}{u\ln\big(\frac{1}{u}+2\big)^p} e^{-\frac{u^2}{2a}}\,du+\frac{1}{a}\int_{a^{1/4}}^t\frac{1}{u\ln\big(\frac{1}{u}+2\big)^p} e^{-\frac{u^2}{2a}}\,du+\frac{1}{a}\int_{t}^\infty\frac{1}{u\ln\big(\frac{1}{u}+2\big)^p} e^{-\frac{u^2}{2a}}\,du\\
\leq & \frac{1}{a}\int_{0}^{a^{1/4}}\frac{1}{u\ln^p\big(\frac{1}{u}\big)}  \,du+\frac{1}{a}\int_{0}^t\frac{1}{u\ln\big(\frac{1}{u}+2\big)^p} e^{-\frac{1}{2a^{1/2}}}\,du+\frac{2}{  t^2\ln^p2}   e^{-\frac{t^2}{2a}}\\
\leq & \frac{4^p}{a(p-1)\ln^{p-1}\frac{1}{a}}+\frac{C}{a}  e^{-\frac{1}{2a^{1/2}}} +\frac{2}{  t^2\ln^p2}   e^{-\frac{t^2}{2a}}.
\end{align*}
Therefore
\begin{align*}
\E^B\Big[\int_{\R^2}&\int_{\R^2}\frac{1}{|x-y|^2\ln\big(\frac{1}{|x-y|}+2\big)^p} \mu_t(dx)\mu_t(dy)\Big]\\
\leq& C\int_{0}^t \int_{0}^t \Big( \frac{1}{|r-s| \ln^{p-1}\frac{1}{|r-s|}}+\frac{1}{|r-s|}  e^{-\frac{1}{2|r-s|^{1/2}}} +    e^{-\frac{t^2}{2|r-s|}}  \Big)\,drds. 
\end{align*}
As this latter quantity is obviously finite, the proof is complete.\qed
%%%%%%%%%%%%%%%%%%%%%%%%%%%%%%%%%%%%%%%%%%%%%%%%%%%%%%%
\subsection{First order expansion of the maximum of the field along Brownian paths}\label{ss.convfield}
%%%%%%%%%%%%%%%%%%%%%%%%%%%%%%%%%%%%%%%%%%%%%%%%%%%%%%%
To begin with, we claim:

\begin{proposition}\label{cv0} 
For all $x\in\R^2$, almost surely in $X$,  the family of random mapping $t\mapsto F^\epsilon(x,t)$ converges to $0$ in the space $C(\R_+,\R_+)$ as $\epsilon\to 0$.
\end{proposition}

\noindent {\it Proof.} Fix $x\in\R^2$.  Observe first that 
$$F^\epsilon(x,t)=\epsilon^{2}\int_0^te^{2X_\epsilon(x+B_r)}\,dr=\int_0^t e^{g_\epsilon(x+B_r)}e^{2X_\epsilon(x+B_r)-2\E^X[X_\epsilon(x+B_r)^2]}\,dr$$
where $$g_\epsilon(u)=2\E^X[X_\epsilon(u)^2]-2\ln\frac{1}{\epsilon}=\int_1^{1/\epsilon}\frac{k(u,x,x)-1}{u}\,du.$$
By assumption [A.4], the function $g_\epsilon$ converges uniformly over the compact subsets of $\R^2$. Furthermore, for each $t>0$, the set $\{x+B_s;s\in[0,t]\}$ is a compact set and $g_\epsilon$ converges uniformly over this compact set. So, even if it means considering $\int_0^t  e^{2X_\epsilon(x+B_r)-2\E^X[X_\epsilon(x+B_r)^2]}\,dr$ instead of $F^\epsilon(x,t)$, we may assume that $\E^X[X_\epsilon(x+B_r)^2]=\ln\frac{1}{\epsilon}$, in which case   $F^\epsilon(x,t)$ is a martingale with respect to the filtration $\mathcal{F}_\epsilon=\sigma\{X_r(x);\epsilon\leq r,x\in\R^2\}$. As this martingale is nonnegative,   it converges almost surely. We just have to prove that the limit is $0$. To this purpose, we use a lemma in \cite{Dur}. Translated into our context, it reads:
\begin{lemma}\label{durrettlem}
The almost sure convergence of the family $(F^\epsilon(x,t))_\epsilon$ towards $0$ as $\epsilon\to 0$ is equivalent to the fact that $\limsup_{\epsilon\to 0} F^\epsilon(x,t)=+\infty$ under the probability measure defined by:
$$Q_{|\mathcal{F}_\epsilon}=t^{-1} F^\epsilon(x,t)\,\Pb^X.$$
\end{lemma}
The main idea of what follows is to prove that, under $Q$, $F^\epsilon(x,t)$ is stochastically bounded from below by the exponential of a Brownian motion so that $\limsup_{\epsilon\to 0} F^\epsilon(x,t)=+\infty$ and we apply Lemma \ref{durrettlem} to conclude . 

To carry out this argument, let us define a new probability measure $\Theta_\epsilon$ on $\mathcal{B}(\R^2)\otimes\mathcal{F}_\epsilon $ by
\begin{equation}
\Theta_{\epsilon |\mathcal{B}(\R^2)\otimes\mathcal{F}_\epsilon}=t^{-1} e^{2X_\epsilon(y)-2\ln\frac{1}{\epsilon}}\,\Pb^X(d\omega)\mu_t(dy),
\end{equation}
where $\mu_t$ stands for the occupation measure of the Brownian motion $B^x$ starting from $x$. We denote by $ \E_{\Theta_\epsilon}$ the corresponding expectation. In fact, since the above definition defines a pre-measure on the ring $\mathcal{B}(\R^2)\otimes\bigcup_\epsilon \mathcal{F}_\epsilon$, one can define a measure $\Theta$ on $\mathcal{B}(\R^2)\otimes\mathcal{F}$ by using Caratheodory's extension theorem. We recover the relation $\Theta_{| \mathcal{B}(A)\otimes\mathcal{F}_\epsilon }=\Theta_\epsilon $.  Similarly, we construct   the probability measure $Q $ on $\mathcal{F}=\sigma\big(\bigcup_{\epsilon}\mathcal{F}_\epsilon\big)$ by setting:
$$Q_{ | \mathcal{F}_\epsilon}=t^{-1}F^\epsilon(x,t)\,d\Pb^X,$$ which is nothing but the marginal law of $(\omega,y )\mapsto \omega$ with respect to  $\Theta_\epsilon$. 
We state a few elementary properties below.
The conditional law of $y$ given $\mathcal{F}_\epsilon$ is given by:
$$\Theta_\epsilon(dy |\mathcal{F}_\epsilon)=\frac{e^{2X_\epsilon(y)-2\ln\frac{1}{\epsilon}}}{F^\epsilon(x,t)}\,\mu_t(dy).$$
If $Y$ is a $\mathcal{B}(\R^2)\otimes\mathcal{F}_\epsilon$-measurable random variable then it has the following conditional expectation given $\mathcal{F}_\epsilon$:
\begin{equation*}
\E_{\Theta_\epsilon}[Y | \mathcal{F}_\epsilon]=\int_{\R^2} Y(y,\omega)\frac{e^{2X_\epsilon(y)-2\ln\frac{1}{\epsilon}}}{F^\epsilon(x,t)}\,\mu_t(dy).
\end{equation*}

Now we turn to the proof of Proposition \ref{cv0} while keeping  in mind this preliminary background. Let us observe that it is enough to prove that the set $$\{\limsup_{\epsilon\to 0}F^\epsilon(x,t)=+\infty\}$$ has probability $1$ conditionally to $y$  under $ \Theta$   to deduce that it satisfies
$$Q\big(\{\limsup_{\epsilon\to 0}F^\epsilon(x,t)=+\infty\}\big)=1.$$
So we have to compute the law of $F^\epsilon(x,t)$ under $\Theta(\cdot |y)$.  Recall the definition of $h$ in \eqref{def:h}. We have:
\begin{align*}
F^\epsilon(x,t)&=\int_{\R^2}e^{2X_\epsilon(u)-2\ln\frac{1}{\epsilon}}\mu_t(du) \geq \int_{\R^2}e^{2X_\epsilon(u)-2\ln\frac{1}{\epsilon}}\ind_{B(y,\epsilon)}\mu_t(du).
\end{align*}
Let us now write
$$X_\epsilon(u)=\lambda_\epsilon(u,y)X_\epsilon(y)+Z_\epsilon(u,y)$$ where
$\lambda_\epsilon(u,y)=\frac{K_\epsilon(u,y)}{\ln \frac{1}{\epsilon}}$ and $Z_\epsilon(u,y)=X_\epsilon(u)-\lambda_\epsilon(u,y)X_\epsilon(y)$. Observe that the process $(Z_\epsilon(u,y))_{u\in\R^2}$ is independent of $X_\epsilon(y)$. Therefore
\begin{align*}
F^\epsilon(x,t)&\geq \int_{\R^2}e^{2X_\epsilon(u)-2\ln\frac{1}{\epsilon}}\ind_{B(y,\epsilon)}\mu_t(du)\\
&=e^{2X_\epsilon(y)-2\ln\frac{1}{\epsilon}}\int_{\R^2}e^{2 Z_\epsilon(u,y) +2(\lambda_\epsilon(u,y) -1) X_\epsilon(y) }\ind_{B(y,\epsilon)}\mu_t(du)\\
&=e^{2X_\epsilon(y)-4\ln\frac{1}{\epsilon}+\ln h(\epsilon)}\times\inf_{u\in B(y,\epsilon)}e^{2(\lambda_\epsilon(u,y) -1) X_\epsilon(y) }\times \frac{1}{\epsilon^{2}h(\epsilon)}\int_{\R^2}e^{2 Z_\epsilon(u,y)}\ind_{B(y,\epsilon)}\mu_t(du).
\end{align*}
Let us define
$$a_\epsilon(y)=\int_{\R^2}\ind_{B(y,\epsilon)}\mu_t(du).$$
With the help of the Jensen inequality, we deduce  
\begin{align*}
F^\epsilon(x,t) \geq &e^{2X_\epsilon(y)-4\ln\frac{1}{\epsilon}+\ln h(\epsilon)} \,\inf_{u\in B(y,\epsilon)}e^{2(\lambda_\epsilon(u,y) -1) X_\epsilon(y) }   \frac{a_\epsilon(y)}{\epsilon^{2}h(\epsilon)}  \exp\Big( \frac{1}{a_\epsilon(y)} \int_{B(y,\epsilon)}2  Z_\epsilon(u,y)\mu_t(du)\Big).
\end{align*}
Let us set
$$Y_\epsilon=a_\epsilon(y)^{-1}\int_{\R^2}2  Z_\epsilon(u,y)\ind_{B(y,\epsilon)}\mu_t(du).$$
Finally, for all $R>0$, we use the independence of $Y_\epsilon$ and $X_\epsilon(y)$ to get
\begin{align}
\Theta\Big(&\{\limsup_{\epsilon\to 0}F^\epsilon(x,t)=+\infty\}|y\Big) \\
\geq &\Theta\Big(\{\limsup_{\epsilon\to 0} e^{2X_\epsilon(y)-4\ln\frac{1}{\epsilon}+\ln h(\epsilon)}\times\inf_{u\in B(y,\epsilon)}e^{2(\lambda_\epsilon(u,y) -1) X_\epsilon(y) }\times  \frac{a_\epsilon(y)}{\epsilon^{2}h(\epsilon)}  \times\exp(-R)=+\infty\}|y\Big)\nonumber\\
&\times\Theta(Y_\epsilon\geq -R|y).\label{mult1}
\end{align}
Now we analyze the behaviour of each term in the above expression. 

First, notice that $\Theta(Y_\epsilon\geq -R|y)= \Pb^X(Y_\epsilon\geq -R)$ and that $Y_\epsilon$ is a centered Gaussian random variable under $\Pb^X$ with variance
\begin{align*}
a_\epsilon(y)^{-2}\int_{B(y,\epsilon)\times B(y,\epsilon)}\E^X\big[&  (X_\epsilon(u)-\lambda_\epsilon(u,y)X_\epsilon(y)) (X_\epsilon(u')-\lambda_\epsilon(u',y)X_\epsilon(y))\big]\mu_t(du)\mu_t(du').
\end{align*}
This quantity may be easily evaluated with assumption [A.2] and proved to be less than some constant $C$, which does not depend on $\epsilon$ and $y\in \{x+B_s;s\in [0,t]\}$. Therefore
\begin{align}
\Theta(Y_\epsilon\geq -R|y)\geq 1-\rho(R)\label{mult2}
\end{align}
for some nonnegative function $\rho$ that goes to $0$ as $R\to \infty$.

Second, from Theorem \ref{th:legall}, there exists a constant $c$ such that $\Pb^B$-almost surely, the set 
$$E=\{z\in \R^2;\limsup_{\epsilon\to 0} \frac{a_\epsilon(z)}{\epsilon^2h(\epsilon)}=c \}$$ has full $\mu_t$-measure. Since $E$ has full $\mu_t$ measure, even if it means extracting a random subsequence (only depending on $B$), we may assume that 
\begin{equation}\label{mult3}
\lim_{\epsilon\to 0} \frac{a_\epsilon(y)}{\epsilon^2h(\epsilon)}=c.
\end{equation}

Third, under $\Theta(\cdot|y)$, the process $X_\epsilon(y)-2\ln\frac{1}{\epsilon}$ is a Brownian motion, call it $\bar{B}$, in logarithmic time, i.e. 
$$ \bar{B}_{\ln \frac{1}{\epsilon}}= X_\epsilon(y)-2\ln\frac{1}{\epsilon}.$$
We further stress that $\bar{B}$ is independent from $B$, and thus from the random sequence $(a_\epsilon(y))_\epsilon$.  From the law of the iterated logarithm, we deduce that $\Theta(\cdot|y)$-almost surely
\begin{equation}\label{mult4}
\limsup_{\epsilon\to 0}  e^{2X_\epsilon(y)-4\ln\frac{1}{\epsilon}+\ln h(\epsilon)}=+\infty.
\end{equation}

Fourth, using assumption [A.2], it is readily seen that there exists some constant $C$ such that for all $\epsilon\in (0,1)$, $y\in\{x+B_s,s\in [0,t]\}$ and all $u\in B(y,\epsilon)$
\begin{equation*}
2|\lambda_\epsilon(u,y) -1|\leq C(-\ln \epsilon)^{-1}.
\end{equation*}
Since  the process $X_\epsilon(y)-2\ln\frac{1}{\epsilon}$ is a Brownian motion in logarithmic time under $ \Theta(\cdot|y)$, we deduce
\begin{equation}\label{mult5}
\liminf_{\epsilon\to 0} \inf_{u\in B(y,\epsilon)}e^{2(\lambda_\epsilon(u,y) -1) X_\epsilon(y) }\geq 1.
\end{equation}
By gathering \eqref{mult3}+\eqref{mult4}+\eqref{mult5}, we deduce
that, under $\Theta(\cdot|y)$:
\begin{equation}\label{mult6}
\limsup_{\epsilon\to 0} e^{2X_\epsilon(y)-4\ln\frac{1}{\epsilon}+\ln h(\epsilon)}\times\inf_{u\in B(y,\epsilon)}e^{2(\lambda_\epsilon(u,y) -1) X_\epsilon(y) }\times  \frac{a_\epsilon(y)}{\epsilon^{2}h(\epsilon)}  \times\exp(-R)=+\infty.
\end{equation}
By plugging \eqref{mult2}+\eqref{mult6} into \eqref{mult1}, we get for all $R>0$:
\begin{align*}
\Theta\Big(&\{\limsup_{\epsilon\to 0}F^\epsilon(x,t)=+\infty\}|y\Big) 
\geq   1-\rho(R).
\end{align*}
By choosing $R$ arbitrarily large, we complete the proof of Proposition \ref{cv0} with the help of 
Lemma \ref{durrettlem}.\qed

\begin{proposition}\label{coro:max} 
Almost surely in $X$, for all $x\in \R^2$, $\Pb^{B^x}$ almost surely, for all $t\geq 0$ we have:
\begin{equation}\label{supadur}
\sup_{\epsilon>0}\sup_{s\in[0,t]}X_\epsilon( B^x_s)-2\ln\frac{1}{\epsilon}+a\ln\ln\frac{1}{\epsilon}<+\infty
\end{equation} for all $a\in ]0,\frac{1}{4}[$.
\end{proposition}

\noindent {\it Proof.} Assume that the kernel $k(u,x,y)$ in \eqref{corrX} is given by $k(u(x-y))$ for some continuous covariance kernel $k$ with $k(0)=1$. It is then proved in \cite{Rnew7} that, for each fixed $a\in ]0,\frac{1}{4}[$, there exists a sequence $(C_n)_n$ of $\Pb^X$-almost surely finite random variables such that, 
\begin{equation}\label{aux:sup}
\sup_{\epsilon>0}\sup_{x\in B(0,n)}X_\epsilon(x )-2\ln\frac{1}{\epsilon}+a\ln\ln\frac{1}{\epsilon}<C_n. 
\end{equation}
Now for each fixed $x$, the Brownian motion $B^x$ has $\Pb^x$-almost surely continuous sample paths. 
Therefore, $\Pb^x$-almost surely, for any $t\geq 0$, we can find $n$ such that $\forall s\in [0,t]$, $B^x_s\in B(0,n)$. Thus the claim follows from \eqref{aux:sup} in this specific case.

One must make some extra effort to extend this result to the more general situation of assumption [A]. It will be convenient to set $t=\ln\frac{1}{\epsilon}$ for $\epsilon\in]0,1]$. For each $x\in\R^2$, we consider the mapping $$t\mapsto \varphi_x(t)=\E[X_{e^{-t}}(x)^2]=\int_1^{e^t}\frac{k(u,x,x)}{u}\,du.$$ It is continuous and strictly increasing and we denote by $\varphi_x^{-1}(t)$ the inverse mapping. Let us then consider the mapping $t\mapsto T_x(t)$ defined by $\varphi_x(T_x(t))=t$. We consider the Gaussian process $Y_t(x)=X_{e^{-T_x(t)}}(x)$, which has constant variance $t$. We have
\begin{align*}
\E[Y_t(x)Y_t(y)]&=\int_1^{e^{T_x(t)\wedge T_y(t)}}\frac{k(u,x,y)}{u}\,du\geq \int_1^{e^{T_x(t) }}\frac{k(u,x,y)}{u}\,du=\int_1^{e^t}\frac{k(e^{\varphi_x^{-1}(\ln v)},x,y)}{k(e^{\varphi_x^{-1}(\ln v)},x,x)}\,\frac{dv}{v}.
\end{align*}
By assumption [A.5], for each compact set $K$, we can find a constant $C_K$ such that $$\frac{k(u,x,y)}{k(u,x,x)}\geq (1-C_Ku^{1/2}|x-y|^{1/2})_+$$ for all $x\in K$ and $y\in\R^2$. We deduce, for all $x\in K$
\begin{align*}
\E[Y_t(x)Y_t(y)]& \geq \int_1^{e^t} \frac{(1-C_Ke^{\varphi_x^{-1}(\ln v)/2}|x-y|^{1/2})_+}{v}\,dv.
\end{align*}
By assumption [A.4], one can check that the mapping $x\mapsto \frac{e^{\varphi_x^{-1}(\ln v)}}{v}$ converges uniformly on $K$ towards a bounded strictly positive function. So even if it means changing the constant $C_K$, we have
\begin{align}\label{comparsup}
\E[Y_t(x)Y_t(y)]& \geq \int_1^{e^t} \frac{(1-C_K|x-y|^{1/2})_+}{v}\,dv.
\end{align}
for all $x\in K$ and $y\in\R^2$. From \cite{Rnew7}, this latter  covariance kernel satisfies the estimate \eqref{supadur}. From \cite{cf:Kah}, the above comparison between covariance kernels \eqref{comparsup} with equal variance entails that the result also holds for the process $Y_t$. It is then plain to conclude by noticing that the function $x\mapsto \frac{T_x(t)}{t}$ converges uniformly as $t\to \infty$ over the compact sets towards a strictly positive limit (this results from [A.4]). \qed

%%%%%%%%%%%%%%%%%%%%%%%%%%%%%%%%%%%%%%%%%%%%%%%%%%%%%%%
\subsection{Limit of the derivative PCAF}\label{ss.deriv}
%%%%%%%%%%%%%%%%%%%%%%%%%%%%%%%%%%%%%%%%%%%%%%%%%%%%%%%
Inspired by the construction of measures at criticality \cite{Rnew7,Rnew12}, it seems reasonable to think that the change of times $F^\epsilon$, when suitable renormalized, should converge towards a random change of times that coincides with the limit of the following  process
\begin{align*}
F^{',\epsilon}(x,t)& :=\int_0^t \big(2\ln \frac{1}{\epsilon}-X_{\epsilon}(x+B_u)\big) e^{ 2 X_\epsilon(x+B_u)- 2 \ln\frac{1}{\epsilon}} \,du .
\end{align*}
Establishing the convergence of the above martingale is the main purpose of this section.
Observe that, $\Pb^B$ almost surely,  the family $(F^{',\epsilon}(x,t))_\epsilon$ is  almost martingale for each $t>0$ (it is when you replace $\ln \frac{1}{\epsilon}$ by $\E^X[X_\epsilon(x)^2]$). Nevertheless, it is not nonnegative and not uniformly integrable. It is therefore not obvious that such a family almost surely converges towards a (non trivial) positive limiting random variable. The following theorem is the main result of this section:

\begin{theorem}\label{mainFderiv}
Assume [A.1-5] and fix $x\in\R^2$. For each $t>0$, the family $(F^{',\epsilon}(x,t))_\epsilon$ converges almost surely in $X$ and in $B$ as $\epsilon\to 0$ towards a  positive random variable denoted by $F'(x,t)$, such that $F'(x,t)>0$ almost surely. Furthermore, almost surely in $X$ and in $B$, the (non necessarily positive) random mapping $t\mapsto F^{',\epsilon}(x,t)$ converges  as $\epsilon\to 0$ in the space $C(\R_+,\R_+)$ towards a  strictly increasing continuous random mapping $t\mapsto F'(x,t)$. 
\end{theorem}

Throughout this section, we will assume that assumptions [A.1-5] are in force. The observation made in the beginning of the proof of Proposition \ref{cv0} remains valid here and, without loss of generality, we may assume that $\ln \frac{1}{\epsilon}=\E^X[X_\epsilon(x)^2]$. In this way, the family $(F^{',\epsilon}(x,t))_\epsilon$ is a martingale. Actually, 
Proposition \ref{coro:max} tells us that it is a positive martingale for $t$ large enough. Therefore it converges almost surely towards a limit $F'(x,t)$. But it is not uniformly integrable so that there are several complications involved in establishing   non-triviality of the limit.  We have to introduce some further tools to study the convergence.  We will introduce a family of auxiliary "truncated" martingales, called below $F^{',\epsilon}_\beta(x,t)$, which are reasonably close to $F^{',\epsilon}(x,t)$ while being square integrable. This will be enough to get the non triviality of 
$F'(x,t)$.

Given  $t>0$, $z,x\in\R^2$  and $\beta>0$, we  introduce the random variables
\begin{align}
f^\beta_\epsilon(z)=&\big(2\ln \frac{1}{\epsilon}-X_{\epsilon}(z)+\beta\big)\ind_{\{\tau^\beta_z< \epsilon\}}  e^{ 2 X_\epsilon(z)- 2 \ln\frac{1}{\epsilon}} \label{fbeta}\\
F^{',\epsilon}_\beta(x,t)=&\int_0^t \big(2\ln \frac{1}{\epsilon}-X_{\epsilon}(B^x_u)+\beta\big)\ind_{\{\tau^\beta_{B^x_u}< \epsilon\}}  e^{ 2 X_\epsilon(B^x_u)- 2 \ln\frac{1}{\epsilon}} \,du=\int_0^t f^\beta_\epsilon(B^x_u) \,du  \nonumber\\
\widetilde{F}^{',\epsilon}_\beta(x,t)=&\int_0^t \big(2\ln \frac{1}{\epsilon}-X_{\epsilon}(B^x_u)\big)\ind_{\{\tau^\beta_{B^x_u}< \epsilon\}}  e^{ 2 X_\epsilon(B^x_u)- 2 \ln\frac{1}{\epsilon}} \,du,\nonumber
\end{align}
where, for each $u\in [0,t]$, $\tau^\beta_u$ is the $(\mathcal{F}_\epsilon)_\epsilon$-stopping time  defined by
$$ \tau^\beta_z=\sup\{r\leq 1,X_r(z)-2\ln\frac{1}{r}>\beta \}.$$
 
In what follows, we will first investigate the convergence of  $(F^{',\epsilon}_\beta(x,t))_{\epsilon\in ]0,1]}$ to deduce first the convergence of $(\widetilde{F}^{',\epsilon}_\beta(x,t))_{\epsilon\in ]0,1]}$ and then the convergence of $(F^{',\epsilon}(x,t))_{\epsilon\in ]0,1]}$. We claim:

\begin{proposition}
We fix $x\in\R^2$ and $t>0$. Almost surely in $B$, the process $(F^{',\epsilon}_\beta(x,t))_{\epsilon \in ]0,1]}$ is a continuous positive $\mathcal{F}_\epsilon$-martingale  and thus converges almost surely in $X$ and $B$ towards a  nonnegative random variable denoted by  $F'_\beta(x,t)$. 
\end{proposition}

This proposition is a direct consequence of the  stopping time theorem and the martingale convergence theorem. Details are thus left to the reader. What is more involved is the study of the uniform integrability of this martingale: 
\begin{proposition}\label{prop:log}
We fix $x\in\R^2$ and $t>0$. Almost surely in $B$, the martingale $(F^{',\epsilon}_\beta(x,t))_{\epsilon >0}$ is uniformly integrable.
\end{proposition}

\noindent {\it Proof.}    Let us first state the following  lemma, which will serve in the forthcoming computations. 
%The proof of the first lemma can be found in \cite{Rnew7} and relies on the stopping time theorem and the fact that $(f^\beta_\epsilon(x))_{\epsilon\in]0,1]}$ is a martingale:

%\begin{lemma}\label{stoppingtime}
%Consider $x\neq w\in\R^2$.  
%$$\E[f^\beta_\epsilon(x)f^\beta_\epsilon(w)]=\E[f^\beta_{(|x-w|\vee \epsilon)\wedge 1}(x)f_{(|x-w|\vee \epsilon)\wedge 1}^\beta(w)].$$
%\end{lemma}

\begin{lemma}\label{occupation}
Let us denote $\mu_t$ the occupation measure of the Brownian motion $B^x_t$. $\Pb^{B^x}$-almost surely, for all $\delta>0$, there are a compact set $K$ and some constant $L$ such that $\mu_t(K^c)\leq \delta$ and for all $y\in K$
$$\sup_{r\leq 1}\frac{\mu_t(B(y,r))}{r^2g(r)}\leq L\quad \text{with }g(r)=\ln\big(\frac{1}{r}+2\big)^3.$$%\quad \text{with }h(r)=\ln\frac{1}{r}\ln\ln\ln\frac{1}{r}.$$ 
\end{lemma}

\noindent{\it Proof.} It is a direct consequence of Lemma \ref{browniancap}.\qed

\medskip

According to this lemma, for each fixed $\delta>0$, we are given a compact $K=K_\delta$ satisfying the above conditions. We denote by $\mu_t^K$ the measure $\mu_t^K(dy)=\ind_K(y)\,\mu_t(dy)$ and 
$$F^{',\epsilon}_{\beta,K}(x,t)=\int_0^t \big(2\ln \frac{1}{\epsilon}-X_{\epsilon}(y)+\beta\big)\ind_{\{\tau^\beta_{y}< \epsilon\}}  e^{ 2 X_\epsilon(y)- 2 \ln\frac{1}{\epsilon}} \ind_K(y)\,\mu_t(dy).$$
To prove Proposition \ref{prop:log}, it suffices to prove that the family $(F^{',\epsilon}_{\beta,K}(x,t))_\epsilon$ is uniformly integrable. Indeed, if true, we have for each $R>0$ ($K^c$ is the complement of $K$ in $\R^2$)
\begin{align*}
\Pb^X(F^{',\epsilon}_\beta(x,t)\ind_{\{F^{',\epsilon}_\beta(x,t)\geq R\}})\leq &\Pb^X(2F^{',\epsilon}_{\beta,K}(x,t)\ind_{\{2 F^{',\epsilon}_{\beta,K}(x,t)\geq R\}})+ \Pb^X(2F^{',\epsilon}_{\beta,K^c}(x,t)\ind_{\{2 F^{',\epsilon}_{\beta,K^c}(x,t)\geq R\}})\\
\leq &\Pb^X(2F^{',\epsilon}_{\beta,K}(x,t)\ind_{\{2 F^{',\epsilon}_{\beta,K}(x,t)\geq R\}})+2\E^X[F^{',\epsilon}_{\beta,K^c}(x,t)].
%=&\Pb^X(2F^{',\epsilon}_{\beta,K}(x,t)\ind_{\{2 F^{',\epsilon}_{\beta,K}(x,t)\geq R\}})+2\int_{K^c}\E^X[f^\beta_\epsilon(y)]\ind_{K^c}(y)\mu_t(dy)\\
%=&\Pb^X(2F^{',\epsilon}_{\beta,K}(x,t)\ind_{\{2 F^{',\epsilon}_{\beta,K}(x,t)\geq R\}})+2\beta\delta.
\end{align*}
We deduce 
$$\limsup_{R\to\infty }\Pb^X(F^{',\epsilon}_\beta(x,t)\ind_{\{F^{',\epsilon}_\beta(x,t)\geq R\}})\leq  2\E^X[F^{',\epsilon}_{\beta,K^c}(x,t)]=2\beta\delta.$$
By choosing $\delta$ arbitrarily small, we prove the uniform integrability of $(F^{',\epsilon}_\beta(x,t))_{\epsilon >0}$. 

So, we just have to focus on the uniform integrability of the family $(F^{',\epsilon}_{\beta,K}(x,t))_{\epsilon >0}$.   We  introduce the annulus  $C(y,\epsilon,1)=B(y,1)\setminus B(y,\epsilon)$ for $\epsilon\in (0,1)$. We get:
\begin{align}
\E^X[F^{',\epsilon}_{\beta,K}(x,t)^2]=&\int_{\R^2}\int_{\R^2}\E^X[f^\beta_\epsilon(y)f^\beta_\epsilon(w)]\,\mu_t^K(dy)\mu_t^K(dw)\nonumber\\
%=&\int_{\R^2}\int_{\R^2}\E^X[f^\beta_{(|y-w|\vee \epsilon)\wedge 1}(y)f^\beta_{(|y-w|\vee \epsilon)\wedge 1}(w)]\,\mu_t^K(dy)\mu_t^K(dw)\nonumber\\
=&\int_{\R^2}\int_{B(y,\epsilon)}\E^X[f^\beta_{ \epsilon }(y)f^\beta_{ \epsilon}(w)]\,\mu_t^K(dy)\mu_t^K(dw)\nonumber\\
&+\int_{\R^2}\int_{C(y,\epsilon,1)}\E^X[f^\beta_{\epsilon }(y)f^\beta_{\epsilon}(w)]\,\mu_t^K(dyx)\mu_t^K(dw)\nonumber\\
&+ \int_{\R^2}\int_{B(y,1)^c}\E^X[f^\beta_{\epsilon }(y)f^\beta_{\epsilon}(w)]\,\mu_t^K(dy)\mu_t^K(dw)\nonumber\\
\stackrel{def}{=}&\Pi^1_\epsilon+\Pi^2_\epsilon+\Pi^3_\epsilon.\label{nc1}
\end{align}
It is not difficult to see that $ \Pi^3_\epsilon\leq Ct^2$ for some constant $C$ independent of $\epsilon$. The main terms are $\Pi^1_\epsilon$ and $\Pi^2_\epsilon$. We begin with $\Pi^1_\epsilon$.
\begin{remark}
Before going further into details, let us just heuristically explain how to complete the proof. On the ball $B(y,\epsilon)$, the process $X_\epsilon(w)$ is very close to $X_\epsilon(y)$. Therefore, with a good approximation, we can replace $X_\epsilon(w)$ by $X_\epsilon(y)$ and get:
\begin{align*}
\Pi^1_\epsilon &\leq C\int_{\R^2}\int_{B(y,\epsilon)}\E^X\Big[(1+(X_\epsilon(y))^2)e^{2 X_\epsilon(y)+2\ln \frac{1}{\epsilon}}(\beta-X_\epsilon(y))\ind_{\{\sup_{s\in [\epsilon,1]}X_s(y)\leq \beta\}})\Big]\,\mu_t^K(dy)\mu_t^K(dw) .
\end{align*}
Let us define a new probability measure on $\mathcal{F}_\epsilon$ by $$\Pb^\beta(A)=\E^X[\ind_A(\beta-X_\epsilon(y))\ind_{\{\sup_{s\in [\epsilon,1]}X_s(y)\leq \beta\}})]$$ and recall that, under $\Pb^\beta$, the process $(\beta-X_s)_{\epsilon\leq s\leq 1}$ has the law of $(\beta_{\ln\frac{1}{s}})_{\epsilon\leq s\leq 1}$
where $(\beta_{u})_{u}$ is a $3$-dimensional Bessel process starting from $\beta$. Hence
\begin{align*}
\Pi^1_\epsilon 
&\leq C\int_{\R^2}\int_{B(y,\epsilon)}\E^\beta\Big[(1+(\beta_{\ln\frac{1}{\epsilon}})^2)e^{-2 \beta_{\ln\frac{1}{\epsilon}}+2\ln \frac{1}{\epsilon}} \Big]\,\mu_t^K(B(y,\epsilon))\mu_t^K(dy)\\
&\simeq C\int_{\R^2}(1+\ln \frac{1}{\epsilon})e^{-\sqrt{ \ln \frac{1}{\epsilon}}+2\ln \frac{1}{\epsilon}}\mu_t^K(B(y,\epsilon))\mu_t^K(dy)
\end{align*}
and this latter quantity goes to $0$ as $\epsilon\to 0$ since  $\mu_t^K(B(y,\epsilon))\leq L \epsilon^2g(\epsilon)$. Similar ideas will allow us to treat $\Pi^2_\epsilon$. Nevertheless, details are a bit more tedious.
\end{remark}

Let us now try to make rigorous the above remark.   Observe that necessarily $K_s(y,w)\leq \ln\frac{1}{s}$. Let us  define the functions $h_1$, $h_2$ and $\bar{h} $ by:
\begin{equation}\label{IBM1}
h_1(s)=\ln\frac{1}{s}-K_s(y,w)=h_2(s),\quad \bar{h}(s)=K_s(y,w).
\end{equation}
By considering $3$ independent Brownian motions $ B^1,B^2,\bar{B}$, we further define 
\begin{equation}\label{IBM2}
P^{y,w}_s=B^1_{h_1(s)},\quad P^{w,y}_s=B^2_{h_2(s)},\quad Z_s=\bar{B}_{\bar{h}(s)}.
\end{equation}
Observe that the process $(X_s(y),X_s(w))_{0\leq s\leq 1}$ has the same law as  the process $(P^{y,w}_s+Z_s,P^{w,y}_s+Z_s)_{0\leq s\leq 1}$. Now we compute $\Pi^1_\epsilon $ and then use a Girsanov transform:
\begin{align*}
\Pi^1_\epsilon 
=&\int_{\R^2}\int_{B(y,\epsilon)}\E^X[(\beta-P^{w,y}_\epsilon-Z_\epsilon-2\ln\frac{1}{\epsilon})\ind_{\{\sup_{r\in [\epsilon,1]}P^{w,y}_r+Z_r-2\ln\frac{1}{r}\leq \beta\}}(\beta-P^{y,w}_\epsilon-Z_\epsilon-2\ln\frac{1}{\epsilon})\times\dots\\
&\dots\ind_{\{\sup_{r\in [\epsilon,1]}P^{y,w}_r+Z_r-2\ln\frac{1}{r}\leq \beta\}}e^{2P^{y,w}_\epsilon+4Z_\epsilon+2P^{w,y}_\epsilon-4\ln\frac{1}{\epsilon}}]\,\mu_t^K(dy)\mu_t^K(dw)\\
=&\int_{\R^2}\int_{B(y,\epsilon)}\E^X[(\beta-P^{w,y}_\epsilon-Z_\epsilon)\ind_{\{\sup_{r\in [\epsilon,1]}P^{w,y}_r+Z_r\leq \beta\}}(\beta-P^{y,w}_\epsilon-Z_\epsilon)\times\dots\\
&\dots\ind_{\{\sup_{r\in [\epsilon,1]}P^{y,w}_r+Z_r\leq \beta\}}e^{2Z_\epsilon-2\ln\frac{1}{\epsilon}+4K_\epsilon(y-w)}]\,\mu_t^K(dy)\mu_t^K(dw).
\end{align*}
Let us set:
$$\beta^{y,w}_\epsilon=\beta-\min_{s\in[\epsilon,1]}P^{y,w}_s.$$
Because of assumption [A.2], we have $$\sup_{\epsilon\in (0,1]}\sup_{w\in B(y,\epsilon)}\sup_{s\in ]\epsilon,1]}h_1(s)+h_2(s)\leq c$$ for some constant $c>0$ only depending on $k$. Therefore we have:
\begin{align*}
\Pi^1_\epsilon 
\leq &C\int_{\R^2}\int_{B(y,\epsilon)}\E^X\Big[\Big(1+(Z_\epsilon)^2\Big)e^{2Z_\epsilon+2\ln\frac{1}{\epsilon} } (\beta^{y,w}_\epsilon-Z_\epsilon)  \ind_{\{\sup_{r\in [\epsilon,1]} Z_r\leq \beta^{y,w}_\epsilon\}}\Big]\,\mu_t^K(dy)\mu_t^K(dw).
\end{align*}
Let us define a new (random) probability measure on $\mathcal{F}_\epsilon$  by $$\Pb^{\beta,y,w}(A)=\frac{1}{\beta^{y,w}_\epsilon}\E^{Z}[\ind_A(\beta^{y,w}_\epsilon-Z_\epsilon(y))\ind_{\{\sup_{s\in [\epsilon,1]}Z_s(y)\leq \beta^{y,w}_\epsilon\}})|\beta^{y,w}_\epsilon]$$ with associated expectation denoted by $\E^{\beta,y,w} $. Recall that, under $\Pb^{\beta,y,w}$, the process $(\beta^{y,w}_\epsilon-Z_s)_{\epsilon\leq s\leq 1}$ has the law of $(\beta_{K_s(y,w)})_{\epsilon\leq s\leq 1}$
where $(\beta_{u})_{u}$ is a $3$-dimensional Bessel process starting from $\beta^{y,w}_\epsilon$. Hence:
\begin{align*}
\Pi^1_\epsilon 
\leq &C\int_{\R^2}\int_{B(y,\epsilon)}\E^X\Big[\beta^{y,w}_\epsilon\E^{\beta,y,w}\Big[\Big(1+(\beta_{K_\epsilon(y,w)})^2\Big)e^{-2\beta_{K_\epsilon(y,w)}+2\ln\frac{1}{\epsilon} }  \Big]\Big]\,\mu_t^K(dy)\mu_t^K(dw).
\end{align*}
Let us compute the quantity $$\E^{\beta,y,w}\Big[\Big(1+(\beta_{K_\epsilon(y,w)})^2\Big)e^{-2\beta_{K_\epsilon(y,w)}  }  \Big].$$
To this purpose, we use the fact that the law of a $3d$-Bessel process starting from $\beta^{y,w}_\epsilon$ is given by $\sqrt{(B^1_t-\beta^{y,w}_\epsilon)^2+(B^2_t)^2+(B^3_t)^2}$ where $B^1,B^2,B^3$ are three independent Brownian motions. Therefore, by using \eqref{property:k} when necessary, we get
\begin{align*}
\E^{\beta,y,w}\Big[&\Big(1+(\beta_{K_\epsilon(y,w)})^2\Big)e^{-2\beta_{K_\epsilon(y,w)} }  \Big]\\=&\int_{\R^3}(1+(u-\beta^{y,w}_\epsilon)^2+v^2+w^2)e^{-2\sqrt{(u-\beta^{y,w}_\epsilon)^2+v^2+w^2}}e^{-\frac{u^2+v^2+w^2}{2K_\epsilon(y,w)}}\,\frac{dudvdw}{(2\pi K_\epsilon(y,w))^{3/2}}\\
\leq &C(1+(\beta^{y,w}_\epsilon)^2)e^{2\beta^{y,w}_\epsilon}\int_{\R^3}(1+u^2+v^2+w^2)e^{-2\sqrt{u^2+v^2+w^2}}e^{-\frac{u^2+v^2+w^2}{2K_\epsilon(y,w)}}\,\frac{dudvdw}{(2\pi K_\epsilon(y,w))^{3/2}}\\
\leq & C(1+(\beta^{y,w}_\epsilon)^2)e^{2\beta^{y,w}_\epsilon}\int_{0}^\infty(1+r^2)e^{-2r}e^{-\frac{r^2}{2\ln\frac{1}{\epsilon}}}\,\frac{r^2 dr}{(  \ln\frac{1}{\epsilon})^{3/2}},
\end{align*}
for some constant $C>0$, which may have changed along lines. Let us set 
$$\forall a\geq 0,\quad H(a)=\int_{0}^\infty(1+r^2)e^{-2r}e^{-\frac{r^2}{2a}}\,\frac{r^2dr}{ a^{3/2}}.$$
It is plain to check that
$$H(a)\leq C(\max(1,a))^{-3/2}$$ for some positive constant $C$. Hence
\begin{align*}
\E^{\beta,y,w}\Big[&\Big(1+(\beta_{K_\epsilon(y,w)})^2\Big)e^{-2\beta_{K_\epsilon(y,w)} }  \Big] 
\leq   C(1+(\beta^{y,w}_\epsilon)^2)e^{2\beta^{y,w}_\epsilon}H(  \ln\frac{1}{\epsilon}).
\end{align*}
We deduce:
\begin{align*}
\Pi^1_\epsilon 
\leq &C(  \ln\frac{1}{\epsilon})^{-3/2} \int_{\R^2}\int_{B(y,\epsilon)}\E^X\Big[\beta^{y,w}_\epsilon  \big(1+(\beta^{y,w}_\epsilon)^2\big)e^{2\beta^{y,w}_\epsilon}\Big]\,\mu_t^K(dy)\mu_t^K(dw)\\
\leq &CH(  \ln\frac{1}{\epsilon}) \epsilon^{-2}\int_{\R^2} \mu_t^K(B(y,\epsilon))\,\mu_t^K(dy).
\end{align*}
Here we have used the fact  $\E^X\Big[ C\beta^{y,w}_\epsilon(1+(\beta^{y,w}_\epsilon)^2)e^{2\beta^{y,w}_\epsilon} \Big]$ is finite and does not depend on $y,w$ because $h_1(s)$ is bounded on $\R_+$ independently of $s,y,w$. Because of Lemma \ref{occupation}, this latter quantity goes to $0$ as $\epsilon\to 0$ $\Pb^{B^x}$-almost surely, and so does  $\Pi^1_\epsilon$.

Now we treat $\Pi^2_\epsilon$. We will follow similar arguments as for $\Pi^1_\epsilon$, though different behaviors are involved. Indeed, in this case, we have to face the possible long range correlations of the kernel $k$. So we adapt the decomposition of the couple $(X_s(y),X_s(w))_{s\in]0,1]}$ into Wiener integrals as follows. Let us consider a  smooth function $\varphi$ with compact support in the ball $B(0,1)$, such that $0\leq \varphi\leq 1$ and $\varphi=1$ over a neighborhood of $0$.
 Let us  define the functions $h_1$,  $h_2$, $\bar{h} $ and $\widehat{h}$ by:
\begin{align*} 
h_1(s)=&\ln\frac{1}{s}-K_s(y,w)=h_2(s), &\bar{h}(s)=&\int_1^{\frac{1}{s}} \frac{k(u,y,w)\varphi(u(y-w))}{u}\,du,\\
\widehat{h}(s)=&\int_1^{\frac{1}{s}} \frac{k(u,y,w)\big(1-\varphi(u(y-w))\big)}{u}\,du.
\end{align*}
By considering $4$ independent Brownian motions $ B^1,B^2,\bar{B},\widehat{B}$, we further define 
\begin{equation}
P^{y,w}_s=B^1_{h_1(s)},\quad P^{w,y}_s=B^2_{h_2(s)},\quad Z_s=\bar{B}_{\bar{h}(s)},\quad \widehat{Z}_s=\widehat{B}_{\widehat{h}(s)}.
\end{equation}
An elementary computation of covariance shows that the process $(X_s(y),X_s(w))_{0\leq s\leq 1}$ has the same law as  the process $(P^{y,w}_s+Z_s+\widehat{Z}_s,P^{w,y}_s+Z_s+\widehat{Z}_s)_{0\leq s\leq 1}$. The process $Z_s$ encodes the short-scale correlations of the two Brownian motions $(X_s(y))_s$ and $(X_s(w))_{s}$, and is the process that will rule the behaviour of   the expectation $\E^X[f^\beta_\epsilon(y)f^\beta_\epsilon(w)]$. The remaining terms are just negligible perturbations that we will have to get rid of in the forthcoming computations. 

We make first a few elementary remarks. Observe  that  $\bar{h}(s)=\bar{h}(|y-w|)$ for all $s\geq |y-w|$ in such a way that 
$Z_\epsilon=Z_{|y-w|}$ for $\epsilon\leq |y-w|$.
We also set $$D:=\sup\{\widehat{h}(s);s\in]0,1]; y,w\in \R^2\}<+\infty.$$ We will often use the elementary relation:
\begin{equation}\label{elem}
\forall a \geq 0,\forall x\in\R,\quad (\beta-a-x)\ind_{\{x\leq \beta-a\}}\leq (\beta-x)\ind_{\{x\leq \beta\}}.
\end{equation}

Now we begin the computations related to $\Pi^2_\epsilon$. 
So we consider $w\in C(y,\epsilon,1)$ and we have by the Girsanov transform and \eqref{elem}:
\begin{align*}
\E^X&[f^\beta_\epsilon(w)f^\beta_\epsilon(y)]\\=&\E^X\Big[(\beta-P^{y,w}_\epsilon-Z_\epsilon-\widehat{Z}_\epsilon+2\ln\frac{1}{\epsilon})\ind_{\{\sup_{u\in [\epsilon,1]}P^{y,w}_u+Z_u+\widehat{Z}_u-2\ln\frac{1}{u}\leq \beta\}}e^{2P^{y,w}_\epsilon+2Z_\epsilon+2\widehat{Z}_\epsilon-2\ln\frac{1}{\epsilon}}\dots\\
&\dots\times(\beta-P^{w,y}_\epsilon-Z_\epsilon-\widehat{Z}_\epsilon+2\ln\frac{1}{\epsilon})\ind_{\{\sup_{u\in [\epsilon,1]}P^{w,y}_u+Z_u+\widehat{Z}_u-2\ln\frac{1}{u}\leq \beta\}}e^{2P^{w,y}_\epsilon+2Z_\epsilon+2\widehat{Z}_\epsilon-2\ln\frac{1}{\epsilon}}\Big]\\
\leq&e^{8D}\E^X\Big[(\beta-P^{y,w}_\epsilon-Z_\epsilon-\widehat{Z}_\epsilon+2\ln\frac{1}{\epsilon})\ind_{\{\sup_{u\in [\epsilon,1]}P^{y,w}_u+Z_u+\widehat{Z}_u-2\ln\frac{1}{u}\leq \beta\}}e^{2P^{y,w}_\epsilon+2Z_\epsilon -2\ln\frac{1}{\epsilon}}\dots\\
&\dots\times(\beta-P^{w,y}_\epsilon-Z_\epsilon-\widehat{Z}_\epsilon+2\ln\frac{1}{\epsilon})\ind_{\{\sup_{u\in [\epsilon,1]}P^{w,y}_u+Z_u+\widehat{Z}_u-2\ln\frac{1}{u}\leq \beta\}}e^{2P^{w,y}_\epsilon+2Z_\epsilon-2\ln\frac{1}{\epsilon}}\Big]\\
\leq&e^{8D}\E^X\Big[(\beta-\min_{s\in]0,1]}\widehat{Z}_s-P^{y,w}_\epsilon-Z_\epsilon +2\ln\frac{1}{\epsilon})\ind_{\{\sup_{u\in [\epsilon,1]}P^{y,w}_u+Z_u-2\ln\frac{1}{u}\leq \beta-\min_{s\in]0,1]}\widehat{Z}_s\}}e^{2P^{y,w}_\epsilon+2Z_\epsilon -2\ln\frac{1}{\epsilon}}\dots\\
&\dots\times(\beta-\min_{s\in]0,1]}\widehat{Z}_s-P^{w,y}_\epsilon-Z_\epsilon+2\ln\frac{1}{\epsilon})\ind_{\{\sup_{u\in [\epsilon,1]}P^{w,y}_u+Z_u-2\ln\frac{1}{u}\leq \beta-\min_{s\in]0,1]}\widehat{Z}_s\}}e^{2P^{w,y}_\epsilon+2Z_\epsilon-2\ln\frac{1}{\epsilon}}\Big].
\end{align*}
The point is now to see that the above expectation reduces to the same expectation with $\epsilon$ replaced by $|y-w|$.  First observe that  $Z_\epsilon=Z_{|y-w|}$ for $\epsilon\leq |y-w|$. Second, from assumption [A.3], we have 
$$\sup_{y\in K,w\in\R^2}\int_{\frac{1}{|y-w|}}^{\infty}\frac{k(u,y,w)}{u}\,du=C_K<+\infty.$$
Therefore we can deduce $$\forall u\leq |y-w|,\quad  \big|\ln \frac{|y-w|}{u}- (h_1(u)-h_1(|y-w|))\big|\leq c$$ for some constant $c$ independent of everything that matters. This means that the quadratic variations of the martingale $(P^{w,y}_u-P^{w,y}_{|y-w|})_{u\leq |y-w|}$ can be identified with $\ln\frac{|y-w|}{\epsilon}$ up to some constant $c$ independent of $y,w,\epsilon$. We further stress that both martingales $P^{y,w}$ and $P^{w,y}$ are independent.
Therefore, by conditioning with respect to $\mathcal{F}_{|y-w|} $,  the integrand in the above expectation essentially reduces to  the product of two independent martingales (recall that, if $X_t=\int_0^tf(r)\,dB_r$ is a Wiener integral, then $(\beta+2\E[X_t^2]-X_t)\ind_{\{\sup_{s\in [0,t]}X_s-2\E[X_t^2]\leq \beta\}}e^{2X_t-2\E[X_t^2]}$  is a martingale). By applying the stopping time theorem and by setting $ \widehat{\beta}=\beta-\min_{s\in]0,1]}\widehat{Z}_s+c$, we get
\begin{align*}
\Pi^2_\epsilon 
\leq &C\int_{\R^2}\int_{C(y,\epsilon,1)}\E^X[( \widehat{\beta}-P^{w,y}_{|y-w|}-Z_{|y-w|}-2\ln\frac{1}{{|y-w|}})\ind_{\{\sup_{r\in [{|y-w|},1]}P^{w,y}_r+Z_r-2\ln\frac{1}{r}\leq  \widehat{\beta}\}}\\
&\dots( \widehat{\beta}-P^{y,w}_{|y-w|}-Z_{|y-w|}-2\ln\frac{1}{|y-w|})\ind_{\{\sup_{r\in [|y-w|,1]}P^{y,w}_r+Z_r-2\ln\frac{1}{r}\leq  \widehat{\beta}\}}\\
&\dots e^{2P^{y,w}_{|y-w|}+4Z_{|y-w|}+2P^{w,y}_{|y-w|}-4\ln\frac{1}{|y-w|}}]\,\mu_t^K(dy)\mu_t^K(dw).
\end{align*}
Recall that $\sup_{y,w\in K,s\in ]|y-w|,1]}\E[(P^{y,w}_{s})^2+(P^{w,y}_{s})^2]\leq c'$   and some constant $c'>0$ only depending on $k$ by assumption [A.2].   Indeed, for $s\in[|y-w|,1]$, we have: 
\begin{align*}
h_1(s)=&\ln\frac{1}{s}-\int_1^{\frac{1}{s}}\frac{k(u (y-w))}{u}\,du= \int_1^{\frac{1}{s}}\frac{1-k(u (y-w))}{u}\,du\\
\leq & \int_1^{\frac{1}{s}}\frac{C_Ku |y-w|}{u}\,du\leq   C_K,
\end{align*}
where $C_K$ is the Lipschitz constant given by assumption [A.2]. 
So, if we use the Girsanov transform and even if it means changing the value of $\widehat{\beta}$ by adding $C_K$, we get 
\begin{align*}
\Pi^2_\epsilon 
\leq &C\int_{\R^2}\int_{C(y,\epsilon,1)}\E^X[(\widehat{\beta}-P^{w,y}_{|y-w|}-Z_{|y-w|})\ind_{\{\sup_{r\in [|y-w|,1]}P^{w,y}_r+Z_r\leq\widehat{\beta}\}}(\widehat{\beta}-P^{y,w}_{|y-w|}-Z_{|y-w|})\\
&\dots\ind_{\{\sup_{r\in [|y-w|,1]}P^{y,w}_r+Z_r\leq \widehat{\beta}\}}e^{2Z_{|y-w|}-2\ln\frac{1}{|y-w|}+4K_{|y-w|}(y-w)}]\,\mu_t^K(dy)\mu_t^K(dw).
\end{align*}
Let us then define 
$$\beta^{y,w}=\widehat{\beta}-\min_{s\in[|y-w|,1]}P^{y,w}_{|y-w|}.$$
We deduce:
\begin{align*}
\Pi^2_\epsilon 
\leq &C\int_{\R^2}\int_{C(y,\epsilon,1)}\E^X\Big[\Big(1+(Z_{|y-w|})^2\Big)e^{2Z_{|y-w|}+2\ln\frac{1}{|y-w|} } (\beta^{y,w} -Z_{|y-w|})\times \dots\\
&\dots\times \ind_{\{\sup_{r\in [|y-w|,1]} Z_r\leq \beta^{y,w} \}}\Big]\,\mu_t^K(dy)\mu_t^K(dw).
\end{align*}
Let us define a new (random) probability measure on $\mathcal{F}_{|y-w|}$ by $$\Pb^{\beta,y,w}(A)=\frac{1}{\beta^{y,w}}\E^{Z}[\ind_A(\beta^{y,w} -Z_{|y-w|}(y))\ind_{\{\sup_{s\in [|y-w|,1]}Z_s(y)\leq \beta^{y,w}\}})|\beta^{y,w} ]$$ and recall that, under $\Pb^{\beta,y,w}$, the process $(\beta^{y,w} -Z_s)_{|y-w|\leq s\leq 1}$ has the law of $(\beta_{K_s(y-w)})_{|y-w|\leq s\leq 1}$
where $(\beta_{u})_{u}$ is a $3$-dimensional Bessel process starting from $\beta^{y,w}$. Hence
\begin{align*}
\Pi^2_\epsilon 
\leq &C\int_{\R^2}\int_{C(y,\epsilon,1)}\E^X\Big[\beta^{y,w} \E^{\beta,y,w}\Big[\Big(1+(\beta_{K_{|y-w|}(y-w)})^2\Big)e^{-2\beta_{K_{|y-w|}(y-w)}+2\ln\frac{1}{|y-w|} }  \Big]\Big]\,\mu_t^K(dy)\mu_t^K(dw).
\end{align*}
From \eqref{property:k}, we have
$$\ln\frac{1}{|y-w|}-C_K\leq K_{|y-w|}(y-w)\leq \ln\frac{1}{|y-w|}+C_K.$$ If we proceed along the same lines as we did previously, we get:
\begin{align}
\E^{\beta,y,w}\Big[&\Big(1+(\beta_{K_{|y-w|}(y-w)})^2\Big)e^{-2\beta_{K_{|y-w|}(y-w)}  }  \Big]\nonumber\\%=&\int_{\R^3}(1+(u-\beta^{y,w})^2+v^2+w^2)e^{-\sqrt{(u-\beta^{y,w})^2+v^2+w^2}}e^{-\frac{u^2+v^2+w^2}{2K_{|y-w|}(y-w)}}\,\frac{dudvdw}{(2\pi K_{|y-w|}(y-w))^{3/2}}\\
%\leq & C(1+(\beta^{y,w})^2)e^{\beta^{y,w}}\int_{\R^3}(1+u^2+v^2+w^2)e^{-\sqrt{u^2+v^2+w^2}}e^{-\frac{u^2+v^2+w^2}{4K_\epsilon(y-w)}}\,\frac{dudvdw}{(2\pi K_\epsilon(y-w))^{3/2}}\\
\leq & C(1+(\beta^{y,w})^2)e^{\beta^{y,w}}\int_{0}^\infty(1+r^2)e^{-r}e^{-\frac{r^2}{2\ln\frac{1}{|y-w|}}}\,\frac{r^2dr}{(  \ln\frac{1}{|y-w|})^{3/2}}\nonumber\\
\leq &C(1+(\beta^{y,w})^2)e^{\beta^{y,w}}H (  \ln\frac{1}{|y-w|}) .\label{nc2}
\end{align}
Therefore, thanks to Lemma \ref{occupation} (or Lemma \ref{browniancap}),
$$\sup_{\epsilon\in]0,1]}\Pi^2_\epsilon\leq C\int_{\R^2}\int_{B(y,1)}\frac{H(\ln\frac{1}{|y-w|})}{|y-w|^2}\mu_t^K(dw)\mu_t^K (dy)<+\infty.$$
%To sum up, we have proved that the martingale $\big(F^{',\epsilon}_{\beta,K}(x,t)\big)_\epsilon$ is uniformly integrable for each compact set $K=K_\delta$, where
%$$F^{',\epsilon}_{\beta,K}(x,t)=\int_{\R^2}f^\beta_\epsilon(u)\mu_t^K(du).$$
% This entails  uniform integrability of the martingale $(F^{',\epsilon}_\beta(x,t))_\epsilon$: for each $R>0$ and any $\delta>0$, we set $K=K_\delta$ and we have:
%\begin{align*}
%\limsup_{R\to\infty}\sup_{\epsilon\in (0,1]}&\E^X[F^{',\epsilon}_\beta(x,t)\ind_{\{F^{',\epsilon}_\beta(x,t)>R\}}]\\
%&\leq \limsup_{R\to\infty}\sup_{\epsilon\in (0,1]}2\E^X[F^{',\epsilon}_{\beta,K}(x,t)\ind_{\{2F^{',\epsilon}_{\beta,K}(x,t)>R\}}]+2\mu_t(K^c)\\
%&\leq 2\delta.
%\end{align*}
%By choosing $\delta$ arbitrarily close to $0$, 
The proof of Proposition \ref{prop:log} is complete.\qed

\vspace{2mm}
Before proceeding with the proof of   Theorem \ref{mainFderiv}, let us first state a few corollaries of the previous computations. For $\beta>0$ and $\epsilon\in]0,1]$, define the random measure $$M_\epsilon^\beta(dx)=f_\epsilon^\beta(x)\,dx.$$
Following \cite{Rnew7}, the family $(M^\beta_\epsilon)_{\epsilon\in]0,1]}$ almost surely converges as $\epsilon\to 0$ in the sense of weak convergence of measure towards a limiting non trivial measure  $M^\beta(dx)$ and  $M^\beta(dx)=M'(dx)$ on compact sets for $\beta$ (random) large enough. The following corollary proves  Theorem \ref{mainderiv} and therefore considerably generalizes the results in \cite{Rnew7}.

\begin{corollary}
Assume [A.1-5]. Consider the random measure $$M_\epsilon^\beta(dx)=f_\epsilon^\beta(x)\,dx.$$
Then for $p<1/2$:
$$\E^X\Big[\int_{B(0,1)}\int_{B(0,1)} \ln^p \frac{1}{|y-w|} M^\beta(dy)M^\beta(dw)\Big]<+\infty.$$
Therefore, the measure $M^\beta$ (and consequently $M'$) is diffuse.
\end{corollary} 

\noindent {\it Proof.} If we just replace the occupation measure of the Brownian motion by the Lebesgue measure  along the lines of the proof of Proposition \ref{prop:log}, we get:
\begin{align*}
\E^X&\Big[\int_{B(0,1)}\int_{B(0,1)}  \ln^p \frac{1}{|y-w|} M^\beta(dy)M^\beta(dw)\Big]\\
\leq &C  \int_{B(0,1)}\int_{B(0,1)} \frac{1}{|y-w|^2} \ln^p \frac{1}{|y-w|}H (  \ln\frac{1}{|y-w|})dy dw .%\\
%=& C  \int_{B(0,1)}\int_{B(0,1)}  \frac{1}{|y-w|^2}\ln^p \frac{1}{|y-w|}  \int_{0}^\infty(1+r^2\ln\frac{1}{|y-w|})e^{-r\ln^{1/2}\frac{1}{|y-w|}}e^{-\frac{r^2}{2}}\,r^2drdy dw \\
%\leq &C\int_0^1\int_0^\infty \frac{1}{\rho^2}\ln^p \frac{1}{\rho} (1+r^2\ln\frac{1}{\rho})e^{-r\ln^{1/2}\frac{1}{\rho}}e^{-\frac{r^2}{2}}\,r^2\,\rho drd \rho\\
%=&C\int_0^\infty\int_0^\infty  r^{-2p}y^{2p+1}  (1+y^2)e^{-y}e^{-\frac{r^2}{2}}\,  drd y.
\end{align*}
This latter quantity is finite for $p<1/2$. The finiteness of such an integral implies that, almost surely, the measure $M^{\beta}$ cannot give mass to singletons.\qed

This result is closely related, though weaker, than \eqref{modulus}. Yet, the setup of our proof is quite general whereas \eqref{modulus} has been proved in \cite{KNW}  for a specific one-dimensional measure that exhibits nice scaling relations. Nonetheless, \eqref{modulus} is expected to hold in greater generality but we do not know to which extent the proofs in \cite{KNW} extend to more general situations. In the same spirit, we claim:
\begin{corollary}\label{coro:cont}
Fix $t>0$ and $p<1/2$. For each $\delta>0$, there is a compact set $K$ such that $\mu_t(K^c)\leq \delta$ and 
$$\E^{X}\Big[\int_0^t\int_0^t  \ln^p_+\frac{1}{|B^x_r-B^x_s|}\ind_{K}(B^x_r)\ind_{K}(B^x_s)F'_\beta(x,dr)F'_\beta(x,ds)\Big]<+\infty.$$
In particular, almost surely in $X$ and in $B^x$, the random mapping $r\mapsto F'_\beta(x,r)$ does not possess discontinuity point on $[0,t]$. 
\end{corollary}

\noindent {\it Proof.}  For each $\delta>0$, we use once again Lemma \ref{occupation} to find  a compact set $K$ and some constant $L$ such that $\mu_t(K^c)\leq \delta$ and for all $y\in K$
\begin{equation}\label{cap:cont}
\sup_{r\leq 1}\frac{\mu_t(B(y,r))}{r^2g(r)}\leq L,\quad \text{with }g(r)=\ln\big(2+\frac{1}{r}\big)^3.
\end{equation} 
We denote by $\mu_t^K$ the measure
$$\mu_t^K(dy)=\ind_K(y)\,\mu_t(dy).$$
Once again, the computations made in proposition \ref{prop:log} show that
\begin{align*}
\E^{X}&\Big[ \int_0^t\int_0^t \ln^p_+\frac{1}{|B^x_r-B^x_s|}\ind_{K}(B^x_r)\ind_{K}(B^x_s)F'_\beta(x,dr)F'_\beta(x,ds)\Big]\\
\leq &C  \int_0^t\int_0^t \frac{\ind_{\{|B^x_r-B^x_s|\leq 1\}}}{|B^x_r-B^x_s|^2}\ln^p\frac{1}{|B^x_r-B^x_s|}H\Big(\ln\frac{1}{|B^x_r-B^x_s|}\Big)\ind_{K}(B^x_r)\ind_{K}(B^x_s)dr ds \\
&+  C \int_0^t\int_0^t   \ind_{\{|B^x_r-B^x_s|\geq 1\}}\ln^p_+\frac{1}{|B^x_r-B^x_s|}\ind_{K}(B^x_r)\ind_{K}(B^x_s)dr ds \\
\leq &C  \int_0^t\int_0^t \frac{\ind_{\{|u-v|\leq 1\}}}{|u-v|^2}\ln^p\frac{1}{|u-v|}H\Big(\ln\frac{1}{|u-v|}\Big)\mu_t^K(du)\mu_t^K(dv).% \\
%&+   \int_0^t\int_0^t   \beta^2\ln^p\frac{1}{|u-v|}\mu_t^K(du)\mu_t^K(dv).
\end{align*}
Because of \eqref{cap:cont}, the above integrals are finite for $p<1/2$. 

It is plain to deduce the continuity property of the mapping $s\mapsto F'_\beta(x,[0,s])$ since the Brownian motion has continuous sample paths. Indeed, the discontinuity points of this mapping corresponds to the set $\mathcal{A}$ of atoms of the measure $F'_\beta(x,ds)$, which are countable. For each $n\in\N^*$, let us denote by $\mathcal{A}_n$ the set of atoms in $[0,t]$ of this measure that are of size strictly greater than $\int_0^t\ind_{K_{1/n}^c}(B^x_s)F'_\beta(x,dr)$. For $s\in \mathcal{A}_n$, we necessarily have $\ind_{K_{1/n}^c}(B^x_s)=0$.  Also, we  necessarily have  $\ind_{K_{1/n}}(B^x_s)=0$ otherwise the integral  $ \int_0^t\int_0^t \ln^p\frac{1}{|B^x_r-B^x_s|}\ind_{K_{1/n}}(B^x_r)\ind_{K_{1/n}}(B^x_s)F'_\beta(x,dr)F'_\beta(x,ds)$ would be infinite. We deduce that $\mathcal{A}_n$ is empty for all $n\in\N^*$, meaning that there is no atom of size greater than $\int_0^t\ind_{K_{1/n}^c}(B^x_s)F'_\beta(x,dr)$ for all $n$. Since 
$$\E^x\Big[\int_0^t\ind_{K_{1/n}^c}(B^x_s)F'_\beta(x,dr)\Big]=\mu_t(K_{1/n}^c)\leq 1/n\to 0 \quad \text{as }n\to \infty,$$
we deduce that the quantity $\int_0^t\ind_{K_{1/n}^c}(B^x_s)F'_\beta(x,dr)$ converges to $0$ in probability as $n\to\infty$. We complete the proof.
\qed
 
\vspace{2mm}

We are now in position to handle the proof of Theorem \ref{mainFderiv}.

\noindent {\it Proof of Theorem \ref{mainFderiv}.} We first observe that the martingale $(F^{',\epsilon}_\beta(x,t))_{\epsilon >0}$ possesses almost surely the same limit as the process $(\widetilde{F}^{',\epsilon}_\beta(x,t))_{\epsilon >0}$
because
\begin{equation}\label{diffz}
|F^{',\epsilon}_\beta(x,t)-\widetilde{F}^{',\epsilon}_\beta(x,t)|=\beta\int_0^t\ind_{\{\tau^\beta_{B^x_u}< \epsilon\}} e^{2X_{\epsilon}(B^x_u)-2\ln\frac{1}{\epsilon}} du\leq \beta F^{\epsilon}(x,t)
\end{equation}
 and the last quantity converges almost surely towards $0$  (see Proposition \ref{cv0}). Using Corollary \ref{coro:max}, we have almost surely in $X$ and in $B$:
$$\sup_{\epsilon>0}\max_{s\in [0,t]}X_\epsilon(B^x_s)-2\ln\frac{1}{\epsilon}<+\infty,$$
which obviously implies
\begin{align*}
\forall \epsilon>0 ,\quad F^{',\epsilon}(t,x)=\widetilde{F}^{',\epsilon}_\beta(x,t)
\end{align*}
for $\beta$  (random) large enough. We deduce that, almost surely in $X$ and in $B$, the family $(F^{',\epsilon}(x,t))_{\epsilon> 0}$ converges towards a positive random variable.

It is plain to deduce the random measures $(F^{',\epsilon}_\beta(x,dt))_\epsilon$ converges in the sense of weak convergence of measures towards a random measure $F^\beta(x,dt)$. To prove convergence in $C(\R_+,\R_+)$, we just have to prove that the mapping  $t\mapsto F^\beta(t,x)$ is continuous. This property is proved in Corollary \ref{coro:cont}. From \eqref{diffz} and Proposition \ref{cv0} again, we deduce that the family of random mappings   $(t\mapsto F^{',\epsilon}(x,t))_\epsilon$ almost surely converges as $\epsilon\to 0$ in  $C(\R_+,\R_+)$ towards a nonnegative nondecreasing mapping $t\mapsto F'(x,t)$.

Let us prove that, almost surely in $X$ and in $B$, the mapping $t\mapsto F'(x,t)$ is strictly increasing. We first write the relation, for $\epsilon'<\epsilon$,
\begin{align}\label{passlim}
F^{',\epsilon'}_\beta(x,dr)=&(2\ln \frac{1}{\epsilon} -X_{\epsilon}(B^x_r)+\beta)\ind_{\{\tau^\beta_{B^x_r}<{\epsilon'}\}} e^{2X_{\epsilon'}(B^x_r)-2\ln\frac{1}{\epsilon'}}\,dr\\&+(2\ln\frac{\epsilon}{\epsilon'}-X_{\epsilon'}(B^x_r)+X_{\epsilon}(B^x_r)+\beta)\ind_{\{\tau^\beta_{B^x_r}<\epsilon'\}} e^{2X_{\epsilon'}(B^x_r)-2\ln\frac{1}{\epsilon'}}\,dr.\nonumber
\end{align}
By using the same arguments as throughout this section, we pass to the limit in this relation as $\epsilon'\to 0$ and then $\beta\to \infty$ to get
\begin{align}\label{equiv}
F' (x,dr)= &  e^{2X_{\epsilon}(B^x_r)-2\ln\frac{1}{\epsilon}}\,F'_\epsilon(x,dr)
\end{align}
where $F'_\epsilon(x,dr)$ is almost surely defined as
$$F'_\epsilon(x,dr)=\lim_{\beta \to \infty}\lim_{\epsilon' \to 0}F^{',\epsilon'}_{\beta,\epsilon} (x,dr)$$
and $F^{',\epsilon'}_{\beta,\epsilon} (x,dr)$ is given by
$$(2\ln\frac{\epsilon}{\epsilon'}-X_{\epsilon'}(B^x_r)+X_{\epsilon}(B^x_r)+\beta)\ind_{\{\tau^\beta_{\epsilon,B^x_r}< \epsilon'\}} e^{2(X_{\epsilon'}-X_\epsilon)(B^x_r)-2\ln\frac{\epsilon}{\epsilon'}}\,dr$$
where
$$\tau^\beta_{\epsilon,z}=\sup\{u\leq 1;X_{u\epsilon}(z)-X_{\epsilon}(z)-2\ln\frac{1}{u}>\beta-X_\epsilon(z)+2\ln\frac{1}{\epsilon}\}.$$
Let us stress that we have used the fact that the measure $$(2\ln\frac{1}{\epsilon}-X_{\epsilon}(B^x_r)+\beta)\ind_{\{\tau^\beta_{B^x_r}<\epsilon'\}} e^{2X_{\epsilon'}( B^x_r)-2\ln\frac{1}{\epsilon'}}\,dr
$$ goes to $0$ (it is absolutely continuous w.r.t. to $F^{\epsilon'}(x,dr)$) when passing to the limit in \eqref{passlim} as $\epsilon'\to 0$.   From \eqref{equiv}, it is plain to deduce that, almost surely in $B$, the event $\{F'(x,[s,t])=0\}$ (with $s<t$) belongs to the asymptotic sigma-algebra generated by the field $\{(X_\epsilon(x))_x;\epsilon>0\}$. Therefore it has probability $0$ or $1$ by the $0-1$ law of Kolmogorov. Since we have already proved that it is not $0$, this proves that almost surely in $B$, $\Pb^X(F'(x,[s,t])=0)=0$ for any $s<t$. By considering a countable family of intervals $[s_n,t_n]$ generating the Borel sigma field on $\R_+$, we deduce that, almost surely in $X$ and in $B$, the mapping $t\mapsto F'(x,t)$ is strictly increasing. \qed

\subsection{Renormalization of the change of times}\label{ss.renorm}
%%%%%%%%%%%%%%%%%%%%%%%%%%%%%%%%%%%%%%%%%%%%%%%%%%%%%%%%
Here we explain the Seneta-Heyde norming for the change of times $F^\epsilon$. Some technical constraints prevents us from claiming that it holds under the only assumptions [A.1-5]. So it is important to  stress here  that the Seneta-Heyde renormalization is not necessary to construct the critical LBM. It just illustrates that the derivative construction of the change of times $F'(x,t)$ also corresponds to a proper renormalization of $F^\epsilon$.

\begin{theorem}\label{th:seneta}
Assume [A.1-5] and either [A.6] or [A.6']. Then the conclusions of Theorem \ref{mainderiv2} hold and almost surely in $B$, we have the following convergence in $\Pb^X$-probability as $\epsilon\to 0$
$$\sqrt{\ln\frac{1}{\epsilon}}F^\epsilon(x,t)\to \sqrt{\frac{2}{\pi}}F'(x,t).$$
\end{theorem}

\noindent {\it Proof.} The proof  is a rather elementary adaptation of the proof in \cite[section D]{Rnew12}. Just reproduce the proof in \cite[section D]{Rnew12} while replacing the Lebesgue measure by the occupation measure of the Brownian motion and use Lemma \ref{occupation} when necessary. Details are  left to the reader. \qed

\begin{remark}
In the case of example  \eqref{ex:GFF}, it is necessary to consider the occupation measure of the Brownian motion killed upon touching the boundary of the domain $D$. 
\end{remark}

\subsection{The LBM does not get stuck}\label{ss.stuck}
%%%%%%%%%%%%%%%%%%%%%%%%%%%%%%%%%%%%%%%%%%%%%%%%%%

In this subsection, we make sure that the LBM does not get stuck in some area of the state space $\R^2$. Typically, this situation may happen over areas where the field $X$ takes large values, therefore having as consequence to slow down the LBM. Mathematically, this can be formulated as follows: check that the mapping $t\mapsto F'(x,t)$  tends to $\infty$ as $t\to\infty$. 

\begin{theorem}\label{cvF}
Assume [A.1-5] and fix $x\in\R^2$. Almost surely in $X$ and in $B$, 
\begin{equation}
\lim_{t\to\infty}F'(x,t)=+\infty.
\end{equation} 
\end{theorem} 

\ni
{\it Proof.} It suffices to reproduce the techniques of \cite[subsection 2.5]{GRV1}.\qed
 
\begin{remark}
Of course, this statement does not hold in the case of a bounded planar domain when one considers a Brownian motion killed upon touching the boundary of $D$. In this case, the Liouville Brownian motion (see below) will run until touching the boundary of $D$.
\end{remark}

%%%%%%%%%%%%%%%%%%%%%%%%%%%%%%%%%%%%%%%%%%%%%%%%%%%%%%%
\subsection{Defining the critical LBM when starting from a given fixed point}\label{ss.convquad}
%%%%%%%%%%%%%%%%%%%%%%%%%%%%%%%%%%%%%%%%%%%%%%%%%%%%%%%

We are now in position to define the critical LBM when starting from one fixed point. Indeed, once the change of times $F'$ has been constructed, the strategy is the same as in \cite[subsection 2.10]{GRV1}.

\begin{definition}\label{def:cvq} Assume [A.1-5]. The  {\bf critical Liouville Brownian motion} is defined  by:
\begin{equation}\label{def:lbm}
\forall t\geq 0,\quad \mathcal{B}^x_t=B^x_{\langle \LB^{x}\rangle_t},\quad  \text{ and }
\quad \forall t\geq 0,\quad \sqrt{2/\pi}\,F'(x,\langle \mathcal{B}^x\rangle_t)=t.
\end{equation}
As such, the mapping $t\mapsto \langle \mathcal{B}^x\rangle_t$ is defined on $\R_+$, continuous and strictly increasing.
\end{definition}

\begin{theorem}\label{th:cvq} Assume [A.1-5] and   either [A.6] or [A.6']. Fix $x\in\R^2$. Almost surely in $B$ and in $\Pb^X$-probability, the family $( B,\langle\mathcal{B}^{\epsilon,x}\rangle, \mathcal{B}^{\epsilon,x})_\epsilon$ converges in the space $C(\R_+,\R^2)\times C(\R_+,\R_+)\times C(\R_+,\R^2)$ equipped with the supremum norm on compact sets  towards the triple $(B,\langle\mathcal{B}^x\rangle,\mathcal{B}^x)$.
%\item Fix $x\in\R^2$. Almost surely in $B$ and in $\Pb^X$-probability, the family $(\mathcal{B}^{\epsilon,x} )_\epsilon$ converges in the space $C(\R_+,\R^2)$ equipped with the supremum norm on compact sets  towards a random process $\mathcal{B}^x$ starting from $x$, called the {\bf critical Liouville Brownian motion},   defined by:
%$$\forall t\geq 0,\quad \mathcal{B}^x_t=B^x_{\langle \LB^{x}\rangle_t},$$ and
%$$\forall t\geq 0,\quad \sqrt{2/\pi}\,F'(x,\langle \mathcal{B}^x\rangle_t)=t.$$
%\end{enumerate}
\end{theorem}
 
\begin{remark}
One may wonder whether the process introduced in Definition \ref{def.LBM1} also converges in probability. Actually, the argument carried out in  \cite[section 2.12]{GRV1} remains true here: almost surely in $X$, the couple of processes $(\bar{B},\LB^\epsilon)_\epsilon$  in Definition \ref{def.LBM1} converges in law towards a couple $(\bar{B},\LB)$ where $\bar{B}$  and   $\LB$    are independent. Since $\LB^\epsilon$ is measurable w.r.t. $\bar{B}$, this shows that the process $\LB^\epsilon$ does not converge in probability. This justifies our approach of studying the convergence via the Dambis-Schwarz representation theorem: it leads to studying a process (definition \ref{d.elbm2}) that  converges in probability.
\end{remark} 

\section{Critical LBM as a Markov process}\label{sec.markov}
%%%%%%%%%%%%%%%%%%%%%%%%%%%%%%%%%%%%%%%%%%%%%%%%%%
%%%%%%%%%%%%%%%%%%%%%%%%%%%%%%%%%%%%%%%%%%%%%%%%%%
\begin{remark}
In this whole section, we assume that assumptions [A.1-5] are in force. Sometimes, we make a statement to relate our results to the metric tensor $g_\epsilon$. For such a connection to be made, the Seneta-Heyde norming is needed and thus assumption [A.6] or [A.6'] are required. To avoid confusion, this will be explicitly mentioned.
\end{remark}

In this section, we will investigate the critical LBM as a Markov process, meaning that we aim at constructing almost surely in $X$  the critical LBM starting from every point. In the previous section, the guiding line was similar to \cite{GRV1} besides technical difficulties. From now on, the difference will be conceptual too: in the subcritical situation, the issue of constructing the LBM starting from every point is possible because the mapping 
$$x\mapsto \int_{\R^2}\ln_+\frac{1}{|x-y|}M_\gamma(dy)$$ is a continuous function of $x$, where $M_\gamma$ stands for the subcritical measure with parameter $\gamma<2$. This idea is somewhat underlying the theory of traces of Dirichlet forms developed in \cite{fuku} for instance. At criticality, the main obstacle is pointed out in  \eqref{modulus}: the best modulus of continuity that one may hope for $M'$ is of the type $$M'(B(x,r))\leq C\frac{1}{\sqrt{\ln(1+r^{-1})}}$$ and cannot be improved as proved in \cite{BKNSW} in the context of multiplicative cascades. The mapping 
\begin{equation}\label{mappingmp}
x\mapsto \int_{\R^2}\ln_+\frac{1}{|x-y|}M'(dy)
\end{equation}
 is thus certainly not continuous and may even take infinite values. The ideas of \cite{GRV1} need be renewed to face the issues of criticality.

On the other hand,  the theory of Dirichlet forms  \cite{fuku,revuzyor} (or potential theory) tells us that one can construct a Positive Continuous Additive Functional (PCAF for short) associated to $M'$ provided that the mapping \eqref{mappingmp} does not take too many infinite values. More precisely, $M'$ is required not to  give mass to polar sets of the Brownian motion. The problem of the theory of Dirichlet forms is that one can guarantee the existence of the PCAF but it cannot be identified and all the information that we get by explicitly constructing the PCAF $F'$ is lost in this approach. Also, while being extremely powerful in the description of the Dirichlet form of the LBM,  the theory of Dirichlet forms gives much weaker results than the coupling approach developed in \cite{GRV1,GRV2} concerning the qualitative/quantitative properties of $F'$. We thus definitely need to gather both of these approaches.

Our strategy will be to identify a large set of points of finiteness of the mapping \eqref{mappingmp} in order to construct a perfectly identified PCAF on the whole space via coupling arguments. Then we will prove that $M'$ does not charge polar sets in order to identify our PCAF with that of the theory of Dirichlet forms in the sense of the Revuz correspondence. Once this gap is bridged, we can apply the full machinery of \cite{fuku} to get a lot of further information about $F'$: mainly, a full description of the Dirichlet form associated to the critical LBM.

\subsection{Background on positive continuous additive functionals and Revuz measures}
%%%%%%%%%%%%%%%%%%%%%%%%%%%%%%%%%%%%%%%%
To facilitate the reading of our results, we  summarize here some basic notions of potential theory applied to the standard Brownian $(\Omega,(B_t)_{t \geq 0},(\mathcal{F}_t)_{t \geq 0},(\Pb^x)_{x \in \R^2})$ in $\R^2$ seen as a Markov process, which is of course reversible for the canonical volume form $dx$ of $\R^2$. These notions can be found with further details in \cite{fuku,revuzyor}. 
One may then consider the classical notion of capacity associated to the Brownian motion. In this context, we have the following definition:
\begin{definition}[Capacity and polar set] 
The capacity of an open set $O\subset \R^2$ is defined by
$${\rm Cap}(O)=\inf\{\int_D|f(x)|^2\,dx+\int_D|\nabla f(x)|^2\,dx;f\in H^1(\R^2,dx),\,\,f \geq 1\text{ over }O \}.$$
The capacity of a Borel measurable set $K$ is then defined as: 
\begin{equation*}
{\rm Cap}(K)= \underset{O \text{open}, K \subset O}{\inf}{\rm Cap}(O). 
\end{equation*}
The set $K$ is said polar when ${\rm Cap}(K)=0$.
\end{definition}

\begin{definition}[Revuz measure] 
A Revuz measure $\mu$ is a Radon measure on $\R^2$ which does not charge the polar sets.
\end{definition}

Then we introduce the notion  of PCAF that we will use in the following (see \cite{fuku,revuzyor}):

\begin{definition}[PCAF] \label{defPCAF}
A Positive Continuous Additive Functional   $(A(x,t,B))_{t \geq 0,x\in\R^2/N}$ of the Brownian motion (with $B_0=0$) on $\R^2$ is defined by:\\
-a polar set $N$  (for the standard Brownian motion),\\ 
- for each $x\in\R^2/N$ and $t\geq 0$, $A(x,t,B)$ is  $\mathcal{F}_t$-adapted and continuous, with values in $[0,\infty]$ and $A(x,0)=0$,\\
- almost surely,  
\begin{equation*}
 A(x,t+s,B)-A(x,t,B)= A (x+B_t,s,B_{t+\cdot}-B_t), \quad s,t \geq 0.
\end{equation*} 
 \end{definition}

In particular, a PCAF  is defined for all starting points $x \in \R^2$ except possibly on a polar set for the standard Brownian motion. One can also work with a PCAF starting from \textbf{all} points, that is when the set $N$ in the above definition can be chosen to be empty. In that case, the PCAF is said {\it in the strict sense}.

Finally, we conclude with the following definition on the support of a PCAF:
\begin{definition}[support of a PCAF] 
Let $(A(x,t,B))_{t \geq 0,x\in\R^2/N}$ be a PCAF with associated polar set $N$. The support of $(A(x,t,B))_{t \geq 0,x\in\R^2/N}$ is defined by:
\begin{equation*}
 \tilde{Y}= \Big\{x \in \R^2  \setminus N: \: \Pb(R(x)=0)=1\Big\},
\end{equation*}
where $R(x)=\inf \lbrace t>0: \: A(x,t,B)>0 \rbrace$.
\end{definition}

From section 5 in \cite{fuku}, there is a one to one correspondence between Revuz measures $\mu$ and PCAFs $(A(x,t,B))_{t \geq 0,x\in\R^2/N}$ under the Revuz correspondence: for any $t>0$ and any nonnegative Borel functions $f,h$:
\begin{equation}\label{revuzcorr}
\E_{h.dx}\Big[\int_0^tf(B_r)\,dA_r\Big]=\int_0^t\int_{\R^2}f(x)P_rh(x)\,\mu(dx)\,dr,
\end{equation}
 where $(P_r)_{r\geq 0}$ stands for the semigroup associated to the planar Brownian motion on $\R^2$.

\subsection{Capacity properties of the critical measure}
%%%%%%%%%%%%%%%%%%%%%%%%%%%%%%%%%%%%%%%%%%%%%%%%%%
The purpose of this section is to establish some preliminary results in order to apply the theory of Dirichlet forms. In particular, we will establish that the critical measure does not charge polar sets. To this purpose, we will need some fine pathwise properties of Bessel processes. So we first recall a result from \cite{Motoo}:
\begin{theorem}\label{th:motoo}
Let $X$ be a $3d$-Bessel process on $\R_+$ starting from $x\geq 0$ with respect to the law $\Pb_x$.
\begin{enumerate}
\item Suppose that $\phi \uparrow \infty$ such that $\int_1^{\infty}\frac{\phi(t)^3}{t}e^{-\frac{1}{2}\phi(t)^2}\,dt<+\infty$. Then
$$\Pb_x\Big(X_t>\sqrt{t}\phi(t)\,\text{ i.o. as } t\uparrow+\infty\Big)=0 .$$
\item Suppose that $\psi \downarrow  0$ such that $\int_1^{\infty}\frac{\psi(t)}{t} \,dt<+\infty$. Then
$$\Pb_x\Big(X_t<\sqrt{t}\psi(t)\,\text{ i.o. as }  t\uparrow+\infty\Big)=0 .$$
\end{enumerate}
\end{theorem}

Recall that we denote by $M_\beta$ the measure $M_\beta(dx)=\lim_{\epsilon\to 0}f^\beta_\epsilon(x)\,dx$. The main purpose of this subsection is to prove the following result:

\begin{theorem}\label{th:cap}
$\bullet$ Let us consider two functions $\phi,\psi$ as in the above theorem and define 
$$Y_\epsilon(x)=  2\ln\frac{1}{\epsilon}- X_\epsilon(x).$$
 We introduce the set 
$$E=\{x\in \R^2; \limsup_{\epsilon\to 0}\frac{Y_\epsilon(x)}{\sqrt{\ln\frac{1}{\epsilon}}\phi(\ln\frac{1}{\epsilon})}\leq 1\,\,\,\text{ and }\,\,\,\liminf_{\epsilon\to 0}\frac{Y_\epsilon(x)}{\sqrt{\ln\frac{1}{\epsilon}}\psi(\ln\frac{1}{\epsilon})}\geq 1\}.$$ 
 Then $M^\beta$ gives full mass to $E$, i.e. $M^\beta (E^c)=0$. \\
 $\bullet$ Furthermore, for each ball $B$, for each $\delta>0$, there is a compact set $K\subset  B$ such that $M_\beta(B\cap K^c)\leq \delta$ and for all $p>0$
$$\int_K\int_B\Big( \ln\frac{1}{|x-y|}\Big)^pM_\beta(dx)M_\beta(dy)<+\infty.$$
\end{theorem}

\noindent {\it Proof.} Let us consider a non empty ball $B$. We introduce the Peyri\`ere probability measure $Q^\beta$ on $\Omega\times B$:
$$\int f(x,\omega)\,dQ^\beta =\frac{1}{\beta|B|}\E^X[\int_Bf(x,\omega)M^\beta(dx)].$$
Let us consider the random process $t\mapsto \bar{X}_t =X_{e^{-t}} $. The main idea is that under $Q^\beta$, the process $(\beta-\bar{X}_t+2t)_{t\geq 0}$ has the law very close to a $3d$-Bessel process $(\beta_{t})_t$ starting from $\beta$. The first claim then follows from Theorem \ref{th:motoo}. We just have to  precise  the notion of "very close": some negligible terms appear because of the difference between $t$ and $\E^X[\bar{X}_t(x)^2]$ and we have to quantify them.

We consider the measure (convergence is established the same way as for $M_\beta$) 
$$\bar{M}_\beta(dx)=\lim_{\epsilon\to 0}(2\E^X[\bar{X}_t(x)^2]-X_t(x)+\beta) \ind_{\{\sup_{u\in [0,t]} \bar{X}_u-2\E^X[\bar{X}_u(x)^2]\leq \beta\}}e^{2\bar{X}_t-2\E^X[\bar{X}_t(x)^2]}\,dx$$ and we set $D=\sup_{t\geq 0}\sup_{x\in B}|\E^X[\bar{X}_t(x)^2]-t|<+\infty$. We set $H(x)=\lim_{t\to \infty}\E^X[\bar{X}_t(x)^2]-t$ (see assumption [A.4]). Observe that
$$M_\beta(dx)\leq e^{H(x)}\bar{M}_{\beta+2D}(dx).$$
Let us set $\hat{\beta}=\beta+2D$. Therefore, we consider the probability measure $\bar{Q}^{\hat{\beta}}$ on $\Omega\times B$:
$$\int f(x,\omega)\,d\bar{Q}^{\hat{\beta}} =\frac{1}{\hat{\beta}|B|}\E^X[\int_Bf(x,\omega)\bar{M}^{\hat{\beta}}(dx)]$$ 
and $Q^\beta$ is absolutely continuous with respect to $\bar{Q}^{\hat{\beta}}$. Under $\bar{Q}^{\hat{\beta}}$ (following arguments already seen), the process $(\hat{\beta}-\bar{X}_t+2\E^X[\bar{X}_t(x)^2])_{t\geq 0}$ has the law of $({\rm Bess}_{T_t(x)})_t$ where $({\rm Bess}_{t})_t$ a $3d$-Bessel process  starting from $\hat{\beta}$ and $T_t(x)=K_{e^{-t}}(x,x)$. The first claim then follows from Theorem \ref{th:motoo}, the fact that $\sup_{x\in B,t\geq 0}|K_{e^{-t}}(x,x)-t|<+\infty$ and the absolute continuity of $Q^\beta$ w.r.t. $\bar{Q}^{\hat{\beta}}$.

 Now we prove  the second statement, which is more technical.  To simplify things a bit, we assume that $\E[X_t(x)^2]=t$. If this not the case, we can apply the same strategy as above, namely considering $\bar{M}_{\beta+2D}$ instead of $M_\beta$. We consider now a couple of functions $\phi(t)=(1+t)^\chi$ and $\psi(t)=(1+t)^{-\chi}$ for some small positive parameter $\chi$ close to $0$. They satisfy the assumptions of Theorem \ref{th:motoo}. Then we consider the random compact sets  for $R>0$ and $t\in [0,+\infty]$,
\begin{align*}
K^1_{R,t}=&\big\{ x\in B;\sup_{u \in [0,t]}\frac{\beta-\bar{X}_u(x)+2u}{(1+u)^{1/2+\chi}}\leq R\big\}\quad &K^2_{R,t}= \big\{ x\in B;\inf_{u \in [0,t]}\frac{\beta-\bar{X}_u(x)+2u}{(1+u)^{1/2-\chi}}\geq \frac{1}{R}\big\}\\
K_{R,t}=&K^1_{R,t}\cap K_{R,t}^2.
\end{align*}
We will write $K_R$ for $K_{R,\infty}$. From Theorem \ref{th:motoo}, we have $\lim_{R\to\infty}Q^\beta(K_R)=1$. Therefore, we have $\lim_{R\to \infty}M^\beta(K_R^c\cap B)=0$ $\Pb^X$ almost surely.  
%%%%%%%%%%%%%%%%%%%%%%%%%%%%

Let us denote $\E^Q$ expectation with respect to the probability measure $Q^\beta$. Without loss of generality, we assume that $\beta|B|=1$: this avoids to repeatedly write the renormalization constant appearing in the definition of $Q^\beta$. To prove the result, we  compute 
\begin{align*}
\E^Q[ \ind_{K_R}(x) M^\beta(B(x,e^{-t}))]=&\lim_{\epsilon\to 0}\E[ \int_B\ind_{K_R}(x)M_\beta(B(x,e^{-t}))M_\beta(dx)]\\
=&\lim_{\epsilon\to 0}\E[\int_B\int_{B(x,e^{-t})}\ind_{K_R}(x)f^\beta_{\epsilon}(w)\,dw  M_\beta(dx)]\\
\leq &\lim_{\epsilon\to 0}\int_B\int_{B(x,e^{-t})}\E[\ind_{K_{R,-\ln|x-w|}}(x)f^\beta_{\epsilon}(w)f^\beta_{\epsilon}(x)]\,dwdx.
\end{align*}
Now we can argue as in the proof of Proposition \ref{prop:log} (treatment of $\Pi_\epsilon^2$)
to see that for $\epsilon\leq |x-w|$
$$\E[\ind_{K_{R,-\ln|x-w|}}(x)f^\beta_{\epsilon}(w)f^\beta_{\epsilon}(x)]\leq C\E[\ind_{K_{R,-\ln|x-w|}}(x)f^{\hat{\beta}}_{|x-w|}(w)f^{\hat{\beta}}_{|w-x|}(x)]$$ 
for some irrelevant constant $C$ and $\hat{\beta}=\beta-\min_{u\geq \ln \frac{1}{|x-w|}} \hat{Z}_u$ where the process $\hat{Z}$ is independent of the sigma algebra $\mathcal{F}_{|x-w|}$ and  $\hat{Z}_s=B_{\hat{h}(s)}-B_{\hat{h}(-\ln|x-w|)}$  for some Brownian motion $B$ and $\hat{h}(s)=\int_{0}^{e^s}\frac{k(u,x,w)(1-\varphi(u(x-w)))}{u}\,du$.
\begin{remark}
The above technical part is due to the presence of possible long range correlations. Though they do not affect qualitatively our final estimates, getting rid of them may appear technical. The reader who wishes to skip this technical part   may instead consider the compact case: $k(u,x,v)=0$ if $u|x-w|\geq 1$. In that case, we just have 
$$\E[\ind_{K_{R,-\ln|x-w|}}(x)f^\beta_{\epsilon}(w)f^\beta_{\epsilon}(x)]=\E[\ind_{K_{R,-\ln|x-w|}}(x)f^{ \beta}_{|x-w|}(w)f^{ \beta}_{|w-x|}(x)]$$ because of the stopping time theorem and the fact that both martingales $(f^\beta_{\epsilon}(w)-f^\beta_{|x-w|}(w))_{\epsilon\leq |x-w|}$ and $(f^\beta_{\epsilon}(x)-f^\beta_{|x-w|}(x))_{\epsilon\leq |x-w|}$ are   independent.
\end{remark}
We are thus left with computing $\E[\ind_{K_{R,-\ln|x-w|}}(x)f^{\hat{\beta}}_{|x-w|}(w)f^{\hat{\beta}}_{|w-x|}(x)]$. To this purpose, we need the following lemma, the proof of which is straightforward as a simple computation of covariance. Thus, details are left to the reader.  
 
\begin{lemma}
We  consider $w\not = x$ such that $|x-w|\leq  1$. The law of the couple $(\bar{X}_t(w),\bar{X}_t(x))_{t\geq 0}$ can be decomposed as:
$$(\bar{X}_t(w),\bar{X}_t(x))_{t\geq 0}=(Z^{x,w}_t+P^{x,w}_t,P^{x,w}_t+D^{x,w}_t)_{t\geq 0}$$
where the process $Z^{x,w},P^{x,w},D^{x,w}$ are independent centered Gaussian and, for some independent Brownian motions $B^1,B^2$ independent of $ \bar{X}$: 
\begin{align*}
P^{x,w}_t=&\int_0^tg'_{x,w}(u)\,d\bar{X}_u(x),& Z^{x,w}_t= \int_0^t(1-g'_{x,w}(u)^2)^{1/2}\,dB^1_u,  & &
 D^{x,w}_t= \int_0^t(1-g'_{x,w}(u)^2)^{1/2}\,dB^2_u.& &
 \end{align*}
 with $g_{x,w}(u)=\int_1^{e^u}\frac{k(y(x-w))}{y}\,dy$.
 Moreover 
\begin{equation}\label{supvar}
\sup_{t\leq \ln\frac{1}{|x-w|} }\E[(Z^{x,w}_t)^2+(D^{x,w}_t)^2]\leq C,
\end{equation} for some constant $C$ independent of $x,w$ such that $|x-w|\leq 1$.
\end{lemma}
Setting $s_0=\ln\frac{1}{|x-w|}$, we use this lemma to get
\begin{align*}
 \E[ \ind_{K_{R,-\ln |x-w|}}&(x)f^{\hat{\beta}}_{|w-x|}(w) f^{\hat{\beta}}_{|w-x|}(x)]\\
\leq &\E\Big[\ind_{K_{R,s_0}}(x)({\hat{\beta}}-P^{x,w}_{s_0}-Z^{x,w}_{s_0}-2s_0)_+e^{2(P^{x,w}_{s_0}+Z^{x,w}_{s_0})-2s_0}({\hat{\beta}}-P^{x,w}_{s_0}-D^{x,w}_{s_0}-2s_0)_+\times\dots \\
& \quad \ind_{\{\sup_{u\in[0,s_0]}P^{x,w}_{u}+D^{x,w}_{u}-2u\leq {\hat{\beta}}\}}e^{2(P^{x,w}_{s_0}+D^{x,w}_{s_0})-2s_0}  \Big]. 
\end{align*} 
We get rid of the process $Z^{x,w}$ by using first the Girsanov transform with the corresponding exponential term, and then by estimating the remaining terms containing  $Z^{x,w}$ with the help of \eqref{supvar}. We get for some constant $C$ that may vary along lines but does not depend on $x,w$:
 \begin{align*}
 \E[&\ind_{K_{R,s_0}}(x)f^{\hat{\beta}}_{|w-x|}(w) f^\beta_{|w-x|}(x)]\\
\leq &C\E\Big[\ind_{K_{R,s_0}}(x)(1+{\hat{\beta}})(1+|P^{x,w}_{s_0}-2\E[(P^{x,w}_{s_0})^2]|)e^{2P^{x,w}_{s_0}-2\E[(P^{x,w}_{s_0})^2]}(\beta-P^{x,w}_{s_0}-D^{x,w}_{s_0}-2s_0)_+\times\dots \\
&\dots\times\ind_{\{\sup_{u\in[0,s_0]}P^{x,w}_{u}+D^{x,w}_{u}-2u\leq {\hat{\beta}}\}}e^{2(P^{x,w}_{s_0}+D^{x,w}_{s_0})-2s_0}  \Big]. 
\end{align*} 
We use the Girsanov transform again to make the term $e^{2(P^{x,w}_{s_0}+D^{x,w}_{s_0})-2\ln\frac{1}{|x-w|}}$ disappear:
  \begin{align*}
 \E[&\ind_{K_{R,s_0}}(x)f^{\hat{\beta}}_{|w-x|}(w) f^{\hat{\beta}}_{|w-x|}(x)]\\
\leq &C\E\Big[\ind_{K_{R,s_0}}(x)(1+{\hat{\beta}})(1+|P^{x,w}_{s_0}|)e^{2P^{x,w}_{s_0}+2\E[(P^{x,w}_{s_0})^2]}({\hat{\beta}}-P^{x,w}_{s_0}-D^{x,w}_{s_0})_+\ind_{\{\sup_{u\in[0,s_0]}P^{x,w}_{u}+D^{x,w}_{u}\leq {\hat{\beta}}\}}  \Big]. 
\end{align*} 
Now we write $P^{x,w}_{s_0}=P^{x,w}_{s_0}+D^{x,w}_{s_0}-D^{x,w}_{s_0}$ and use \eqref{supvar} to see that we can replace $\E[(P^{x,w}_{s_0})^2] $ by $s_0$ even if it means modifying the constant $C$ (which still does not depend on $x,w$):  
  \begin{align*}
 \E[&\ind_{K_{R,s_0}}(x)f^{\hat{\beta}}_{|w-x|}(w) f^{\hat{\beta}}_{|w-x|}(x)]\\
\leq &C\E\Big[\ind_{K_R}(x)(1+{\hat{\beta}})(1+|P^{x,w}_{s_0}+D^{x,w}_{s_0}|+|D^{x,w}_{s_0}|)e^{2(P^{x,w}_{s_0}+ D^{x,w}_{s_0})+2s_0}e^{-2D^{x,w}_{s_0}}\times\dots \\
&\quad ({\hat{\beta}}-P^{x,w}_{s_0}-D^{x,w}_{s_0})_+\ind_{\{\sup_{u\in[0,s_0]}P^{x,w}_{u}+D^{x,w}_{u}\leq {\hat{\beta}}\}}  \Big]. 
\end{align*} 
Now we use the fact that for $x\in K_{R,s_0}$:
 \begin{align*}
\sup_{u\leq -\ln|x-w|} \frac{\beta-P^{x,w}_{u}-D^{x,w}_{u}}{(1+u)^{1/2+\chi}}\leq R & &\text{and} & & 
\inf_{u\leq -\ln|x-w|}  \frac{\beta-P^{x,w}_{u}-D^{x,w}_{u}}{(1+u)^{1/2-\chi}}\geq \frac{1}{R} ,
\end{align*}
 to get
   \begin{align*}
 \E[&\ind_{K_{R,s_0}}(x)f^{\hat{\beta}}_{|w-x|}(w) f^{\hat{\beta}}_{|w-x|}(x)]\\
\leq &C\E\Big[\ind_{K_{R,s_0}}(x)(1+{\hat{\beta}})\big(1+ R(1+s_0)^{1/2+\chi}+|D^{x,w}_{s_0}|\big)e^{-2R^{-1} (1+s_0)^{1/2-\chi}    +2s_0}e^{-2D^{x,w}_{s_0}}\times\dots\\
&\quad ({\hat{\beta}}-P^{x,w}_{s_0}-D^{x,w}_{s_0})_+\ind_{\{\sup_{u\in[0,s_0]}P^{x,w}_{u}+D^{x,w}_{u}\leq {\hat{\beta}}\}}  \Big]\\
\leq & \frac{C}{|x-w|^2}e^{-2R^{-1} s_0^{1/2-\chi}}(1+ Rs_0^{1/2+\chi})\times\dots 
\\&\dots\times \E\Big[(1+{\hat{\beta}})(1+ |D^{x,w}_{s_0}|)e^{-2D^{x,w}_{s_0}}  (\beta-P^{x,w}_{s_0}-\min_{u\in[0,s_0]}D^{x,w}_{u})_+\ind_{\{\sup_{u\in[0,s_0]}P^{x,w}_{u} \leq \beta-\min_{u\in[0,s_0]}D^{x,w}_{u}\}}  \Big]\\
\leq & \frac{C}{|x-w|^2}e^{-2R^{-1} s_0^{1/2-\chi}}(1+ Rs_0^{1/2+\chi}) \E\Big[(1+\hat{\beta})(1+ |D^{x,w}_{s_0}|)e^{-2D^{x,w}_{s_0}}  ({\hat{\beta}}-\min_{u\in[0,s_0]}D^{x,w}_{u})   \Big].
\end{align*} 
Because of \eqref{supvar}, the last expectation is finite and bounded independently of $x,w$. Indeed, $\hat{\beta}$, $(D^{x,w}_{u})_{u\leq s_0}$ are independent Wiener integrals with bounded variance (independently of $x,w$). Therefore we can find $\alpha>0$ such that $\sup_{|x-w|\leq 1}\E[e^{\alpha (\min_{u\leq s_0} D^{x,w}_{u})^2}+e^{\alpha \hat{\beta}^2}]<+\infty$. With these estimates, it is plain to see that the above expectation is finite.
 To sum up, we have proved
 \begin{align*}
\E^Q[&\ind_{K_R}(x) M^\beta(B(x,e^{-t}))]\\
\leq &\int_{B(x,e^{-t})}\frac{C}{|x-w|^2}e^{-2R^{-1}(\ln\frac{1}{|x-w|})^{1/2-\chi}}(1+ R(\ln\frac{1}{|x-w|})^{1/2+\chi}) \,dw\\
=&C\int_0^{e^{-t}}\rho^{-1}e^{-2R^{-1}(\ln\frac{1}{\rho})^{1/2-\chi}}(1+ R(\ln\frac{1}{\rho})^{1/2+\chi} )\,d\rho\\
=&C\int_t^{\infty} e^{-2R^{-1}y^{1/2-\chi}}(1+ Ry^{1/2+\chi}) \,dy.
\end{align*}
Finally we have
 \begin{align*}
\E\Big[ \int_{K_R}\int_{B}&\Big( \ln\frac{1}{|x-y|}\Big)^pM_\beta(dx)M_\beta(dy)\Big]\\
=&\E^Q\Big[\ind_{K_R}(x)\int_{B} \Big( \ln\frac{1}{|x-y|}\Big)^pM_\beta(dy) \Big]\\
 \leq &\sum_{n=1}^\infty \E^Q\Big[\ind_{K_R}(x)\int_{B\cap \{2^{-n-1}<|x-y|\leq 2^{-n}\}} \Big( \ln\frac{1}{|x-y|}\Big)^pM_\beta(dy) \Big]\\
  \leq &\sum_{n=1}^\infty (n+1)^p\ln^p2\, \E^Q\Big[\ind_{K_R}(x)M_\beta(B(x,2^{-n})) \Big]\\
 \leq &C\sum_{n=1}^\infty (n+1)^p \, \int_{n\ln 2}^{\infty} e^{-2R^{-1}y^{1/2-\chi}}(1+ Ry^{1/2+\chi}) \,dy.
  \end{align*}
This last series is easily seen to be finite. The proof of Theorem \ref{th:cap} is complete.\qed

\vspace{1mm}  
It is then a routine trick to deduce 
\begin{corollary}\label{th:capsupbeta}
For each ball $B$, for all $\delta>0$, there is a compact set $K_\delta\subset B$ such that $$M_\beta(B\cap K^c_\delta)\leq \delta$$ and for all $p>0$, for all $x\in K_\delta$
$$ \sup_{\substack{x\in \mathcal{O}\text{ open }\subset B,\\ {\rm diam}(\mathcal{O})\leq 1}}M_\beta( \mathcal{O})\big(-\ln{\rm diam}(\mathcal{O})\big)^{p}<+\infty.$$
\end{corollary}

\begin{corollary}\label{th:capsup}
\begin{enumerate}
\item For each ball $B$, for each $\delta>0$, there is a compact set $K\subset  B$ such that $M'(B\cap K^c)\leq \delta$ and for all $p>0$
$$\int_K\int_B\Big( \ln\frac{1}{|x-y|}\Big)^pM'(dx)M'(dy)<+\infty.$$
\item For each ball $B$, for all $\delta>0$, there is a compact set $K_\delta\subset B$ such that $M'(B\cap K^c_\delta)\leq \delta$  and for all $p>0$, for all $x\in K_\delta$
$$ \sup_{\substack{x\in \mathcal{O}\text{ open }\subset B,\\ {\rm diam}(\mathcal{O})\leq 1}}M'(  \mathcal{O})\big(-\ln{\rm diam}(\mathcal{O})\big)^{p}<+\infty.$$
\item Almost surely in $X$, the Liouville measure $M'$ does not charge the polar sets of the (standard) Brownian motion.
\end{enumerate}
\end{corollary} 

\vspace{1mm}
\noindent {\it Proof.} Both statements results from the fact that $M_\beta$ coincides with $M'$ over bounded sets for $\beta$ (random) large enough.\qed 

\vspace{1mm}
Though we will not need such a strong statement in the following , we state the optimal modulus of continuity bound that we can get with our methods
\begin{corollary}\label{th:optimal}
For each ball $B$, for all $\chi\in]0,1/2[$ and for all $\delta>0$, there is a compact set $K_\delta\subset B$ such that $M'(B\cap K^c_\delta)\leq \delta$  and for all $x\in K_\delta$
$$ \sup_{\substack{x\in \mathcal{O}\text{ open }\subset B,\\ {\rm diam}(\mathcal{O})\leq 1}}M'(  \mathcal{O})\exp\big((-\ln{\rm diam}(\mathcal{O}))^{\frac{1}{2}-\chi}\big)<+\infty.$$
\end{corollary} 

\begin{remark}
Observe that the question of the capacity properties of the measure $M'$, i.e. Corollary \ref{th:capsup}, was initially raised in \cite{KNW} (see at the end of the first section).
\end{remark}

\subsection{Defining $F'$ on the whole of  $\R^2$ }\label{whole}
%%%%%%%%%%%%%%%%%%%%%%%%%%%%%%%%%
In this subsection, we have two main objectives: to construct the PCAF $F^{'}(x,\cdot) $ on the whole of $\R^2$ and prove the convergence of $(F^{',\epsilon}_\beta(x,\cdot))_\epsilon$ towards $F'(x,\cdot)$. The main difficulty here is the following: in section \ref{ss.deriv} we have  proved the almost sure (in $X$ and $B$) convergence of $(F^{',\epsilon}_\beta(x,\cdot))_\epsilon$ towards $F$' when the starting point $x$ is fixed. Of course, we can deduce that almost surely in $X$, for a countable collection of given starting points, this convergence holds. The main difficulty is to prove that this convergence holds for {\bf all} possible starting points and this definitely requires some further arguments.

For $\mathbf{q}=(q_1,q_2)\in \Z^2$, we denote by $C_\mathbf{q}$ the cube $[q_1,q_1+1]\times[q_2,q_2+1]$. We fix $p>1$ and for each $\mathbf{q}\in \Z^2$ and $\delta>0$, we denote by $K^\mathbf{q}_{\delta,L}$ the compact set  
\begin{equation}
K^\mathbf{q}_{\delta,L}=K_\delta\cap\big\{x\in C_\mathbf{q};  \sup_{\substack{x\in\mathcal{O}\text{ open },\\ {\rm diam}(\mathcal{O})\leq 1}}M_\beta( \mathcal{O})\big(-\ln{\rm diam}(\mathcal{O})\big)^{p}\leq L\big\},
\end{equation}
where $K_\delta$ is the compact set given by Corollary \ref{th:capsupbeta} applied with $B=C_\mathbf{q}$. Then we set 
\begin{equation}
K^\mathbf{q}_{\delta}=\bigcup_{L>0}K^\mathbf{q}_{\delta,L},\quad K_{\delta,L}=\bigcup_{\mathbf{q}\in\Z^2}K^\mathbf{q}_{\delta,L}\quad \text{ and }\quad S=\bigcup_{\delta>0,\mathbf{q}\in\Z^2,L>0}K^\mathbf{q}_{\delta,L}.
\end{equation}
From Corollary \ref{th:capsupbeta}, we have $M_\beta((K^\mathbf{q}_{\delta})^c\cap C_\mathbf{q})\leq \delta$. Therefore, $M_\beta((\bigcup_{\delta>0}K^\mathbf{q}_{\delta})^c\cap C_\mathbf{q})=0$ for each $\mathbf{q}\in\Z^2$ and thus 
\begin{equation}
M_\beta(S^c)=0.
\end{equation}
 
 We consider a Brownian motion $B$ starting from $0$ and define a Brownian motion starting from $x$ by $B^x=x+B$ for each point $x\in\R^2$.  Following Section \ref{sec.starLBM}, we may
assume that $(F^{',\epsilon}_\beta(x,\cdot))_\epsilon$ converges almost surely in $X$ and $B$ in $C(\R_+;\R_+)$ towards $F'_\beta(x,\cdot)$ for each rational points $x\in\Q^2$.  

For each $\delta,L>0$ and $x\in\Q^2$, we define the adapted continuous random mapping
\begin{equation}
F^{',\delta,L}_{\beta}(x,t)=\int_0^t\ind_{K_{\delta,L}}(B^x_r)F'_\beta(x,dr).
\end{equation}

\begin{proposition}\label{propadapted}
Almost surely in $X$, for each   $\delta,L>0$ and $x\in \R^2$, there exists a $B^x$-adapted continuous random mapping, still  denoted by $F^{',\delta,L}_{\beta}(x,\cdot)$ such that for all sequence of rational points $(x_n)_n$ converging towards $x$, %we can extract a subsequence $(x_{\phi(n)})_n$ such that 
the sequence $(F^{',\delta,L}_{\beta}(x_{n},\cdot))_n$ converges in $\Pb^B$-probability in $C(\R_+,\R_+)$ towards $F^{',\delta,L}_{\beta}(x,\cdot)$.  
\end{proposition}

\noindent {\it Proof.} Let us fix $x\in\R^2$. We consider a sequence $(x_n)_n$ of rational points converging towards $x$. We first establish the convergence in law under $\Pb^B$ of the sequence $(B^{x_n},F^{',\delta,L}_{\beta}(x_n,\cdot))_n$ in $C(\R_+,\R^2)\times C(\R_+,\R_+)$. To this purpose, the main idea is an adaptation of \cite[section 2.9]{GRV1} with minor modifications. Yet, we outline the proof because we will play with this argument throughout this section. For all $0<s<t$, we write $F^{',\delta,L}_{\beta}(x,]s,t])$ for $F^{',\delta,L}_{\beta}(x, t)-F^{',\delta,L}_{\beta}(x,s)$. The proof relies on two arguments: a coupling argument and the estimate:
\begin{equation}\label{logzero}
\sup_{y\in B(0,R)\cap\Q^2}\E^B[F^{',\delta,L}_{\beta}(y,t)]\to 0,\quad \text{ as }t\to 0.
\end{equation}
We begin with explaining \eqref{logzero}. For all $y\in\Q^2$, we have:
\begin{align*}
\E^B[F^{',\delta,L}_{\beta}(y,t)]=\int_{\R^2}\int_0^tp_r(y,z)\,dr\ind_{K_{\delta,L}}(z)M_\beta(dz).
\end{align*}
Furthermore we have for each $p>1$
\begin{equation}\label{decay}
\sup_{r\in ]0,1/2]}\sup_{x\in\R^2}(-\ln r)^{p}M_\beta(B(x,r)\cap K_{\delta,L})\leq 2^pL.
\end{equation}
Indeed, observe first that $\sup_{r\in ]0,1/2]}\sup_{x\in K_{\delta,L}}(-\ln r)^{p}M_\beta(B(x,r)\cap K_{\delta,L})\leq  L$ by definition. To extend this formula to $x\not \in K_{\delta,L}$, take $r\leq 1/2$ and observe that we have two options: either $B(x,r)\cap K_{\delta,L}$ is empty in which case $M_\beta(B(x,r)\cap K_{\delta,L})=0$ or   we can find $y\in B(x,r)\cap K_{\delta,L}$ and then $M_\beta(B(x,r)\cap K_{\delta,L})\leq M_\beta(B(y,2r)\cap K_{\delta,L})\leq L(-\ln 2r)^{-p}\leq 2^pL(-\ln r)^{-p}$.

We deduce that for all $R>0$,
\begin{equation}\label{logzerobis}
\sup_{y\in B(0,R)}\int_{\R^2}\int_0^tp_r(y,z)\,dr\ind_{K_{\delta,L}}(z)M_\beta(dz)\to 0,\quad \text{ as }t\to 0.
\end{equation}
Indeed, the decay of the size of balls induced by \eqref{decay}, used with $p>2$, is enough to overcome the $\ln$-singularity produced by the heat kernel integral: $\int_0^tp_r(y,z)\,dr$ (to be exhaustive, one should adapt the argument in \cite[section 2.7]{GRV1} but this is harmless).  Hence \eqref{logzero}.

Let us now prove that  the family $(F^{',\delta,L}_{\beta}(x_n,\cdot))_n$ is tight in $C(\R_+;\R_+)$. This consists in checking that for all $T,\eta>0$:
\begin{equation}\label{tight1}
\lim_{\delta\to 0}\limsup_{n\to \infty}\Pb^{B}\big(\sup_{\substack{0\leq s,t\leq T\\|t-s|\leq \delta}}|F^{',\delta,L}_{\beta}(x_n,]s,t])|>\eta\big)=0.
\end{equation}
The control of this supremum  uses two arguments: a control of the involved quantities when $s,t$ are closed to $0$ via \eqref{logzero} and a control of this supremum via a coupling argument when $s,t$ are far enough from $0$. To quantify the proximity to $0$, we introduce a parameter $\theta>0$. We have for $\delta<\theta$:
\begin{align}\label{tard1}
\Pb^{B}\big(&\sup_{\substack{0\leq s,t\leq T\\|t-s|\leq \delta}}|F^{',\delta,L}_{\beta}(x_n,]s,t])|>\eta\big)\nonumber\\
\leq &\Pb^{B}\big(  F^{',\delta,L}_{\beta}(x_n,2\theta)>\eta/4\big)+\Pb^{B}\big( \sup_{\substack{\theta\leq s,t\leq T\\|t-s|\leq \delta}}|F^{',\delta,L}_{\beta}(x_n,]s,t])|>\eta/2\big).
\end{align}
To establish \eqref{tight1}, it is thus enough to prove that the $\limsup_{\theta\to 0}\limsup_{\delta\to 0}\limsup_{n\to \infty}$ of each term  in the  right-hand side of the above expression vanishes.
The first term  is easily   treated with the help of the Markov inequality and \eqref{logzero} so that we now focus on the second term. To this purpose, we recall the  following coupling lemma:

\begin{lemma}\label{lem:coupling}
Fix $x\in\R^2$ and let us start a Brownian motion $B^{x}$ from $x$. Let us consider another independent Brownian motion $B'$ starting from $0$ and denote by $B^y$, for a rational $y\in\Q^2$, the Brownian motion $B^y=y+B'$. Let us denote by $\tau^{x,y}_1$ the first time at which the first components of $B^{x}$ and $B^y$ coincide:
$$\tau_1^{x,y}=\inf\{u>0;B^{1,x}_u=B^{1,y}_u\}$$ 
and by $\tau_2^{x,y}$ the first time at which the second components  coincide after $\tau_1^{x,y}$:
$$\tau_2^{x,y}=\inf\{u>\tau_1^{x,y};B^{2,x}_u=B^{2,y}_u\}$$ 
The random process $\overline{B}^{x,y}$ defined by 
$$\overline{B}^{x,y}_t=\left\{ 
\begin{array}{lll}
(B^{1,x}_t,B^{2,x}_t)  & \text{if}  &  t\leq \tau_1^{x,y} \\
(B^{1,y}_{t}, B^{2,x}_t)  & \text{if}  &  \tau_1^{x,y}<t\leq \tau_2^{x,y}\\
(B^{1,y}_{t}, B^{2,y}_t)  & \text{if}  &  \tau_2^{x,y}<t .
\end{array}
\right.$$ is  a new Brownian motion on $\R^2$ starting from $x$, and coincides with $B^y$ for all times $t>\tau_2^y$. Furthermore, if $y-x\to 0$, we have for all $\eta>0$:
$$\Pb(\tau_2^{x,y}>\eta)\to 0 .$$
\end{lemma}

We choose $y\in \Q^2$. We can consider the couple  $(F^{',\delta,L,x_n}_{\beta}(y,\cdot),F^{',\delta,L}_{\beta}(x_n,\cdot))$ where  $F^{',\delta,L}_{\beta}(x_n,\cdot)$ is the same as that considered throughout this section and  $F^{',\delta,L,x_n}_{\beta}(y,\cdot)$ is constructed as $F^{',\delta,L}_{\beta}(y,\cdot)$ but we have used  the Brownian motion $\bar{B}^{y,x_n}$ of Lemma \ref{lem:coupling} instead of the Brownian motion $B^y$. The important point to understand is that this couple does not have the same law as the couple
$(F^{',\delta,L}_{\beta}(y,\cdot),F^{',\delta,L}_{\beta}(x_n,\cdot))$ but it has the same $1$-marginal. Furthermore we have $F^{',\delta,L,x_n}_{\beta}(y,]s,t])=F^{',\delta,L}_{\beta}(x_n,]s,t])$ for $\tau_2^{y,x_n}\leq s<t$. We deduce:
\begin{align*}
\Pb^{B}\big(  \sup_{\substack{\theta\leq s,t\leq T\\|t-s|\leq \delta}} &|F^{',\delta,L}_{\beta}(x_n,]s,t])|>\eta/2\big)\\
\leq & \Pb^{B}\big(  \sup_{\substack{\theta\leq s,t\leq T\\|t-s|\leq \delta}}|F^{',\delta,L}_{\beta}(y,]s,t])|>\eta/2\big)+\Pb( \tau_2^{y,x_n}>\theta).
\end{align*}
It is then obvious to get:
\begin{equation}\label{tard2}
\limsup_{\delta\to 0}\limsup_{n\to \infty}\Pb^{B}\big(  \sup_{\substack{\theta\leq s,t\leq T\\|t-s|\leq \delta}}  |F^{',\delta,L}_{\beta}(x_n,]s,t])|>\eta/2\big)\leq \Pb( \tau_2^{y,x}>\theta).
\end{equation}
Since the choice of $y$ was arbitrary, we can now choose $y$ arbitrarily close to $x$ to make this latter term as close to $0$ as we please for a fixed $\theta$. Hence \eqref{tight1} and the family $ (F^{',\delta,L}_{\beta}(x_n,\cdot))_n$ is tight in $C(\R_+,\R_+)$.

We can also use the coupling argument to prove that there is only one possible limit in law for all subsequences $(x_n)_n$ such that $x_n\to x$ as $n\to\infty$, thus showing the convergence in law in $C(\R_+,\R_+)$ of the family $ (F^{',\delta,L}_{\beta}(x_n,\cdot))_n$. Here we have only dealt with the convergence of the family $ (F^{',\delta,L}_{\beta}(x_n,\cdot))_n$ but it is straightforward to adapt the argument to the family $ (B^{x_n},F^{',\delta,L}_{\beta}(x_n,\cdot))_n$. 

Now, we come to the convergence in $\Pb^B$-probability of $ (F^{',\delta,L}_{\beta}(x_n,\cdot))_n$. We fix $t>0$, we consider the mapping
$$(x,y)\in\Q^2\times\Q^2\mapsto \E^{B}\big[ (F^{',\delta,L}_{\beta}(x,t)-F^{',\delta,L}_{\beta}(y,t))^2\big].$$
We can expand the square and, following  the ideas in \cite[section 2.7]{GRV1}, use \eqref{decay} to control the $\ln$-singularities to see that the above mapping extends to a continuous function on $\R^2\times \R^2$, which vanishes on the diagonal $\{(x,x);x\in\R^2\}$.
%\begin{align*}
%2\int_{\R^2\times\R^2}\int_0^t\int_s^t p_s(x,u)p_{r-s}(y+u,z+x)\,dr\,dsM_{\epsilon'}^\beta(du)M_{\epsilon}^\beta(dz)
%\end{align*}
Therefore, if $(x_n)_n$ is a sequence in $\Q^2$ converging towards $x$, the sequence $(F^{',\delta,L}_{\beta}(x_n,t))_n$ converges in $L^2$ under $\Pb^B$. We can thus extract a subsequence $(x_{\phi(n)})_n$ such that for all $t\in\Q\cap\R_+$, the sequence $(F^{',\delta,L}_{\beta}(x_{\phi(n)},t))_n$ converges   $\Pb^B$-almost surely. Also, we have seen that this subsequence converges in law in $C(\R_+,\R_+)$ towards a continuous random mapping. The Dini theorem implies that $(F^{',\delta,L}_{\beta}(x_{\phi(n)},\cdot))_n$ converges   $\Pb^B$-almost surely  in $C(\R_+,\R_+)$ towards a limit denoted by $F^{',\delta,L}_{\beta}(x,\cdot)$.\qed

We are now in position to construct $F'_\beta$ on the whole of $\R^2$.
\begin{theorem}\label{fprimebeta}
Almost surely in $X$, for all $x\in\R^2$, the random measure $F^{',\delta,L}_{\beta}(x,dt)$ converges $\Pb^B$-almost surely as $\delta\to 0$ and $L\to\infty$ in the sense of weak convergence of measures towards a random measure denoted by $F'_\beta(x,dt)$. Furthermore:
\begin{enumerate}
\item for $x\in\Q^2$, $F'_\beta(x,\cdot)$ coincides with  the limit of the family $(F^{',\epsilon}_\beta(x,\cdot))_\epsilon$ as $\epsilon\to 0$  defined in subsection \ref{ss.deriv}.
\item for $x\in S\cup\Q^2$, convergence of the mapping $t\mapsto F^{',\delta,L}_{\beta}(x,t)$ towards $t\mapsto F^{'}_{\beta}(x,t)$ holds $\Pb^B$-almost surely in $C(\R_+,\R_+)$ as $\delta\to 0$, $L\to\infty$.
\item for $x\not \in S$, for all $s>0$, convergence of the mapping $t\mapsto F^{',\delta,L}_{\beta}(x,]s,t])$ towards $t\mapsto F^{'}_{\beta}(x,]s,t])$ holds $\Pb^B$-almost surely in $C([s,+\infty[,\R_+)$.
\end{enumerate}
\end{theorem}

\noindent {\it Proof.} Observe that $K_{\delta,L}\subset K_{\delta',L'}$ if $\delta'\leq \delta$ and $L'\geq L$. Therefore, for all $x\in \R^2$ and for all $t\geq 0$, $F^{',\delta,L}_{\beta}(x,t)\leq F^{',\delta',L'}_{\beta}(x,t)$ if $\delta'\leq \delta$. We can thus define the almost sure limit:
$$F^{',\infty}_\beta(x,t)=\lim_{\delta\to 0,L\to\infty}F^{',\delta,L}_{\beta}(x,t).$$
Actually, this also implies the weak convergence as $\delta\to 0$ and $L\to\infty$ of the measure $F^{',\delta,L}_{\beta}(x,dt)$ towards a random measure, still denoted by $F^{',\infty}_\beta(x,dt)$.

For $x\in \Q^2$, let us identify $F^{',\infty}_\beta(x,\cdot)$ with the limit $F^{'}_\beta(x,\cdot)$ of the family $(F^{',\epsilon}_\beta(x,\cdot))_\epsilon$ as $\epsilon\to 0$. By construction  we have:
$$F^{',\delta,L}_\beta(x,t)=\int_0^t\ind_{K_{\delta,L}}(B^x_r)F^{'}_\beta(x,dr)\leq F^{'}_\beta(x,t),$$
in such a way that $F^{',\infty}_\beta(x,t)\leq F^{'}_\beta(x,t)$. Second, by the dominated convergence theorem we get
\begin{align*}
\E^{B}\big[|F^{'}_\beta(x,t)-F^{',\infty}_\beta(x,t)|\big]=&\lim_{\delta\to0,L\to\infty}\E^{B }\big[\int_0^t\ind_{K_{\delta,L}^c}(B^x_r)F^{'}_\beta(x,dr)\big]\\
=&\lim_{\delta\to0,L\to\infty}\int_{\R^2}\int_0^tp_r(x,y)\,dr\ind_{K_{\delta,L}^c}(y)M_\beta(dy) \\
=&\int_{\R^2}\int_0^tp_r(x,y)\,dr\ind_{S^c}(y)M_\beta(dy) \\
=&0 
\end{align*}
because $M_\beta(S^c)=0$. Now that we have identified $F^{'}_\beta(x,\cdot)$ with $F^{',\infty}_\beta(x,\cdot)$ on $x\in\Q^2$, we skip the distinction made with the superscript $\infty$ and write  $F^{'}_\beta$ for the limit of $F^{',\delta,L}_\beta$.

For $x\in \Q^2$, the continuity of the mapping $t\mapsto F^{'}_\beta(x,t)$ together with the Dini theorem implies the $\Pb^B$-almost sure convergence of $F^{',\delta,L}_{\beta}(x,\cdot)$ towards $ F^{'}_{\beta}(x,\cdot)$   in $C(\R_+,\R_+)$. Let us now complete the proof of items 2 and 3. We fix $s>0$. The coupling argument established in the proof of Proposition \ref{propadapted} shows that the mapping $t\mapsto F^{'}_{\beta}(x,]s,t]) $ is continuous on $[s,+\infty[$ (it coincides in law with $t\mapsto F^{'}_{\beta}(y,]s,t])$ with $y\in\Q^2$ as soon as $B^x$ and $B^y$ are coupled before time $s$, which happens with probability arbitrarily close to $1$ provided that $y$ is close enough to $x$). The Dini theorem again implies item 3. Let us stress that it is not clear that we can take $s=0$ because we need to control the decay of balls at $x\in S^c$ and this decay may happen to be very bad on $S^c$. 

To prove item 2, we also use  the Dini theorem but we further need to prove that the mapping $t\mapsto F^{'}_{\beta}(x,t) $ is continuous at $t=0$ with $F^{'}_{\beta}(x,0) =0$. This can be done by computing
\begin{align*}
\E^{B}\big[F^{'}_\beta(x,t)\big]=&\int_{\R^2}\int_0^tp_r(x,y)\,dr\ind_{S}(y)M_\beta(dy).
\end{align*}
Following \cite{GRV1},  we observe that the mapping $y\mapsto \int_0^tp_r(x,y)\,dr$ possesses a logarithmic singularity at $y=x$. Furthermore, for $x\in S$, we have 
$$\sup_{\substack{x\in \mathcal{O}\text{ open }\subset B,\\ {\rm diam}(\mathcal{O})\leq 1}}M'(  \mathcal{O})\big(-\ln{\rm diam}(\mathcal{O})\big)^{p}<+\infty.$$
Therefore
$$\int_{\R^2}\int_0^1p_r(x,y)\,dr\ind_{S}(y)M_\beta(dy)<+\infty$$ in such a way that the dominated convergence theorem implies
\begin{align*}
\E^{B}\big[F^{'}_\beta(x,t)\big]=\int_{\R^2}\int_0^tp_r(x,y)\,dr\ind_{S}(y)M_\beta(dy)\to 0,\quad \text{ as }t\to 0.
\end{align*}
To sum up, we have proved that the mapping $t\mapsto F^{'}_{\beta}(x,t) $ is continuous at $t=0$ with $F^{'}_{\beta}(x,0) =0$ for $x\in S$. We complete the proof of item 2 with the help of the Dini theorem.\qed
 
\begin{remark}\label{wrongPCAF}
It is important here to stress that the above Proposition shows that for all $x$, we have defined a  mapping $F^{'}_\beta(x,\cdot)$, which is $\Pb^B$-almost surely continuous with continuous sample paths and satisfies a variant of the additivity property of a PCAF, i.e. for all $s,t\geq 0$ we have almost surely
\begin{equation*}
F^{'}_\beta(x,t+s)= F^{'}_\beta(x,t)+\bar{F}^{'}_\beta(x,s),
\end{equation*}
where, conditionally to $\mathcal{F}_t$, the variable $\bar{F}^{'}_\beta(x,s)$ is distributed as $F^{'}_\beta(x+B_t,s)$ under $\Pb^{x+B_t}$ (measure of a Brownian motion starting from $x+B_t$). Also, we have not so far proved that $F^{'}_\beta$ is a PCAF because it is  defined for all $x$ $\Pb^B$-almost surely whereas we need to define it $\Pb^B$-almost surely for all $x$. Yet, we will see that this problem for the construction of a proper PCAF is not too serious.
\end{remark} 
We claim:

\begin{theorem}\label{PCAF}
Almost surely in $X$, we define
$$F'(x,t)=\lim_{\beta\to\infty}F'_\beta(x,]0,t])$$
where convergence holds $\Pb^B$-almost surely in $C(\R_+,\R_+)$. $F'$ is some form of PCAF in the strict sense in $\R^2$: it is defined for all starting points and satisfies the following variant of the additivity property
\begin{equation}\label{weakadditivity}
F^{'}_\beta(x,t+s)= F^{'}_\beta(x,t)+\bar{F}^{'}_\beta(x,s), \quad s,t \geq 0
\end{equation}
where, conditionally to $\mathcal{F}_t$, the variable $\bar{F}^{'}_\beta(x,s)$ is distributed as $F^{'}_\beta(x+B_t,s)$ under $\Pb^{x+B_t}$ (measure of a Brownian motion starting from $x+B_t$).

 Furthermore, 
\begin{enumerate}
\item for all $x\in\R^2$, $F'$ is continuous, strictly increasing and goes to $\infty$ as $t\to\infty$.
\item $F'$ coincides outside a set of zero capacity with a PCAF of Revuz measure $M'$.
\end{enumerate}
\end{theorem}

\noindent {\it Proof.} For each $x\in\R^2$, for each $t\geq 0$, the mapping  $\beta\mapsto F'_\beta(x,]0,t])$ is increasing. We can thus define  $F'(x,t)=\lim_{\beta\to\infty}F'_\beta(x,]0,t])$. Furthermore, for each ball $B$ containing $x$, $F'(x,t)$ coincides with 
$F'_\beta(x,]0,t])$ for $t<\tau_B(x)=\inf\{u>0; B^x_u\not\in B\}$ and $\beta$ (random) large enough (more precisely for $\beta$ large enough to make $\sup_{x\in B}\sup_{\epsilon\in]0,1]} X_\epsilon(x)-2\ln\frac{1}{\epsilon}<\beta$, see Proposition \ref{coro:max}). It is obvious to check that $F'$ satisfies the additivity \eqref{weakadditivity}.

%satisfies the properties of a PCAF excepted the content of Remark \ref{wrongPCAF}, namely that the defining set of the PCAF %depends on the point $x$. We will come back to this point below.

Now we prove item 1. This results from the coupling argument detailed in Proposition \ref{propadapted} as $F'(x,\cdot)$ is strictly increasing and goes to $\infty$ as $t\to \infty$ for $x\in\Q^2$ (see also \cite[Proposition 2.24]{GRV1}).

Finally, we prove item 2, more precisely we establish the relation \eqref{revuzcorr} for $M'$ and $F'$. The construction of $F'$ entails that, for any $0<s<t$
$$
\E^{B^x}\Big[\int_s^tf(B^x_r)\,F'(x,dr)\Big]=\int_s^t\int_{\R^2}f(y)p_r(x,y)\,M'(dy)\,dr.
$$
Therefore, for any nonnegative Borel functions $f,h$ ($P_r$ stands for the semigroup associated to the planar Brownian motion):
\begin{align*}
\int_{\R^2}h(x)\E^{B^x}\Big[\int_s^tf(B^x_r)\,F'(x,dr)\Big]\,dx=&\int_s^t\int_{\R^2}\int_{\R^2}f(y)p_r(x,y)\,M'(dy)h(x)\,dx\,dr\\
=&\int_s^t \int_{\R^2}f(y)P_rh(y)\,M'(dy)\,dr.
\end{align*}           
It suffices to let $s\to 0$  to conclude.

Now, we use \cite[Theorem 5.1.3]{fuku} to prove that there exists a PCAF $A$ associated to $M'$ because $M'$ is a smooth measure in the sense of \cite{fuku} thanks to Corollary \ref{th:capsup}.  Now we have at our disposal a PCAF $ A$ with Revuz measure $M'$ and an "almost" PCAF $F'$. 
%("almost" means: up to the objection made in  Remark \ref{wrongPCAF}) satisfying  \eqref{revuzcorr}. 
The reader may check that the uniqueness part of \cite[Theorem 5.1.4]{fuku} can be reproduced to prove that $F'$ and $A$ coincide for $x$ outside a set of capacity $0$  (just observe that this proof does not use the fact that the set where the PCAF is defined does not depend on $x$).
\qed

\subsection{Definition and properties of the critical LBM}
%%%%%%%%%%%%%%%%%%%%%%%%%%%%%%%%%

\begin{definition}{\bf (critical Liouville Brownian motion).}\label{LBMdefnew}
Almost surely in $X$, for all $x\in \R^2$, the  law of the  LBM  at criticality, starting from $x$, is defined by:
$$\mathcal{B}^x_t=x+B_{\langle \mathcal{B}^x\rangle_t}$$ where $\langle \mathcal{B}^x\rangle$ is defined by
$$\sqrt{2/\pi} F'(x,\langle \mathcal{B}^x\rangle_t)=t.$$
\end{definition}
We stress that $\LB^x$ is a local martingale. 

\begin{proposition}
The critical LBM is a strong Markov process with continuous sample paths.
\end{proposition}

\noindent {\it Proof.} Strong Markov property results from  \cite[sect. 6]{fuku}. Continuity of sample paths results from the fact that $F'$ is strictly increasing.\qed

\begin{theorem}\label{th.LBMy}
Assume further [A.6] or [A.6']. Almost surely in $X$, for all $x\in S$, the $\epsilon$-regularized Brownian motion $( \mathcal{B}^{\epsilon,x})_\epsilon$ defined by Definition \ref{d.elbm2} converges in law in the space $C(\R_+,\R^2)$ equipped with the supremum norm on compact sets  towards $\mathcal{B}^x$.
\end{theorem}

\vspace{1mm}
\noindent {\it Proof.} This is just a consequence of Theorems \ref{fprimebeta} and \ref{PCAF} as explained in \cite{GRV1}. \qed
 
%\subsubsection{Dirichlet form}
%%%%%%%%%%%%%%%%%%%%%%%%%%

\vspace{1mm}
From \cite[Th. 6.2.1]{fuku}, we claim:
 
\begin{theorem}{\bf (Dirichlet form).}\label{dirichlet}
The critical Liouville Dirichlet form $(\Sigma,\mathcal{F})$ takes on the following explicit form on $L^2(\R^2,M')$:
\begin{equation}\label{formesimple}
\Sigma(f,g)=\frac{1}{2}\int_{\R^2}\nabla f(x)\cdot \nabla g(x)\,dx
\end{equation} 
with domain
$$\mathcal{F}=\Big\{f\in L^2(\R^2,M')\cap H^1_{loc}(\R^2,dx); \nabla f\in L^2(\R^2,dx)\Big\},$$
Furthermore, it is strongly local and regular.
\end{theorem}

%\subsubsection{Semi-group} 
%%%%%%%%%%%%%%%%%%%%%%%%%% 

Let us denote by $P^X_t$ (for $t\geq 0$) the mapping
\begin{equation}\label{semigroup}
f\in C_b(\R^2)\mapsto \big(x\in \R^2\mapsto P^X_tf(x)=\E^B[f(\LB^{x}_t)]\big).
\end{equation}
Similarly we define $P^\epsilon$ as the semigroup generated by the Markov process $\LB^{\epsilon}$. From \cite[sect. 6]{fuku}, we claim:
\begin{theorem}{\bf (Semigroup).}\label{th:semi}
The linear operator $P^X_t$, restricted to $C_c(\R^2)$, extends to a linear contraction on $L^p(\R^2, M')$ for all $1\leq p <\infty$, still denoted $P^X_t$. Furthermore:
\begin{itemize}
\item $(P^X)_{t\geq 0}$ is a Markovian strongly continuous semigroup on $L^p(\R^2, M')$ for $1\leq p <\infty$.
\item Assume further [A.6] or [A.6']. Almost surely in $X$, the $\epsilon$-regularized semigroup $(P^\epsilon)_\epsilon$ converges pointwise for $x\in S$ towards the critical Liouville semigroup. More precisely, for all bounded continuous function $f$, we have: 
$$\forall x\in S,\forall t\geq 0, \quad \lim_{\epsilon\to 0}P^\epsilon_t f(x)=P^X_t f(x).$$
\item $P^X$ is self-adjoint in $L^2(\R^2, M')$.
\item the measure $M'$ is invariant for $P^X_t$.
\end{itemize}
\end{theorem}

%\subsubsection{Critical Liouville Laplacian}
%%%%%%%%%%%%%%%%%%%%%%%%
The {\bf critical Liouville Laplacian} $\Delta_X$ is defined as the generator of the critical Liouville semigroup times the usual extra factor $2$. The critical Liouville Laplacian  corresponds to an operator which can formally be written as
$$ \Delta_X=X^{-1}(x)e^{-2 X (x)}\Delta$$
and can be thought of as the Laplace-Beltrami operator of $2d$-Liouville quantum gravity at criticality (of course when $X$ is a free field).

%\subsection{Analysis of the critical Liouville resolvent and consequences}
%%%%%%%%%%%%%%%%%%%%%%%%%%
One may also consider the resolvent family $(R^X_\lambda)_{\lambda>0}$ associated to the semigroup $(P^X_\t)_\t$. In a standard way, the resolvent operator  reads:
\begin{equation}\label{def:resol}
\forall f\in C_b(\R^2),\forall x\in \R^2,\quad R^X_\lambda f(x)=\int_0^{\infty}e^{-\lambda t}P_t^Xf(x)\,dt.
\end{equation}
Furthermore, the resolvent family $(R^X_\lambda)_{\lambda>0}$ extends to $L^p(\R^2,M')$ for $1\leq p<+\infty$, is strongly continuous for $1\leq p<+\infty$ and is self-adjoint in $L^2(\R^2,M')$.  This results from the properties of the semi-group. As a consequence of Theorem \ref{th:semi}, it is straightforward to see that:
\begin{proposition}\label{cv:resolvent}
Assume further [A.6] or [A.6']. Almost surely in $X$, the $\epsilon$-regularized resolvent family $(R^\epsilon_\lambda)_\lambda$  converges towards
the critical Liouville resolvent $(R^X_\lambda)_\lambda$ in the sense that for all function $f\in C_b(\R^2)$:
$$\forall x \in S,\quad \lim_{\epsilon\to 0}R^\epsilon_\lambda f(x)=R^X_\lambda f(x).$$
\end{proposition}
Also and similarly to \cite{fuku,GRV2}, it is possible to get an explicit expression for the resolvent operator:
\begin{proposition}\label{expr:resol1}
Almost surely in $X$, the resolvent operator takes on the following form for all measurable bounded function $f$ on $\R^2$: 
 $$\forall x\in \R^2,\quad R_\lambda^Xf(x)=\sqrt{2/\pi}\E^B\big[\int_0^{\infty}e^{-\lambda \sqrt{2/\pi}F'(x,t)}f(B^x_t)\,F'(x,dt)\big].$$
\end{proposition}
 
The main purpose of this section is to prove the following  structure result on the resolvent family:
\begin{theorem}{\bf (massive Liouville  Green kernels at criticality)}.\label{th:green}
For every $x\in\R^2$, the resolvent family $(R^X_\lambda)_{\lambda>0}$ is absolutely continuous with respect to the critical Liouville measure. Therefore there exists a family $(\r^X_\lambda(\cdot,\cdot))_\lambda$, called the family of massive critical Liouville Green kernels, of jointly measurable functions such that:
$$\forall x\in \R^2,\forall f \in B_b(\R^2),\quad R^X_\lambda f(x)=\sqrt{2/\pi}\int_{\R^2}f(y)\r^X_\lambda(x,y)\,M'(dy)$$ and such that:\\
1) (strict-positivity) for all $\lambda>0$ and $x\in \R^2$,  $M'(dy)$ a.s., $\r^X_\lambda(x,y)>0,$ \\
2) (symmetry) for all $\lambda>0$ and $x,y\in \R^2$: $\r^X_\lambda(x,y)=\r^X_\lambda(y,x),$\\
3) (resolvent identity) for all $\lambda,\mu>0$, for all $x,y \in \R^2$, $$\r^X_\mu(x,y)-\r^X_\lambda(x,y)=(\lambda-\mu)\sqrt{2/\pi}\int_{\R^2}\r^X_\lambda(x,z)\r^X_{\mu}(z,y)\,M'(dz).$$  
4) ($\lambda$-excessive) for every $y$: $e^{-\lambda t}P_t^X(\r_\lambda(\cdot,y))(x)\leq \r_\lambda(x,y)$ for $M'$-almost every $x$ and for all $t>0$.
\end{theorem}
 
\vspace{1mm}
\noindent {\it Proof.}  We have to show absolute continuity of the resolvent for $x\in \R^2$. Though inspired by \cite{GRV1,GRV2}, we have to adapt the proof because we do not have "uniform convergence" of the PCAF towards $0$ as $t\to 0$. In particular, it is not clear that the resolvent be strong Feller.
 For $\delta>0$, we define for $f\in B_b(\R^2)$
$$R_{\lambda,\delta}^Xf(x)= \sqrt{2/\pi}\E^B\big[\int_\delta^{\infty}e^{-\lambda \sqrt{2/\pi}F'(x,]\delta,t])}f(B^x_t)\,F'(x,dt)\big],$$
where $F'(x,dr)$ stands for the random measure associated to the increasing function $t\mapsto F'(x,t)$.
Once again, the coupling argument of Proposition   \ref{propadapted}, it is plain to see that the mapping 
$$x\in\R^2\mapsto R_{\lambda,\delta}^Xf(x)$$ is continuous. Now we claim 
\begin{lemma}\label{mass1}
For every $x\in \R^2$, $\delta>0$ and all nonnegative bounded Borelian function $f$, we have
$$R_{\lambda}^Xf(x)=0 \Longrightarrow R_{\lambda,\delta}^Xf(x)=0.$$
\end{lemma}
\begin{lemma}\label{mass2}
For every $x\in \R^2$ and all nonnegative bounded Borelian function $f$, we have
$$\forall \delta> 0,\quad R_{\lambda,\delta}^Xf(x)=0\Longrightarrow R_{\lambda}^Xf(x)=0.$$
\end{lemma}

We postpone the proofs of the  above two lemmas. If $A$ is a measurable set such that $M'(A)=0$ then by invariance of $M'$ for the resolvent family, we deduce that $R_\lambda\ind_A(x)=0$ for $M'$-almost every $x\in\R^2$. Since $S$ has full $M'$-measure and because $M'$ has full support, we deduce that   $R_{\lambda}^Xf(x)=0$ for $x$ belonging to a dense subset of $\R^2$. From Lemma \ref{mass1}, $R_{\lambda,\delta}^Xf(x)=0$ for $x$ belonging to a dense subset of $\R^2$. Continuity of $R_{\lambda,\delta}^Xf$ entails that this function identically vanishes on $\R^2$ for every $\delta>0$. With the help of Lemma \ref{mass2}, we deduce   that  $R_\lambda\ind_A(x)=0$ for $x\in \R^2$. Therefore, for all $x\in \R^2$, $R_\lambda \ind_A(x)=0$, thus showing that the measure $A\mapsto R_{\lambda}^X \ind_A(x)$ is absolutely continuous with respect to $M'$.\qed

\vspace{2mm}
\noindent {\it Proof of Lemmas \ref{mass1} and \ref{mass2}.} For $x\in \R^2$ and all bounded nonnegative Borelian function $f$, we have (see Proposition \ref{expr:resol1})
\begin{align}
R_\lambda^Xf(x)=&\sqrt{2/\pi}\E^{B^x}\big[\int_0^{\infty}e^{-\lambda \sqrt{2/\pi}F'(x,t)}f(B^x_t)\,F'(x,dt)\big]\nonumber\\
=&\sqrt{2/\pi}\E^{B^x}\big[\int_0^{\delta}e^{-\lambda \sqrt{2/\pi}F'(x,t)}f(B^x_t)\,F'(x,dt)\big]\nonumber\\
  &+\sqrt{2/\pi}\E^{B^x}\big[e^{-\lambda F'(x,\delta)}\int_\delta^{\infty}e^{-\lambda \sqrt{2/\pi}F'(x,]\delta,t])}f(B^x_t)\,F'(x,dt)\big].\label{majresol}
\end{align}
Therefore
$$R_\lambda^Xf(x)\geq  \sqrt{2/\pi}\E^{B^x}\big[e^{-\lambda F'(x,\delta)}\int_\delta^{\infty}e^{-\lambda \sqrt{2/\pi}F'(x,]\delta,t])}f(B^x_t)\,F'(x,dt)\big].$$
For $x\in \R^2$, we have $e^{-\lambda F'(x,\delta)}>0$ $\Pb^{B^x}$-almost surely. The proof of Lemma \ref{mass1} follows.

Furthermore, for $x\in \R^2$, we have  
\begin{align*}
\sqrt{2/\pi}\E^{B^x}\big[\int_0^{\delta}e^{-\lambda \sqrt{2/\pi}F'(x,t)}&f(B^x_t)\,F'(x,dt)\big]\\
\leq &\lambda^{-1}\|f\|_\infty\E^{B^x}\big[1-e^{-\lambda \sqrt{2/\pi}F'(x,\delta)} \big].
\end{align*}
Since $F'(x,\delta)$ converges in law towards $0$ as $\delta \to 0$, we deduce that
$$\lim_{\delta\to 0}\sqrt{2/\pi}\E^{B^x}\big[\int_0^{\delta}e^{-\lambda \sqrt{2/\pi}F'(x,t)} f(B^x_t)\,F'(x,dt)\big]=0.$$
From \eqref{majresol}, we deduce:
$$R_\lambda^Xf(x)\leq R_{\lambda,\delta}^Xf(x)+\sqrt{2/\pi}\E^{B^x}\big[\int_0^{\delta}e^{-\lambda \sqrt{2/\pi}F'(x,t)} f(B^x_t)\,F'(x,dt)\big].$$
The proof of Lemma \ref{mass2} follows.\qed

 \vspace{2mm}

%\subsubsection{Liouville Green function at criticality}
%%%%%%%%%%%%%%%%%%%%%%%%%%

As prescribed in \cite[section 1.5]{fuku}, let us define the Green function for $f\in L^1(D,M')$  by
$$ Gf(x)=\lim_{t\to\infty} \int_0^tP^X_rf(x)\,dr.$$
We further denote $g$ the standard Green kernel on $\R^2$.  Following \cite{fuku}, we say that the semi-group $(P^X_t)_t$, which is symmetric w.r.t. the measure $M'$,   is   {\it irreducible} if  any $P^X_t$-invariant set $B$ satisfies $M'(B)=0$ or $M'(B^c)=0$. We say that $(P_t^X)$ is recurrent if, for any $f\in L^1_+(D,M')$, we have $Gf(x)=0$ or $Gf(x)=+\infty$ $M'$-almost surely.

\begin{theorem}{\bf (Liouville Green function at criticality).} The critical Liouville semi-group is irreducible and recurrent.The critical Liouville Green function, denoted by $G^X$, is given for every $x\in S$ by
$$G^Xf(x)=\sqrt{2/\pi}\int_{\R^2}\frac{1}{\pi}\ln\frac{1}{|x-y|}f(y)\,M'(dy)$$ for all functions $f\in L^1(\R^2,M')$ such that
$$\int_{\R^2}f(y)\,M'(dy)=0.$$
\end{theorem}

\vspace{1mm}
\noindent {\it Proof.} Irreducibility is a straightforward consequence of Theorem \ref{th:green} and the remaining part of the statement is a straightforward adaptation of \cite{GRV1,GRV2} for $x\in S$.  \qed

\vspace{2mm} 
%\subsubsection{Critical Liouville heat kernel}
%%%%%%%%%%%%%%%%%%%%%%%%%%
We investigate now the existence of probability densities of the critical Liouville semi-group with respect to the critical Liouville measure.
\begin{theorem}{\bf (Critical Liouville heat kernel).}\label{th:heat}The critical Liouville semigroup $(P^X_\t)_{\t>0}$ is absolutely continuous with respect to the critical Liouville measure. There exists a family of nonnegative   functions $(\p^X_t(\cdot,\cdot))_{t\geq 0}$, which we call the critical Liouville heat kernel,  such that:
$$\forall x\in\R^2,dt \,\text{a.s.},\forall f \in B_b(\R^2),\quad P^X_t f(x)=\sqrt{2/\pi}\int_{\R^2}f(y)\p^X_t(x,y)\,M'(dy).$$ 
\end{theorem}

\vspace{1mm}
\noindent {\it Proof.} From Theorems \ref{th:green} and \cite[Theorems 4.1.2 and 4.2.4]{fuku}, the Liouville semi-group is absolutely continuous with respect to the Liouville measure.  \qed

\begin{remark}
Though we call the family $(\p^X_t(\cdot,\cdot))_{t\geq 0}$ heat kernel, we are not in position to establish most of the regularity properties expected from a heat kernel. Furthermore,  a weak form of the notion of spectral dimension is obtained in \cite{spectral}, which is $2$. We do not know how to adapt the argument in the critical case because  we cannot prove the continuity of the mapping $(x,y)\mapsto \p(t,x,y)\,dt$.
\end{remark}

%\subsubsection{Critical LBM and thick points of $X$}
%%%%%%%%%%%%%%%%%%%%%%%%%%%%%%%%%%%%%%%%%%%%%%%%%
Let us consider the set $E$ defined in Theorem \ref{th:cap}. Recall that $M'(E^c)=0$. As a consequence of Theorem \ref{th:green}, we obtain the following result where $\lambda$ is the Lebesgue measure: 
\begin{corollary}
Almost surely in $X$, for all starting points $x\in\R^2$, the critical LBM spends Lebesgue-almost all the time in the set $E$: 
 \begin{equation*}
\text{a.s. in }X, \forall x \in\R^2, \text{ a.s. under }\Pb^{B^x},\quad \lambda \lbrace t \geq 0; \:  \LB^x_t \in    E^c \rbrace =0.
\end{equation*}
\end{corollary}
 
If one applies Theorem \ref{th:heat} instead, one obtains the similar but different result:

\begin{corollary}
Almost surely in $X$, for all $t>0$
$$\Pb^{B^x} \text{a.s.},\quad \LB^x_t\in E.$$
\end{corollary}

\section{Further remarks and GFF on other domains}\label{other}
%%%%%%%%%%%%%%%%

So far, we constructed in detail the LBM on (subdomains of) $\R^2$. This construction may be adapted to other geometries like  the sphere $\S^2$ or torus $\T^2$ (equipped with a standard Gaussian Free Field (GFF) with vanishing average for instance). Actually, our techniques can be adapted to other $2$-dimensional Riemannian manifolds with a scalar metric tensor. The main reason is that a Riemannian manifold is locally isometric to $\R^2$.   We will not detail the proofs since the whole machinery works essentially the same as in the plane: it suffices to have at our disposal a white noise decomposition of the underlying Gaussian distribution and to adapt properly our assumptions. 

We rather give here further details in the case of the GFF on planar domains as the associated Brownian motion possesses important conformal invariance properties. We consider a bounded   planar domain $D$. The Liouville Brownian motion on $D$ is defined as follows:\\
$\bullet$ consider a white noise cut-off approximation $(X_\epsilon)_\epsilon$ of the GFF on $D$ with Dirichlet boundary conditions: see equation \eqref{eq:GFFcutoff}.\\
$\bullet$ first define the time change as the limit
$$F'(x,t)=\lim_{\epsilon\to 0}\epsilon^2\int_0^t(2\E[X_\epsilon(x+B_u)^2]-X_\epsilon(x+B_u))e^{2X_\epsilon(x+B_u)}\,du$$ for all $t<\tau^D_x$ 
where $B$ is a standard planar Brownian motion  and $\tau^D_x$ its first exit time out of $D$. This limit turns out to be the same as 
$$\int_0^tC(x+B_u,D)^2(2\E[X(x+B_u)^2]-X(x+B_u))e^{2X(x+B_u)-2\E[X(x+B_u)^2]}\,du$$
where $C(x,D)$ is the conformal radius at $x$ in the domain $D$ (see \cite{complex} or \cite{cf:DuSh} in a different context).\\
$\bullet$ extend this construction to all possible starting points as in section \ref{sec.markov}. Define the exit time of this LBM out of $D$ by $\hat{\tau}^D_x=\sqrt{2/\pi}F'(x,\tau^D_x)$ and then define the Liouville Brownian motion as in Definition \ref{LBMdefnew} for all time $t<\hat{\tau}^D_x$.\\
$\bullet$ Observe that this Liouville Brownian motion is invariant under conformal reparametrization. This means that for all conformal map $\psi:D'\to D$ the process $\psi^{-1}(\LB^x)$ has the law of the Liouville Brownian motion on $D'$ where in the construction of $F'$ we use the standard reparametrization rule of Liouville field theory $X \rightarrow X \circ \psi + Q \ln |\psi'| $ where $Q=\frac{2}{\gamma}+ \frac{\gamma}{2}$ for a subdomain of $\mathbb{C}$ (or $Q=\frac{2}{\gamma}$ for a GFF on the sphere or the torus with vanishing mean). 
 
%%%%%%%%%%%%%%%%%%%%%%%%%%%%%%%%%%%%%%%%%%%%%%%%%%%%
\appendix
%%%%%%%%%%%%%%%%%%%%%%%%%%%%%%%%%%%%%%%%%%%%%%%%%%%%%%%%%%%%%%% 

%%%%%%%%%%%%%%%%%%%%%%%%%%%%%%%%%%%%%%%%%%%%%%%%%%%%%%%%%%%%%%%%%%%%%%%%%%%%%%%%%%%%%%%%%%%%
%\section{Petites questions a resoudre}\label{questions}
%%%%%%%%%%%%%%%%%%%%%%%%%%%%%%%%%%%%%%%%%%%%%%%%%%%%%%%%%%%%%%%%%%%%%%%%%%%%%%%%%%%%

% 

 \end{document}